\input amstex
\input epsf.tex
\documentstyle{amsppt}
\document
\magnification=1100
\NoBlackBoxes

\hsize=15cm
\vsize=20cm
\baselineskip=12pt

\def\Insert #1#2{\epsfysize = #2 \epsfbox{#1}}

\def\Q{\bold{Q}}
\def\C{\bold{C}}
\def\Z{\bold{Z}}

\let\g\gamma

\def\s{\sigma}
\def\a{\alpha}
\def\b{\beta}

\def\O{{\Cal O}}
\def\Q{\Bbb Q}

\def\*{*}
\let\g\gamma

\def\s{\sigma}
\def\a{\alpha}
\def\b{\beta}

\def\={\equiv}
\def\mod#1{~(\hbox{\rm mod}~#1)}

\def\Pcal{{\Cal P}}

\def\C{\Bbb C}

\def\a{\alpha}

\def\({\big(}
\def\){\big)}

\def\({\left(}          \def\){\right)}

\topmatter
\title
MULTIPLE DIRICHLET SERIES AND MOMENTS\\
OF ZETA AND L--FUNCTIONS\endtitle
\author
Adrian Diaconu\\
Dorian Goldfeld\\
Jeffrey Hoffstein
\endauthor

\thanks
Adrian Diaconu would like to thank AIM for its generous support of this research
in the summer of 2001. The second two authors are partially supported by the
National Science Foundation.
\endthanks

\address
Adrian Diaconu, Columbia University Department of Mathematics, New York, NY
10027
\endaddress
\email
cad\@math.columbia.edu
\endemail

\address
Dorian Goldfeld, Columbia University Department of Mathematics, New York, NY
10027
\endaddress
\email
goldfeld\@columbia.edu
\endemail

\address
Jeffrey Hoffstein, Brown University Department of Mathematics, Providence, RI
02912
\endaddress
\email
jhoff\@math.brown.edu
\endemail

\abstract
This paper  develops an analytic theory of Dirichlet
series in several complex variables which possess
sufficiently many functional equations.   In the first
two sections it is shown how straightforward conjectures
about the meromorphic continuation and polar divisors
of certain such series imply, as a consequence, precise
asymptotics (previously conjectured via random matrix theory)
for moments of zeta functions and quadratic $L$-series.
As
an application of the theory, in a third section, we obtain
the current best known error term for mean values of cubes
of central values of Dirichlet $L$-series. The methods utilized
to derive this result are the convexity  principle for functions
of several complex variables  combined with a knowledge of
groups of functional equations for certain multiple Dirichlet
series.
\endabstract

\endtopmatter

\document

\noindent
{\bf \S 1. Introduction}
\vskip 10pt

A Dirichlet series of type
$$\sum_{m_1=1}^\infty \cdots \sum_{m_n=1}^\infty\;
\frac{1}{m_1^{s_1}\cdots m_n^
{s_n}} \;
\int_{0}^\infty \cdots \int_{0}^\infty\; a\left(m_1,
\ldots, m_n, \; t_1, \ldots, t_{\ell}\right)t_{1} ^{-w_{1}}\cdots
t_{\ell} ^{-w_
{\ell}}\,dt_1\cdots dt_{\ell}$$
(where $a\left(m_1,
\ldots, m_n, \; t_1, \ldots, t_{\ell}\right)$ is a complex valued smooth
function) will be called a {\bf multiple\break
Dirichlet series}. It can be viewed as a Dirichlet series in one variable
whose coefficients are again Dirichlet series in several other variables.
One
of the
simplest examples of a multiple Dirichlet series of more than one variable
is
given by
$$\sum_{d }^\infty \frac{ L(s, \chi_d)}{ |d|^w},$$
where the sum ranges over fundamental discriminants of quadratic fields,
$\chi_d$ is the quadratic character associated to
these fields, and $$L(s, \chi_d) = \sum_{n=1}^\infty
  \frac{\chi_d(n)}{n^s}$$
  is the classical Dirichlet L-function.
This type of double Dirichlet series and a method to obtain its analytic
continuation first appeared in a paper of Siegel \cite{S} in 1956. More
generally, one may consider
$$
Z(s_1, s_2, \ldots, s_{m}, w) = \sum_d \frac{L(s_1, \chi_d)\cdot
L(s_2, \chi_d) \cdots L(s_{m}, \chi_d)}{|d|^{w}}.
\tag 1.1$$

Multiple Dirichlet series arise naturally in many contexts and have
been the subject of a
number of papers in the recent past.  See, \cite{B--F--H--2} for an
overview and
references.  The  reason
for their interest is most apparent when they take the form (1.1).
It is easy to see that if, for fixed $s_1, s_2, \ldots, s_m,$ the
analytic continuation of $Z(s_1, s_2,
\ldots, s_{m}, w)$ could be obtained to all $w \in \Bbb C$
then standard Tauberian
arguments could be used to obtain information about the behavior of
$L(s_1, \chi_d)\cdot
L(s_2, \chi_d) \cdots L(s_{m}, \chi_d)$ as $d$ varies.  For example,
mean values could be
obtained if there is a pole at $w =1$.  The situation becomes
even more interesting
when it is noted that quadratic twists of the $L$-series of
automorphic forms on $GL(m)$
can be viewed as special cases of the product $L(s_1, \chi_d)\cdot
L(s_2, \chi_d) \cdots L(s_{m}, \chi_d)$.  The first example of this
type of application
that we are aware of is \cite{G--H} in the case $m=1$. Here mean
value results are obtained
for quadratic Dirichlet $L$-series. Similar results over a function field
are
obtained in \cite{H--R}, and recently, over more general function
fields, in \cite{F--F}.
Examples of the cases $m=2,3$ when the numerator is the
$L$-series associated to a $GL(m)$ cusp form are given in
\cite{B--F--H--2},
\cite{B--F--H--1}.

In all these examples (except for \cite{F--F}), the analytic
continuation of $(1.1)$ was
obtained by treating the variable $w$ separately.  The fact
that the $L$-series or
products of $L$-series in the numerator occurred in the Fourier
coefficients of certain
metaplectic Eisenstein series was exploited, and analytic
continuation in $w$ was
achieved by the application of Rankin-Selberg transforms.

It later became apparent, however, that there were many advantages to
viewing multiple
Dirichlet series as functions of several complex variables.  In
particular, consider (1.1)
but ``improve" it by redefining the $L$-series in such a way that
$\prod_{i=1}^m L(s_i,\chi_d)$ is the
usual product of $L$-series if $d$ is (the square free part of) a fundamental
discriminant, and is
$\prod_{i=1}^m L(s_i,\chi_{d_0})$ times a correction factor if $d$ is a square
multiple of the
square free part $d_0$.  The correction factors are Dirichlet
polynomials with functional
equations and will be discussed further in Section 4.

The improved, or ``perfect" series, $Z^*(s_1, s_2, \ldots, s_{m},
w)$, then possesses
some unexpected properties.  In particular, in addition to the
obvious functional equations
sending $s_i\rightarrow 1-s_i$, $i = 1,\dots,m$, there are some
``hidden" functional
equations that correspond to some surprising structure when the order
of summation in $Z^*$
is altered.

The fact that such a phenomenon can occur was first observed by Bump
and Hoffstein in the case of $m=1$ and a rational function field, and is
mentioned in \cite{H}. It was first observed and applied  in the case
$m=2$ in
\cite{F--H}.  The possibility of using these extra functional
equations as a basis for
obtaining the analytic continuation of double Dirichlet series was
then discussed in
\cite{B--F--H--2}. It was observed there that in the cases where the
numerator is an
$L$-series of an automorphic form on $GL(m)$, if $m=1,2$ or $3$  then
the functional equations of the corresponding  perfect double
Dirichlet series generate a
finite group. It was also noted that by applying these functional
equations to the region
of absolute convergence a collection of overlapping regions was
obtained whose convex hull
was $\Bbb C^2$.  Thus by appealing to a well known theorem in the theory
of functions of
several complex variables, the complete analytic continuation of $Z^*$
could be
obtained.

In later work, \cite{B--F--H--1}, it was observed that a uniqueness
principle operated in the
cases $m=1,2,3$ and the correction factors were determined by, and
could be computed from,
the functional equations of $Z^*$.  Curiously, for $m\ge 4$ the group of
functional equations becomes infinite and simultaneously the
uniqueness principle
fails.  The space of local solutions becomes 1 dimensional in the
case $m=4$, and higher for
$m>4$.  This appears to correspond to an inability to analytically
continue the double
Dirichlet series past a curve of essential singularities.  See
\cite{B--F--H--1,2} for further
details.  The paper of \cite{F--F}, in addition to providing a completely
general analysis
of the case $m=1$ over a function field, contains some further insights into
this curious phenomenon.

We shall call a multiple Dirichlet series (of $n$ complex variables) {\bf
perfect} if it has meromorphic continuation to $\Bbb C^n$ and, in
addition, it
satisfies a group of functional equations. The case $m=3$ is thus of great
interest as the last instance in which the perfect multiple Dirichlet
series
(for the family of quadratic Dirichlet L--functions) are understood
completely.
In
\cite{B--F--H--1} a  description of the "good"
correction factors was obtained for the case of $m=3$ and an
arbitrary automorphic form
$f$ on $GL(3)$. These are the factors corresponding to primes not
dividing 2 or the level of $f$.  This information was then used to
obtain the analytic
continuation of the associated perfect double Dirichlet series.  As a
consequence,
non-vanishing results for quadratic twists of $L(1/2,f,\chi_d)$ were
obtained.  Also,
after taking a residue at $w=1$, a new proof was obtained for the
analytic continuation of
the symmetric square of an automorphic form on $GL(3)$.

One purpose of this paper is to apply the ideas of \cite{B--F--H--1} to
obtain the meromorphic continuation of the series $Z^*(s,s,s,w)$.
After obtaining this and
developing a sieving method analogous to that used in \cite{G--H} we
reconstruct the
unimproved series of (1.1).  Applying the analytic properties of this we
prove
the following
\vskip 10pt
\proclaim{\bf Theorem 1.1}  For $d$ summed over fundamental
discriminants, and any $\epsilon
>0$
$$\sum_{|d| \le x} L\left({\scriptstyle \frac12}, \;\chi_d\right)^3
\left( 1 - \frac{|d|}{x}\right)\;\; = \;\;
\frac{1}{2} \cdot \frac{6}{\pi^2}\, a_3 \cdot \,\frac{1}{2880}\cdot x
\,(\log x)^6 + \sum_{i=0}^5 c_i x (\log x)^i + \O_\epsilon\left (x
^{\frac{4}{5} +
\epsilon}\right).
$$
The constants $c_i$ are effectively computable.
The following unweighted estimate also holds:
$$\sum_{|d| \le x} L\left({\scriptstyle \frac12}, \;\chi_d\right)^3
\;\; = \;\;
\frac{6}{\pi^2}\, a_3 \cdot \,\frac{1}{2880}\cdot x
\,(\log x)^6 + \sum_{i=0}^5 d_i x (\log x)^i + \O_\epsilon
\left(x^{\theta
+ \epsilon}\right),
$$
where the constants $d_i$ are also effectively computable and
$$\theta = \frac{1}{36}\left(47 - \sqrt{265}\right) \sim 0.853366...$$
\endproclaim
\noindent
This  improves on Soundararajan's
\cite{So}, bound of
$O\left(x^{\frac{11}{12}+\epsilon}\right)$.  The weight
$\left(1-\frac{|d|}{x}\right)$ is included in the first part to show the optimal error
term obtainable by this method. It will be shown in \S 4.4, Proposition 4.12, that we
expect the multiple Dirichlet series $Z(\frac12, \frac12,
\frac12, w)$ to have an additional simple pole at $w = \frac34$ with
non--zero residue.
 Accordingly, we conjecture:

\vskip 10pt
\proclaim{\bf Conjecture 1.2}  For $d$ summed over fundamental
discriminants, and any $\epsilon > 0,$
$$\sum_{|d| \le x} L\left({\scriptstyle \frac12}, \;\chi_d\right)^3
\;\; = \;\;
\frac{6}{\pi^2}\, a_3 \cdot \,\frac{1}{2880}\cdot x
\,(\log x)^6 + \sum_{i=0}^5 d_i x (\log x)^i + bx^{\frac34} + \O_\epsilon
\left(x^{\frac{1}{2}
+ \epsilon}\right),
$$
for effectively computable constants $b \ne 0$ and $d_i \, (i = 0, \ldots, 5).$
\endproclaim

\noindent
{\bf Remark:}
In general, for higher moments, we expect additional terms of lower order in the
full moment conjecture besides the terms coming from the multiple pole at $w = 1.$
This is an interesting problem which we hope to return to at a future time.

\vskip 10pt

The major
objective of this paper is to, at least conjecturally, pass the barrier of
$m
\ge 4$.   The first obstacle to accomplishing this is our incomplete
understanding of the
correct form of the class of perfect multiple Dirichlet series for
$m \ge 4$.  There is an
infinite family of choices, every member of which possesses the
correct functional
equations.  However, for any one of these choices, if an analytic
continuation could be
obtained to a neighborhood including the point
$(1/2,1/2,\dots,1/2,1)$ then a sieving
argument could be applied and a formula analogous to Theorem 1.1
could be proved. In particular, this would imply  the truth of
Conjecture
3.1 of Conrey, Farmer, Keating, and Snaith giving the precise asymptotics
for
the moments of
$$\sum_{|d| \le x} L(1/2, \chi_d)^m$$ for $m = 1, 2, 3, \ldots.$ In
\cite{B--F--H--2} it is explained how if the variables are
specialized to $s = s_1
= \dots = s_{m},$ then any multiple Dirichlet series
possessing the correct
functional equations must hit a certain curve of essential
singularities.  A similar hypercurve is encountered for $m \ge 4$ when the
variables are not specialized.  However, the
point
$(1/2,1/2,\dots,1/2,1)$ lies well inside the boundary of this curve.
Another way of saying
this is that by taking the area of absolute convergence of a
corrected analog of (1.1) and
applying the infinite group of functional equations a region of
analytic continuation is
obtained.  For $m \ge 4$ the point $(1/2,1/2,\dots,1/2,1)$  lies
outside this region, but
inside the region contained by the curve of essential singularities.
The case $m=4$ is
particularly intriguing, as $(1/2,1/2,1/2,1/2,1)$ lies right on the
edge of the open hyperplane of analytic continuation that
can be obtained.

In Section 3 we make the reasonable assumption that an analytic
continuation exists past
the point $(1/2,1/2,\dots,1/2,1)$ for a corrected analog of (1.1).
We then calculate the
contribution of the $2^m$ polar divisors  of (1.1) that pass through
this point.  This
gives us a description of the whole principle part in the Laurent
expansion of (1.1) around
this point.  This description is then translated into Conjecture 3.1.

   As far as the present authors are aware, the first examples of multiple
Dirichlet series involving integrals
appear in the paper of A. Good \cite{G} first announced in 1984. Let
$f(z)$ be
a holomorphic cusp form of even weight $k$ for the
modular group $\Gamma = SL(2, \Bbb Z).$ By developing an ingenious
generalization of the Rankin--Selberg convolution in
polar coordinates Good obtained the meromorphic continuation  of the
multiple
Dirichlet series
$$\int_1^\infty \left | L_f\left( \frac{k}{2} + it\right) \right |^2
t^{-w} \,
dt,$$
where $L_f(s)$ is the Hecke L--function associated to $f$ by Mellin
transform. This function has simple poles at $w =
\frac12 + ir$ where $\frac14 + r^2$ is an eigenvalue associated to a
Maass
form on $\Gamma.$ Good \cite{G} even showed how to
introduce weighting factors into the integral which gave a functional
equation
in
$w$. His method can also be extended to  obtain the meromorphic
continuation of
$$\int_1^\infty L_f(s_1 + it) L_f(s_2 - it) t^{-w} dt.$$

In section 2, we develop the theory of multiple Dirichlet series
associated to
moments of the Riemann zeta function. In this case,
the perfect object has been found for $m = 2$ (using theta functions) and
for
$m = 4$ (using Eisenstein series) by Good \cite{G}, but
his theory has never been fully worked out.   We consider the multiple
Dirichlet series
$$Z(s_1, \ldots, s_{2m}, w) = \int_1^\infty \zeta(s_1 + it) \cdots
\zeta(s_m +
it)\cdot \zeta(s_{m+1} - it) \cdots \zeta(s_{2m} -
it)
\; t^{-w} \; dt$$
and show that it has meromorphic continuation (as a function of $2m + 1$
complex variables) slightly beyond the region of
absolute convergence given by $\Re(s_i) > 1, \Re (w) > 1$ ($i = 1, 2,
\ldots,
2m)$ with a polar divisor at
$w = 1.$ We also show that $Z(s_1, \ldots, s_{2m}, w)$ satisfies certain
quasi--functional equations (see
section 2.2) which allows one to meromorphically continue the multiple
Dirichlet series to an
even larger region. It is proved (subject to Conjecture 2.7) that $Z\left(\frac12,
\ldots,
\frac12, w\right)$ has a multiple pole at the point $w = 1,$ and the leading
coefficient in the Laurent expansion is computed explicitly in Proposition 2.9.
Under the assumption that
$Z\left(\frac12, \ldots, \frac12, w\right)$ has holomorphic continuation to the
region
$Re(w)
\ge 1$ (except for the multiple pole at $w = 1,$ we derive the
Conrey--Ghosh--Keating--Snaith conjecture (see \cite{Ke--Sn--1} and
\cite{C--Gh--2}) for the $(2m)^{\text{th}}$ moment of the zeta function  as
predicted by random matrix theory.

Recently \cite{CFKRS} have presented a heuristic method via approximate
functional equations for obtaining moment conjectures for integral
as well as real and complex moments for general families of zeta and
L-functions. Their method is related to ours in that it uses a group
of approximate functional equations in several complex variables.

\vskip 10pt\noindent
{\bf  Acknowledgments:} The authors would like to extend their warmest
thanks to D. Bump, B. Conrey, S. Friedberg,  P. Sarnak for
many very helpful conversations.

\vskip 20pt\noindent
{\bf \S 2. Moments of the Riemann Zeta--Function}
\vskip 10pt
     For $\Re(s) > 1$, let
$$\zeta(s) \; = \; \sum_{n=1}^\infty \frac{1}{n^s} \; = \; \prod_p \left(1
-
\frac{1}{p^s}\right)^{-1}$$ denote the Riemann zeta function which has
meromorphic
continuation to the whole complex plane with a single simple pole at $s =
1$
with residue $1$. It is well known (see Titchmarsh \cite{T}) that $\zeta$
satisfies the
functional equation
$$\zeta(s) = \chi(s) \zeta(1 - s)$$
where
$$\chi(1 - s) = \frac{1}{\chi(s)} = 2 (2\pi)^{-s} \cos\left(\frac{\pi
s}{2}
\right)
\Gamma(s). \tag 2.1$$
In 1918 Hardy and Littlewood \cite{H--L} obtained the second moment
$$\int_0^x \big | \zeta\left({\scriptstyle\frac12 }+ it\right)\big |^2 \,
dt
\sim x
\log x,$$ and in 1926 Ingham \cite{I} obtained the fourth moment
$$\int_0^x \big |\zeta\left({\scriptstyle\frac12 } + it\right)\big |^4 \,
dt
\sim
\frac{1}{2\pi^2} x (\log x)^4.$$
This result has not been significantly improved until the recent work  of
Motohashi \cite{Mot1} in $1993$ where it was shown that
$$\int_0^x \big |\zeta\left({\scriptstyle\frac12 } + it\right)\big |^4 \,
dt =
x\cdot
P_4(\log x) + O\big( x^{\frac23 + \epsilon}\big),$$
where $P_4$ is a certain polynomial of degree four. Motohashi's work was
based
on
an earlier remarkable observation of Deshouiller and Iwaniec \cite{D--I}
that
integrals
of the Riemann zeta function along the critical line occurred in both the
Selberg and
Kuznetsov trace formula, and that the trace formulae could, therefore, be
used
to
obtain new information about the Riemann zeta function.  By a
careful analysis of the Kuznetsov trace formula, Motohashi
\cite{Mot2} introduced  and was able to obtain the meromorphic continuation
(in
$w$) of the function
  $$\int_1^\infty \zeta\left(s + it\right)^2
\zeta\left(s - it\right)^2  t^{-w} \, dt. \tag 2.2$$
Motohashi pointed out that it is, therefore, possible to view the Riemann
zeta
function
as a generator of Maass wave form L--functions.

  There has been a longstanding folklore
conjecture that
$$\int_0^x \big |\zeta\left({\scriptstyle\frac12 } + it\right)\big |^{2k}
\, dt
\sim
c_k x (\log x)^{k^2}. \tag 2.3$$
In 1984 Conrey and Ghosh \cite{C--Gh--2} gave the more precise conjecture that
$$c_k = \frac{g_k a_k}{\Gamma(1 + k^2)} \tag 2.4$$
where $$a_k = \prod_p \left(1 - \frac{1}{p}\right)^{k^2} \sum_{j =
0}^\infty
\frac{d_k\left(p^j\right)^2}{p^j} \tag 2.5$$
is the arithmetic factor and $g_k$, an integer, is a geometric factor.
Here,
$d_k(n)$
denotes the number of representations of $n$ as a product of $k$ positive
integers. In
this notation, the result of Hardy and Littlewood states that
$g_1 = 1$, while Ingham's result is that
$g_2 = 2.$ In 1998, Conrey and Ghosh \cite{C--Gh--1} conjectured that $g_3 =
42$,
and more recently in 1999, Conrey and Gonek \cite{C--G}  conjectured
that $g_4 =
24024.$ Up to this point, using classical techniques based on
approximating
$\zeta(s)$
  by Dirichlet polynomials,  there seemed to be no way to conjecture the
value
of $g_k$
in general.

  In accordance with the philosophy of Katz and Sarnak \cite{K--S} that
one may
associate probability spaces over compact classical groups to families of
zeta
and
L--functions, Keating and Snaith \cite{Ke--Sn--2} (see also \cite{B--H}) computed
moments of characteristic
polynomials of matrices in the unitary group $U(n)$ and formulated the
conjecture
that $$g_k = k^2! \prod_{j = 0}^{k - 1} \frac{j!}{(j + k)!} \tag 2.6$$
for any positive integer $k$. This conjecture agreed with all the known
results
and was
strongly supported by numerical computations.

We show in the next sections that there exists a multiple Dirichlet series
of
several
complex variables of the type
$(2.2)$  previously introduced by Motohashi, with a polar divisor at $w =
1$,
whose
residue is simply related to the constants  $(2.4), (2.5), (2.6)$. We
further
show that if
one could holomorphically continue this multiple Dirichlet series slightly
beyond this
polar divisor, a proof of the Conrey--Ghosh--Keating--Snaith conjecture
would
follow.

\pagebreak

\noindent
{\bf \S 2.1  The Multiple Dirichlet Series for the Riemann Zeta Function}
\vskip 10pt
Let $s_1, s_2, \ldots, s_{2m}, w$ denote complex variables, $k$ be
an integer, and $\epsilon_i = \pm 1$ for $i = 1, 2, \ldots, 2m.$
We shall consider  multiple Dirichlet series of  type
$$Z_{\epsilon_1, \ldots, \epsilon_{2m}, k}(s_1, \ldots, s_{2m}, w) =
\int_1^\infty \zeta(s_1+\epsilon_1it) \cdots
\zeta(s_{2m} +
\epsilon_{2m}it)\left(\frac{2\pi e}{t}\right)^{k i t} \; t^{-w}\, dt.\tag
2.7$$

It is easy to see that the integral in $(2.7)$ converges absolutely for
$\Re(w)
>
1$ and $\Re(s_i) > 1,$ $(i = 1, 2, \ldots, 2m)$,   and defines (in this
region)
a
holomorphic function of
$2m + 1$ complex variables.  These series are more general than the series
$(2.2)
$ introduced by Motohashi in that they contain the factor
$\left(\frac{2\pi
e}{t}\right)^{k i t}$. It will be shortly seen that this factor
occurs naturally
because of the asymptotic formulae
  \cite{T}
$$\align &\chi(s + it) = e^{\frac{i\pi}{4}}
\left(\frac{2\pi}{t}\right)^{s-
\frac12}
\left(\frac{2\pi e}{t}\right)^{it} \Big\{1 +
O\left(\frac{1}{t}\right)\Big\},\\
&\chi(s - it) = e^{-\frac{i\pi}{4}}
\left(\frac{2\pi}{t}\right)^{s-\frac12}
\left(\frac{2\pi e}{t}\right)^{-it} \Big\{1 +
O\left(\frac{1}{t}\right)\Big\}
\qquad
\text{(for fixed $s$ and $t \to \infty$) },\tag 2.8\endalign$$
for $\chi$, the function occurring in the functional equation $(2.1)$ for
the
Riemann
zeta function.

\vskip 10pt
\proclaim{Proposition 2.1}  For $\sigma > 0$, the function $Z_{\epsilon_1,
\ldots,
\epsilon_{2m}, k}(s_1, \ldots, s_{2m}, w)$ can be holomorphically
continued to
the
domain $\Re(s_i) > -\sigma$ (for $i = 1,\ldots, 2m$) and $\Re(w) > 1
+2m\left(\frac12 + \sigma\right).$ Furthermore, for $k \ne 0$,
$Z_{\epsilon_1,
\ldots,
\epsilon_{2m}, k}$ can be holomorphically continued for $\Re(w) > 0$ and
$\Re\left(s_i + \frac{w}{|k|} \right) > 1 + |k| ^{-1}$ ($i = 1,\ldots,
2m$),
and for $k=0$, it can be meromorphically continued for $\Re(w) > 0$ and
$\Re
(s_i) > 1$ ($i = 1,\ldots, 2m$) with a single simple pole at $w = 1$ with
residue
$$\underset {\ell_1^{\epsilon_1}\cdots\ell_{2m}^{\epsilon_{2m}} =
1}\to{\sum_{\ell_1, \ldots, \ell_{2m}}} \ell_1^{-s_1} \cdots
\ell_{2m}^{-s_{2m}}
$$
\endproclaim

{\bf Proof:} The first part of the Proposition follows immediately from
the
well known
convexity bound
$$|\zeta(s + it)| \ll_{s}(1 + |t| )^{\frac12 + \sigma},$$
for $\Re(s)  > -\sigma$, where the implied constant depends at most on
$s$. For
the second part, we need the following lemma.

\vskip 10pt

\proclaim{Lemma 2.2}   Let $B  > 0$ and $k \in \Bbb R$ be fixed. For
$\Re(w) > 1$ the integral
$$|I_{B,k}(w)| = \int_1^\infty B^{it} \left(\frac{2\pi e}{t}\right)^{kit}
t^{-
w} \, dt$$
converges absolutely and defines a holomorphic function of $w$. Further,
for
$\{B,
k\}
\ne \{1, 0\}$, the function
$I_{B,k}(w)$ may be holomorphically continued to $\Re(w) > 0,$ and for $0
< \Re
(w) \le 1,$ it satisfies the
bound
$$|I_{B,k}(w)| \ll_{k, w} \cases \frac{1}{|\log B|} &\;
\text{if $k = 0$},\\
1 + B ^{\frac{1 - \Re(w)}{k}}\left(1 + |\log B| \right)\; & \text{if $k\ne
0.$}
\endcases$$
Finally, when $B = 1, k = 0$, we have $I_{1,0}(w) = \frac{1}{w - 1}.$
\endproclaim

{\bf Proof:} First, a simple computation shows that $I_{1,0}(w) =
\frac{1}{w -
1}.$ Also, integrating by parts, it can easily be seen that $I_{B,0}(w)$
is a
holomorphic function for $\Re(w) > 0.$ In this case, we have the estimate
$$|I_{B,0}(w)|\ll_{w} \frac{1}{|\log B|}.$$
For $k \ne 0$ and $B ^{\frac{1}{k}}\geq (2\pi) ^{-1},$ we split the
integral
defining $I_{B,k}$ into two parts
$$I_{B,k}(w) = \; \int_{1} ^{\frac{A+1}{e}} \left(\frac{A}{t}\right)^{kit}
t^{-
w} \, dt \; + \; \int_{\frac{A+1}{e}}^\infty
\left(\frac{A}{t}\right)^{kit} t^{-
w} \, dt,$$
where $A = 2\pi e\cdot B^\frac{1}{k}.$ We estimate the first integral
trivially, so, for $0 < \Re(w) \le 1,$
$$\align \Bigg|\int_{1} ^{\frac{A+1}{e}} \left(\frac{A}{t}\right)^{kit}
t^{-w}
\, dt \Bigg|< \frac{\left( \frac{A + 1}{e} \right) ^{1 - \Re(w)} - 1}{1 -
\Re
(w)} &< \left(\frac{A+1}{e} \right) ^{1 - \Re(w)}\log \left( \frac{A+1}{e}
\right) \\
& \ll_{k, w}B ^{\frac{1 - \Re(w)}{k}}\left(1 + |\log B| \right).
\endalign$$
Now, integrating by parts, we have
$$\align &\int_{\frac{A+1}{e}}^\infty \left(\frac{A}{t}\right)^{kit}
t^{-w} \,
dt = \int_{\frac{A+1}{e}}^\infty \left(\frac{A}{t}\right)^{kit}
ik (\log A - \log t - 1)\cdot  \frac{1}{ik (\log A - \log t - 1)t^w}\, dt
\\
&= \frac{1}{ik}\cdot \frac{\left( \frac{eA}{A + 1} \right) ^{\frac{ki(A +
1)}
{e} } }{\log \left( 1 + \frac{1}{A} \right)}\cdot \frac{e ^{w}}{(A+1)
^{w}} \;
- \; \frac{1}{ik}\int_{\frac{A+1}{e}}^\infty
\left(\frac{A}{t}\right)^{kit}
\cdot \frac{1}{(\log A - \log t - 1) ^{2}t^{w + 1}}\, dt \\
&\hskip145pt +\; \frac{w}{ik}\int_{\frac{A+1}{e}}^\infty \left(\frac{A}{t}
\right)^{kit}\cdot \frac{1}{(\log A - \log t - 1)t^{w + 1}}\,
dt.\endalign$$
It follows that the last two integrals converge absolutely for $\Re(w) >
0,$
and hence, the function $I_{B,k}$ is holomorphic in this region. Moreover,
we
have the estimate
$$\align \Bigg|\int_{\frac{A+1}{e}}^\infty \left(\frac{A}{t}\right)^{kit}
t^{-
w} \, dt \Bigg| &\ll \frac{e ^{\Re(w)}A ^{1 - \Re(w)}}{|k|} \; \\
&+ \; \Big|\frac{w}{k} \Big|\int_{\frac{A+1}{e}}^\infty \frac{1}{\log
\left
(\frac{te}{A} \right)}\cdot \frac{1}{t^{1 + \Re(w)}}\, dt \; + \;
\frac{1}{|k|}
\int_{\frac{A+1}{e}}^\infty \frac{1}{\log ^2 \left( \frac{te}{A}
\right)}\cdot
\frac{1}{t^{1 + \Re(w)}}\, dt \\
&\ll \frac{|w|}{\Re(w)}\cdot \frac{e ^{\Re(w)}A ^{1 - \Re(w)}}{|k|}\;
+ \; \frac
{e ^{\Re(w)} }{|k|A ^{\Re(w)}}\int_{1 + \frac{1}{A}}^\infty \frac{1}{\log
^2 u}
\cdot \frac{1}{u^{1 + \Re(w)}}\, du \\
&\ll \frac{|w|}{\Re(w)}\cdot \frac{e ^{\Re(w)}A ^{1 - \Re(w)}}{|k|}\ll_{k,
w} B
^{\frac{1 - \Re(w)}{k}},
\endalign$$
which combined with the previous one gives the required bound for the
function
$I_{B,k}.$
For the remaining case, $B ^{\frac{1}{k}} < (2\pi) ^{-1},$ we split once
again
the integral into two parts
$$I_{B,k}(w)= \int_{1}^{1+\frac{1}{e}}\left(\frac{A}{t}\right)^{kit}
t^{-w}\,
dt\, +\, \int_{1+\frac{1}{e}} ^\infty \left(\frac{A}{t}\right)^{kit}
t^{-w}\,
dt.$$
A similar argument implies that the second integral converges absolutely
for
$\Re(w) > 0,$ and that
$$|I_{B,k}(w)|\ll_{k, w} 1.$$

We now return to the proof of Proposition 2.1.
For $\Re(s_i)  > 1$   ($i =
1,\ldots, 2m$),
$$Z_{\epsilon_1, \ldots,
\epsilon_{2m}, k}(s_1, \ldots, s_{2m}, w) = \sum_{\ell_1, \ldots,
\ell_{2m}
}\ell_1^{-s_1}\cdots \ell_{2m}^{-s_{2m}} \int_1^\infty
\big(\ell_1^{\epsilon_1}
\cdots
\ell_{2m}^{\epsilon_{2m}}\big)^{it} \left(\frac{2\pi e}{t}\right)^{k i t}
\; t^
{-w}\,
dt, \tag 2.9$$
where the sum ranges over  all $2m$--tuples $\{\ell_1, \ldots,
\ell_{2m}\}$ of
positive
integers. For $k \ne 0$ and $0 < \Re(w) \le 1,$ it is clear that the
series on
the right side of $(2.9)$ is absolutely convergent provided $\Re(s_i)$ ($i
=
1,\ldots, 2m$) are sufficiently large. In fact, the estimates from Lemma
2.2
imply that we have absolute convergence even for $\Re \left(s_i +
\frac{w}{|k|}
\right) > 1 + |k| ^{-1}$ ($i = 1,\ldots, 2m$). For $k=0,$ we break the sum
on
the right side of $(2.9)$ into two parts
$$\sum_{\ell_1, \ldots, \ell_{2m}} \;\; = \; \underset
{\ell_1^{\epsilon_1}\cdots\ell_{2m}^{\epsilon_{2m}} = 1}\to{\sum_{\ell_1,
\ldots,
\ell_{2m}}} \; + \; \underset
{\ell_1^{\epsilon_1}\cdots\ell_{2m}^{\epsilon_
{2m}} \ne
1}\to{\sum_{\ell_1, \ldots, \ell_{2m}}}.\tag 2.10$$
By Lemma 2.2 it immediately follows that the first sum in $(2.10)$
will contribute a pole at $w = 1$ with residue precisely as stated in
Proposition 2.1. It is also clear from Lemma 2.2 that the second sum
in $(2.10)$
will give a holomorphic contribution to $(2.9)$ provided $\Re(s_i)$ ($i =
1,
\ldots,
2m$) are sufficiently large so that the sum over
$\ell_1, \ldots ,\ell_{2m}$ converges absolutely. To show convergence for
$\Re(s_i)   > 1$ ($i = 1, \ldots,
2m$) is more delicate and we give the details.

   It follows from Lemma 2.2 that for $\Re(s_i) = \sigma  > 1$, ($i = 1,
\ldots,
2m$),
$$\align\underset {\ell_1^{\epsilon_1}\cdots\ell_{2m}^{\epsilon_{2m}} \ne
1
}\to{\sum_{\ell_1, \ldots, \ell_{2m}}}   \ell_1^{-s_1}\cdots
\ell_{2m}^{-s_{2m}}
&\int_1^\infty \big(\ell_1^{\epsilon_1}
\cdots
\ell_{2m}^{\epsilon_{2m}}\big)^{it} \; t^{-w}\,
dt \;\; \\  &\ll_{w} \; \underset
{\ell_1^{\epsilon_1}\cdots\ell_{2m}^{\epsilon_
{2m}}
\ne 1}\to{\sum_{\ell_1, \ldots, \ell_{2m}}}
\frac{1}{(\ell_1\cdots\ell_{2m})^\sigma}
\frac{1}{\big |\log \ell_1^{\epsilon_1}\cdots\ell_{2m}^{\epsilon_{2m}}\big
|}.
\tag 2.11\endalign$$
We now break the sum on the right side of $(2.11)$ into two parts
$$\underset {\ell_1^{\epsilon_1}\cdots\ell_{2m}^{\epsilon_{2m}}
\ne 1}\to{\sum_{\ell_1, \ldots, \ell_{2m}}} \; = \;
\underset {\ell_1^{\epsilon_1}\cdots\ell_{2m}^{\epsilon_{2m}}
\; \in \; (0,\frac12] \, \cup \, [2,\infty)}\to{\sum_{\ell_1, \ldots,
\ell_
{2m}}} \; + \;
\underset {\ell_1^{\epsilon_1}\cdots\ell_{2m}^{\epsilon_{2m}}
\; \in \; (\frac12,1) \, \cup \, (1,2)}\to{\sum_{\ell_1, \ldots, \ell_
{2m}}}.\tag 2.12$$
The first series on the right side of $(2.12)$ is obviously convergent for
$\sigma > 1.$
We shall show that the second one is also convergent.

  Without loss of generality, let us write
$$\ell_1^{\epsilon_1}\cdots\ell_{2m}^{\epsilon_{2m}} =
\frac{\ell_1\cdots\ell_r}{\ell_{r+1}\cdots\ell_{2m}}.$$
It follows, upon setting $\ell_1\cdots\ell_r = k, \;
\ell_{r+1}\cdots\ell_{2m}
= k \pm a$,
that
$$\align &\underset {\frac{\ell_1\cdots\ell_r}{\ell_{r+1}\cdots\ell_{2m}}
\; \in \; (\frac12,1) \, \cup \, (1,2)}\to{\sum_{\ell_1, \ldots,
\ell_{2m}}}
\hskip -20pt \frac{1}{(\ell_1\cdots\ell_{2m})^\sigma}
\frac{1}{\big |\log
\frac{\ell_1\cdots\ell_r}{\ell_{r+1}\cdots\ell_{2m}}\big |}
\; = \; \sum_{k=2}^\infty \frac{d_r(k)}{k^\sigma} \;\sum_{a=1}^{k-1}
\frac{d_{2m-r}(k+a)}{(k+a)^\sigma}\cdot \frac{1}{\log\left(1 + \frac{a}{k}
\right)} \\
& - \; \sum_{k=3}^\infty \frac{d_r(k)}{k^\sigma}
\;\sum_{a=1}^{\left[\frac{k}{2}
\right]}
\frac{d_{2m-r}(k-a)}{(k-a)^\sigma}\cdot \frac{1}{\log\left(1 - \frac{a}{k}
\right)} \ll \; \sum_{k=2} ^\infty \frac{d_r(k)}{k^\sigma} \;
\sum_{a=1}^{k-1}
\frac{d_{2m-r}(k+a)}{(k+a)^\sigma}\cdot \frac{k}{a} \\
& + \; \sum_{k=3} ^\infty \frac{d_r(k)}{k^\sigma} \;
\sum_{a=1}^{\left[\frac{k}
{2}\right]} \frac{d_{2m-r}(k - a)}{(k - a)^\sigma}\cdot \frac{k}{a}
\ll_{m, r, \epsilon} \; \sum_{k=2}^\infty
\frac{1}{k^{\sigma-\epsilon}} \; \sum_
{a=1}^{k-1}
\frac{k}{a(k+a)}\\
& + \;  \sum_{k=3}^\infty \frac{1}{k^{\sigma-\epsilon}} \; \sum_{a=1}
^{\left
[\frac{k}{2}\right]}\frac{k}{a(k-a)}\ll \; \sum_{k=2}^\infty \frac{\log
k}{k^
{\sigma-\epsilon}},
\endalign$$
for some arbitrarily small $\epsilon > 0.$ Clearly, the last sum converges
if
$\sigma >
1.$ This completes the proof of Proposition 2.1.

    We now deduce a more precise form of the residue given in Proposition
2.1.
This is
given in the next proposition.

\vskip 10pt
\proclaim{Proposition 2.3} Fix $\epsilon > 0.$ Let $\Re(s_i) >
2+\epsilon,$
$\epsilon_i =
\pm 1, (i = 1, \ldots, 2m),$ and define $r$ to be the number of
$\epsilon_i =
1,  (i = 1,
\ldots, 2m).$ If $Z_{\epsilon_1, \ldots, \epsilon_{2m},
k}$ denotes the multiple Dirichlet series defined in $(2.7)$, then we have
$$\underset {w = 1}\to{\text{\rm Res}}\Big[ Z_{\epsilon_1, \ldots,
\epsilon_
{2m},
0}(s_1, \ldots ,s_{2m}, w)\Big ] = R_r(s_1, \ldots, s_{2m}) \cdot
\underset r+1
\le j \le 2m\to{\prod_{1 \le  i \le r}}
\zeta(s_i + s_j),$$
where  $R_r(s_1, \ldots, s_{2m})$ can be holomorphically continued to
the region
$\Re(s_i) > \frac12 - \epsilon.$  Further,
$$R_r\left(\frac12, \ldots,\frac12\right) = \prod_p\left(1 -
\frac{1}{p}\right)^
{m^2}
\left(\sum_{\mu = 0}^\infty d_r\left(p^\mu\right)
d_{2m-r}\left(p^\mu\right)
p^{-\mu}\right),$$
and in particular,
$$R_m\left(\frac12, \ldots,\frac12\right) = a_m,$$
the constant defined in $(2.5).$
\endproclaim

{\bf Proof:} Define
$$U_r(s_1, \ldots, s_{2m}) = \underset {\ell_1\cdots \ell_r \; = \;
\ell_{r+1}\cdots\ell_{2m}}\to{\sum_{\ell_1,
\ldots, \ell_{2m}}} \ell_1^{-s_1}\cdots \ell_{2m}^{-s_{2m}}.$$
It follows from Proposition 2.1, that up to a permutation of the variables
$s_1, \ldots, s_m,$ the function $U_r$ is precisely the residue of
$Z_{\epsilon_1, \ldots, \epsilon_{2m},
0} (s_1, \ldots ,s_{2m}, w)$ at $w = 1.$

If $f(n)$ is a multiplicative function for which the sum
$\overset \infty\to{\underset {n=1}\to{\sum}} f(n)$ converges absolutely,
then
  we have the Euler product identity
$$\overset \infty\to{\underset {n=1}\to{\sum}} f(n) = \prod_p (1 + f(p) +
f
(p^2) +
f(p^3) +
\cdots ).\tag 2.13$$
It  follows from $(2.13)$ that
$$U_r(s_1, \ldots, s_{2m}) = \prod_p\left(\sum_{\mu=0}^\infty \;
\; \underset e_i \ge 0, \; (i = 1, \ldots, 2m) \to{\underset \; e_{r+1} +
\cdots + e_{2m}
=
\mu\to{\sum_{e_1+\cdots
+e_r \; = \; \mu}}} p^{-(e_1s_1 + \cdots + e_{2m}s_{2m})}\right ).$$
Let us now define
$$R_r(s_1, \ldots, s_{2m}) =  U_r(s_1, \ldots, s_{2m}) \cdot \underset r+1
\le j \le 2m\to{\prod_{1 \le  i \le r}}
\zeta(s_i + s_j)^{-1}. \tag 2.14$$
By carefully examining the Euler product for the right hand side of
$(2.14)$,
one
sees that $R_r(s_1, \ldots, s_{2m})$ is holomorphic for $\Re(s_i) >
\frac12 -
\epsilon,
\; (i = 1, \ldots, 2m)$.

Now,
$$\underset e_i \ge 0, \; (i = 1, \ldots, 2m) \to{\underset \; e_{r+1} +
\cdots
+ e_{2m}
=
\mu\to{\sum_{e_1+\cdots
+e_r \; = \; \mu}}} 1 \; = \; d_r\left(p^\mu\right) d_{2m-r}\left
(p^\mu\right).$$
Consequently, if we specialize the variables to $s_1 = s_2 = \cdots =
s_{2m} =
s$, we
obtain
$$R_r(s, \ldots, s) = \prod_p\left(1 - \frac{1}{p^{2s}}\right)^{m^2}
\left(\sum_{\mu = 0}^\infty d_r\left(p^\mu\right)
d_{2m-r}\left(p^\mu\right)
p^{-2\mu s}\right).$$
The proof of Proposition 2.3 immediately follows upon letting $s \to
\frac12.$

\vskip 20pt\noindent
{\bf 2.2 Quasi--Functional Equations}
\vskip 10pt
Fix variables $s_1, s_2, \ldots, s_{2m}, \, w$.
Let $\Cal D_{s_1, \ldots, s_{2m}, w}$ denote the infinite dimensional vector space,
defined over the field
$$K_{s_1, \ldots, s_{2m}} = \Bbb
C\Big((2\pi)^{s_1}, \; \ldots,  \;  (2\pi)^{s_{2m}}\Big ),
$$
 generated by the multiple Dirichlet series
$$Z_{\epsilon_1, \ldots, \epsilon_{2m}, k} (S_1, \ldots, S_{2m}, W),$$
where the variables $\epsilon_j, \; k, \; S_j, $ and $W$ range over the values:
  $$\epsilon_j \in \{\pm 1\}, \; (j = 1, \ldots, 2m)$$
 $$\; k \in \Bbb Z,$$
$$S_j \in  \{s_j, 1-s_j\}, \;\; (j = 1, \ldots, 2m),$$
$$W = w +
\sum_{j=1}^{2m}
\delta_j\left(s_j-\frac12\right)$$
 with $\delta_j \in \{0, 1\}, \;\; (j = 1, \ldots, 2m).$

For
$j = 1, 2,
\ldots ,2m,$ we will define
involutions
   $\gamma_j : \Cal D_{s_1, \ldots, s_{2m}, w} \to \Cal D_{s_1, \ldots, s_{2m}, w}.$
\vskip 10pt
\proclaim{Definition 2.4}  For $j = 1, 2, \ldots ,2m,$ we define an action
$\gamma_j$
on
  $$Z_{\epsilon_1, {\scriptstyle \ldots}, \epsilon_{2m},k}\big(S_1, \ldots,
S_
{2m},
W\big) \in \Cal D_{s_1, \ldots, s_{2m}, w}$$ (the
action denoted by a right superscript) as follows:
$$Z_{\epsilon_1, {\scriptstyle \ldots}, \epsilon_{2m},k}\big(S_1, \ldots,
S_
{2m},
W\big)^{\gamma_j} =  e^{\frac{i\pi \epsilon_j}{4}} (2\pi)^{S_j -\frac12}
\; \,
Z_{\epsilon_1,\ldots , -\epsilon_j,\ldots ,
\epsilon_{2m}, \, k+\epsilon_j}\big(S_1, \,{\scriptstyle \ldots} \, ,
1-S_j, \,
{\scriptstyle
\ldots}\, ,  S_{2m}, W+S_j-{\scriptstyle\frac12}\big).$$
The involutions $\gamma_j, (j = 1, \ldots, 2m)$ generate a finite abelian
group
$G_{2m}$ of
$2^{2m}$ elements which, likewise, acts on $\Cal D_{s_1, \ldots, s_{2m}, w}.$
\endproclaim
We will also denote by $\gamma_j$ $(j = 1, 2, \ldots, 2m)$, the affine
transformations induced by this action
$$(s_1, \ldots, s_{2m}, w)\;\; \overset\gamma_{j}\to{\longrightarrow}\;\;
\big
(s_1, \ldots, 1 - s_j, \ldots, s_{2m}, s_j + w - 1/2 \big).$$

By
Proposition 2.1, we know that $Z_{\epsilon_1,\ldots, \epsilon_{2m},k}
\big(s_1,
\ldots, s_{2m},
w\big)$ has
holomorphic continuation to the region  $$0 < \Re(s_i) < 1,\; (i = 1,
\ldots,
2m), \qquad
\Re(w) > 1 + m. \tag 2.15$$
We would like to use the functional equation $(2.1)$ to obtain a
functional
equation
for the multiple Dirichlet series $Z_{\epsilon_1,\ldots, \epsilon_{2m},k}
\big
(s_1,
\ldots, s_{2m}, w\big)$. To abbreviate notation, we let $$Z(s_1, \ldots,
s_
{2m}, w) =
Z_{\epsilon_1,\ldots, \epsilon_{2m}, k} \big(s_1, \ldots, s_{2m},
w\big).$$

   We shall need an asymptotic expansion of Stirling type \cite{T}
$$\align &\chi(s + it) = e^{\frac{i\pi}{4}}
\left(\frac{2\pi}{t}\right)^{s-
\frac12}
\left(\frac{2\pi e}{t}\right)^{it} \Bigg\{1 + \sum_{n=1}^{N} c_n
t^{-n} + O\left
(t ^{-N-1} \right)\Bigg\},\\
&\chi(s - it) = e^{-\frac{i\pi}{4}}
\left(\frac{2\pi}{t}\right)^{s-\frac12}
\left(\frac{2\pi e}{t}\right)^{-it} \Bigg\{1 + \sum_{n=1} ^{N} \bar c_n
t^{-n}
+ O\left(t ^{-N-1} \right)\Bigg\}
\qquad
\text{(for fixed $s$ and $t \to \infty$) },\tag 2.16\endalign$$
where $c_n$ are certain complex constants. Such expansions are not
explicitly
worked
out in \cite{T}, but they are not hard to obtain.

   It now follows from Definition 2.4, Stirling's asymptotic expansion
$(2.16)$,
and the functional equation $(2.1)$, that in the region $(2.15)$, we have
for
$\gamma \in G_{2m}$, the quasi--functional equation
$$Z\big(s_1, \ldots, s_{2m},
w\big)\;\sim\; Z\big(s_1, \ldots, s_{2m},
w\big)^{\gamma} \; + \; \sum_{n=1}^\infty c'_{n}(\gamma)\; Z\big(s_1,
\ldots, s_
{2m},
w + n\big)^{\gamma},\tag 2.17$$
where $c'_{n}(\gamma_{j}) = c_n$ if $\epsilon_j = +1$ and
$c'_{n}(\gamma_{j}) =
\bar c_n$ if $\epsilon_j = -1$, for $j = 1, 2,\ldots, 2m$, and in general,
$c'_
{n}(\gamma)$ is a linear combination of $c_{n'}$ and $\bar c_{n''}$ with
$n',\;n''\le n.$

  We shall be mainly interested in  $\gamma \in G_{2m}$ for which the
action
given in Definition 2.4
$$Z_{\epsilon_1,\ldots, \epsilon_{2m}, 0} \big(s_1, \ldots, s_{2m},
w\big)  \longrightarrow Z_{\epsilon_1,\ldots, \epsilon_{2m}, 0} \big(s_1,
\ldots,
s_{2m}, w\big)^\gamma \tag 2.18$$
stabilizes $k = 0.$ An element $\gamma \in G_{2m}$ is said to stabilize $k$
relative to $\{\epsilon_1, \ldots, \epsilon_{2m}\}$ provided
$$Z_{\epsilon_1,\ldots, \epsilon_{2m}, k}(s_1, \ldots, s_{2m}, w)^\gamma =
C(s_1,\ldots,s_{2m})\cdot Z_ {\epsilon'_1,\ldots, \epsilon'_{2m}, k'}(s_1', \ldots,
s_{2m}', w')$$
for some $ C(s_1,\ldots,s_{2m}) \in K_{s_1,\ldots, s_{2m}}$ with $k =
k'.$

\vskip 10pt
\proclaim{Definition 2.5} Fix $\epsilon_i = \pm 1, (i = 1, \ldots, 2m).$
We
define
$G_{2m}(\epsilon_1, \ldots, \epsilon_{2m})$ to be the subset of
$G_{2m}$ (defined in Definition 2.4) consisting of all $\gamma \in
G_{2m}$ which
stabilize
$0$ relative to $\{\epsilon_1, \ldots, \epsilon_{2m}\}$.\endproclaim

\vskip 10pt
\proclaim{Proposition 2.6}  Let $1 \le r \le 2m,$ and
$$\epsilon_{i_1} = \epsilon_{i_2} = \cdots = \epsilon_{ir} = + 1,
\qquad \epsilon_{i_{r+1}} =  \epsilon_{i_{r+2}} = \cdots =
\epsilon_{i_{2m}} = -1.$$
Then $G_{2m}(\epsilon_1, \ldots, \epsilon_{2m})$ is  the subgroup of
$G_{2m}$ which is generated by the elements
$\gamma_{i_\mu}\cdot\gamma_{i_\nu}$ with $1 \le \mu \le r, \; r+1
\le \nu \le 2m.$\endproclaim

{\bf Proof:} Note that if we write $\gamma = \gamma_i \cdot
\gamma_j$  (with $i \ne j$) then under the action $(2.18)$ we see that
$$\{k = 0\}\;\;  \overset \gamma \to\longrightarrow \;\; \{k =\epsilon_i +
\epsilon_j\}.$$ So if we choose $i$ from the set $\{i_1, \ldots, i_r\}$
and $j$
from the set $\{i_{r+1}, \ldots, i_{2m}\}$ then we see that $\{k = 0
\}$ is stabilized. It easily follows that these elements generate a
group and every element of this group stabilizes $0$ relative to
$\{\epsilon_1, \ldots, \epsilon_{2m}\}$. Furthermore, every element
which stabilizes $0$ relative to
$\{\epsilon_1, \ldots, \epsilon_{2m}\}$ must lie in this group.
\vskip 10pt\noindent
{\bf Remark:}  We introduced the group $G_{2m}(\epsilon_1, \ldots,
\epsilon_
{2m})$
because it is precisely this group which gives the reflections of the
polar
divisor at $w
= 1$ of the multiple Dirichlet series $Z_{\epsilon_1,\ldots,
\epsilon_{2m}, 0}
\big(s_1, \ldots, s_{2m},
w\big).$ This will be further explained in the next section.

\vskip 20pt\noindent
\S 2.3 {\bf A Fundamental Conjecture for the Riemann Zeta Function}
\vskip 10pt

We observed in Proposition 2.1 that the hyperplane $w-1=0$ belongs to the
polar
divisor of the multiple Dirichlet series $Z_{\epsilon_1, \ldots,
\epsilon_{2m}, k}$ if and only if $k = 0.$ It was also seen that this
hyperplane is the only possible pole in the region $\Cal{F}$ defined by
$$\align \Cal{F} &= \big\{(s_{1}, \ldots, s_{2m}, w) \in \Bbb
C^{2m+1}\;|\; \Re
(s_i)  > 0 \;  (i =
1,\ldots, 2m), \Re(w) > 1 + m \big\}\\
&\cup \big\{(s_1, \ldots, s_{2m}, w) \in \Bbb C^{2m+1}\;|\; \Re(w) > 0,
\Re
(s_i) > 2 \; (i =
1,\ldots, 2m) \big\}. \endalign$$
Now, the set $\bigcap_{\gamma \in G_{2m}}\gamma \big(\Cal{F}\big)$ is
nonempty,
since it contains points for which $\Re(s_i)\sim 1/2$ $(i = 1,\ldots, 2m)$
and
$\Re(w)$ is sufficiently large. It follows from the quasi--functional
equation
$(2.17)$ that the multiple Dirichlet series $Z_{\epsilon_1, \ldots,
\epsilon_
{2m}, 0}$ have meromorphic continuation to the convex closure of the
region
$$\bigcup_{\gamma \in G_{2m}}\gamma \big(\Cal{F}\big)$$
with poles, precisely, at the reflections of the hyperplane $w - 1 = 0$
under
$G_{2m}(\epsilon_1, \ldots, \epsilon_{2m}).$ In order to obtain the
continuation, it is understood that we first multiply $Z_{\epsilon_1,
\ldots,
\epsilon_{2m}, 0}$ by certain linear factors in order to cancel its poles.
We
propose the following conjecture.

\vskip 10pt
\proclaim{Conjecture 2.7} The functions $Z_{\epsilon_1, \ldots,
\epsilon_{2m}, 0}$ have meromorphic continuation to a tube domain in $\Bbb
C^
{2m+1}$ which contains the point $\left(\frac12, \ldots, \frac12, \,
1\right).$
All these functions have the same polar divisor passing through
this point consisting of all the reflections of the hyperplane $w - 1 = 0$
under the group $G_{2m}(\epsilon_1, \ldots, \epsilon_{2m}).$ Moreover, the
functions $$Z_{\epsilon_1, \ldots, \epsilon_{2m},
0}\left(\frac{1}{2},\cdots,
\frac{1}{2},w \right)$$
are holomorphic for $\Re(w) > 1.$
\endproclaim
\vskip 10pt
\proclaim{Theorem 2.8} Conjecture 2.7 implies the
Keating--Snaith--Conrey--
Farmer conjecture $(2.3)$.\endproclaim
\vskip 10pt
{\bf Proof:}
  From now on, we fix
$$\epsilon_{1} = \epsilon_{2} = \cdots = \epsilon_{m} = + 1,
\qquad \epsilon_{m + 1} =  \epsilon_{m + 2} = \cdots =
\epsilon_{2m} = -1,$$
and let $G'_{2m}$ denote the group $G_{2m}(\epsilon_1, \ldots,
\epsilon_{2m}).$
The reflections of the hyperplane $w - 1 = 0$ under the group $G'_{2m}$
are
given by
$$\delta_1 s_1 + \cdots + \delta_{2m} s_{2m} + w - \frac{\delta_1 + \cdots
+ \delta_{2m} + 2}{2}=0, \tag 2.19$$
where $\delta_{i}=0$ or $1$ and
$\delta_1 + \cdots + \delta_m =
\delta_{m+1} + \cdots + \delta_{2m}.$

In this and the next section we require a version of the
Wiener--Ikehara Tauberian theorem.  Stark has proved a vast
generalization of this theorem, \cite{St}.  We will quote here a
limited a case of his result which is sufficient for our needs.
\vskip 10pt
\proclaim{Tauberian theorem (Stark)} Let $S(x)$ be a non-decreasing
function of $x$ and let
$$Z(w) = \int_{1}^\infty S(t) \cdot t^{-w} \, \frac{dt}{t}.$$
Let $P(w)= \gamma_M + \gamma_{M-1}(w-1) + \dots + \gamma_0 (w-1)^M$,
($M \ge 0$) be a polynomial with $\gamma_M \ne 0$ such that $Z(w) -
P(w)(w-1)^{-M-1}$ is holomorphic for  $\Re(w) > 1$ and continuous for
$\Re(w) =  1$.
Then
$$S(x) \;\; \sim \;\; \frac{\gamma_M}{M!} \cdot
\; x
(\log x)^M,
\qquad \text{(as $x \to \infty ).$}$$
\endproclaim
\vskip 10pt

We now let $z(t) = \zeta(1/2 + it) ^{m}$ and
$S(x) = \int_{0} ^x |z(t)| ^{2} \, dt$ in the Tauberian theorem.  It follows
by integration by parts that
$$\int_1^\infty S(t)\cdot t^{-w}\, \frac{dt}{t} = \frac{1}{w} \int_1^\infty |z(t)|^2
t^{-w} \, dt.$$
Consequently,
it is
enough
to show that
$$\lim_{w\to 1} (w - 1)^{m ^{2} + 1} \; Z_{\epsilon_1, \ldots,
\epsilon_{2m}, 0}\big(s_1, \ldots, s_{2m},
w\big)\;=\; g_{2m} \, a_{2m}\, m ^{2}!,$$
where $$g_{2m} = \prod_{\ell =
0}^{m-1}
\frac{\ell !}{(\ell+m)!},$$
and $a_{2m}$ is the constant given in $(2.5)$.

   Let $U(s_1, \ldots, s_{2m}, w)$ denote the function defined by
$$\frac{1}{w - 1} \; R_{m}(s_1, \ldots, s_{2m})\; \prod_{i=1}^m\;
\prod_{j=m+1}^
{2m} \zeta(s_i + s_j). \tag 2.20$$
Then Conjecture 2.7 implies that
$$Z_{\epsilon_1, \ldots,
\epsilon_{2m}, 0}\big(s_1, \ldots, s_{2m},
w\big) - \sum_{\gamma\in G'_{2m}}
U\big(\gamma(s_1, \ldots, s_{2m}, w)\big)\tag 2.21$$
is holomorphic around $\left(\frac12, \ldots, \frac12, \, 1\right).$
The proof of theorem 2.8 is an immediate consequence of the following
proposition.

\proclaim{Proposition 2.9} For $m = 1, 2, \ldots,$ let $G'_{2m}$ denote
the
subgroup of $G_{2m}$ generated by the involutions $\gamma_{ij} =
\gamma_i\cdot\gamma_j,$ ($i=1, \ldots, m$ and $j=m+1, \ldots, 2m$).
Then we have
$$\lim_{w\to 1} \;\;\lim_{(s_1, \ldots, s_{2m}) \to (\frac12,\ldots,
\frac12)}
\Bigg[ (w - 1)^{m ^2+1}\sum_{\gamma\in G'_{2m}}  U\big(\gamma(s_1,
\ldots, s_{2m}, w)\big) \Bigg]
\;\; = \;\; a_{2m}\; g_{2m}\; m ^2! $$
where $$g_{2m} = \prod_{\ell =
0}^{m-1}
\frac{\ell !}{(\ell+m)!},$$
and $a_{2m}$ is the constant given in $(2.5)$.\endproclaim
\vskip 10pt
{\bf Proof:} We start by taking the Taylor expansion of
$$U(s_1, \ldots, s_{2m}, w) \; = \;  a_{2m}\;\cdot \frac{f^*(s_1, \ldots,
s_
{2m})}{(w-1)\overset
m\to{\underset {i = 1}\to{\prod } }\; \overset
2m\to{\underset {j = m+1}\to{\prod } }(s_i+s_j-1)}\tag 2.22 $$
around $(s_1, \ldots, s_{2m}) = \left(\frac12, \ldots,
\frac12\right).$ Here
$$f^*(s_1, \ldots, s_{2m}) = 1 \; + \; \underset \nu_1 + \cdots +
\nu_{2m}\ge 1
\to{\sum_{\nu_1=0}^\infty \cdots
\sum_{\nu_{2m}=0}^\infty} \kappa(\nu_1, \ldots, \nu_{2m})\; \left( s_1 -
\frac{1}{2} \right)  ^{\nu_1}
\cdots \left( s_{2m} - \frac{1}{2} \right)  ^{\nu_{2m}} ,$$
(with $\kappa(\nu_1, \ldots, \nu_{2m}) \in \Bbb C$),
  will be a holomorphic function
which is symmetric separately with respect to the variables $s_1, \ldots,
s_m$
and
$s_{m+1}, \ldots, s_{2m}.$

Now,  make the change of variables $s_i = \frac12 + u_i$ for $i
= 1, 2, \ldots, 2m,$ and $w = v + 1.$ Then, for $i = 1, \ldots, m$ and $j
=
m+1, \ldots, 2m,$ the involutions $\gamma_{ij}$ are transformed to
  $$ (u_1,  \ldots, u_i, \ldots, u_m, \ldots, u_j, \ldots, u_{2m}, v)\;\;
\overset\gamma_{ij}\to{\longrightarrow}\;\; (u_1,
\ldots, -u_j, \ldots, u_m, \ldots, -u_i, \ldots, u_{2m},
u_i + u_j +v).$$
Henceforth, we denote by $G'_{2m}$ the group
generated by the above involutions.

Then by $(2.22)$, it is enough to prove that
$$\lim_{v\to 0} \;\;\lim_{(u_1, \ldots, u_{2m}) \to (0,\ldots, 0)}
\Bigg[ v^{m ^2+1}\sum_{\gamma\in G'_{2m}}  H_f\big(\gamma(u_1,
\ldots, u_{2m}, v)\big) \Bigg]
\;\; = \;\; g_{2m}\; m ^2!, \tag 2.23$$
where
$$H_f(u_1, \ldots, u_{2m}, v) \;\; = \;\; \frac{1}{v}\;\cdot
\frac{f(u_1,
\ldots, u_{2m})}{\overset m\to{\underset {i = 1}\to{\prod} }
\; \overset 2m\to{\underset {j = m+1}\to{\prod} }(u_i + u_j)},$$
and $f$ (which is simply related to $f^*$) is a certain holomorphic
function
and symmetric
separately with respect to the variables $u_1, \ldots, u_m$ and
$u_{m+1}, \ldots, u_{2m}.$ It also satisfies $f(0, \ldots, 0) = 1.$

   The proof of the Proposition is an immediate consequence of the
following
lemma.
\vskip 15pt\noindent
\proclaim{Lemma 2.10} The limit $(2.23)$ exists.\endproclaim
\vskip 10pt
{\bf Proof:}
  Let
$$f = \sum_{k \ge 0} f_k$$
where $f_k$ (for $k = 0, 1, 2, \ldots$) is a homogeneous polynomial of
degree
$k$ and which is also symmetric separately with respect to the variables
$u_1,
\ldots, u_m$ and
$u_{m+1}, \ldots, u_{2m}.$ Here $f_0 = 1.$ It follows that

$$H_f = \sum_{k \ge 0} H_{f_k}.$$
Since the action of the group $G'_{2m}$ commutes with permutations of the
variables $u_1, \ldots, u_{2m},$ it easily follows that
$$\sum_{\gamma\in G'_{2m}} H_{f_k}\big(\gamma(u_1, \ldots, u_{2m},
v)\big)$$ is also symmetric separately with respect to the variables $u_1,
\ldots, u_m$ and $u_{m+1}, \ldots, u_{2m}.$

   Define
$$N_{f_k}(u_1, \ldots, u_{2m}, v) =
\left[\underset{ \delta_1 + \cdots + \delta_m =
\delta_{m+1} + \cdots + \delta_{2m}}\to{\prod_{\delta_1=0}^1
\cdots \prod_{\delta_{2m}=0}^1}\left(v + \delta_1 u_1 + \cdots
+ \delta_{2m} u_{2m}\right) \right]\sum_{\gamma \in G'_{2m}}
H_{f_k}\big(\gamma(u_1, \ldots, u_{2m},
v)\big).$$
Then $N_{f_k}$ is invariant under the group $G'_{2m},$ and it is symmetric
separately in the variables $u_1, \ldots,
u_m,$ and $u_{m+1}, \ldots, u_{2m}.$ Moreover, by checking the action of
the
group $G'_{2m}$ on the product
$$\overset m\to{\underset {i = 1}\to{\prod} }\;
\overset 2m\to{\underset {j = m+1}\to{\prod} }
\left(u_i + u_j \right),$$
it follows that $N_{f_k}$ is a rational function
$$N_{f_k} = \frac{N^*_{f_k}}{D^*_{f_k}}\tag 2.24$$
with denominator
$$D^*_{f_k}(u_1, \ldots, u_{2m}, v) = \overset m\to{\underset
{i = 1}\to{\prod} }\; \overset 2m\to{\underset {j = m+1}\to{\prod} }
\left(u_i + u_j \right)\;
\prod_{1 \le i < j \le m} \left(u_i  - u_j\right)
\prod_{m+1 \le i < j \le 2m}\left(u_i  - u_j\right). \tag 2.25$$

The function $N_{f_k}$ is, in fact, a polynomial in the variables $u_1,
\ldots,
u_{2m}, v.$ To see this, we first observe that, for $1 \le i < j \le m$ or
$m+1 \le i < j \le 2m,$
$$N^*_{f_k}(\ldots, u_i, \ldots, u_j, \ldots , v)
= -N^*_{f_k}(\ldots, u_j, \ldots, u_i, \ldots, v).\tag 2.26$$
This implies that
$$N^*_{f_k}(\ldots, u_i, \ldots, u_i, \ldots, v) = 0$$
which gives
$$(u_i - u_j)\; \big | \; N^*_{f_k}(u_1, \ldots,
u_{2m},  v),
  \tag 2.27$$
for $1 \le i < j \le m$ or $m+1 \le i < j \le 2m.$
On the other hand, it can be observed that
$$D^*_{f_k}(u_1, \ldots, u_{2m}, v) =
-D^*_{f_k}\big(\gamma_{ij}(u_1, \ldots, u_{2m}, v)\big),
\tag 2.28$$
for $i = 1, \ldots, m$ and $j = m+1, \ldots, 2m.$ Since the function
$N_{f_k}$
is invariant under the group $G'_{2m},$ it follows from $(2.24),$ and
$(2.28)$
that
$$N^*_{f_k}(u_1, \ldots, u_{2m}, v) =
-N^*_{f_k}\big(\gamma_{ij}(u_1, \ldots,
u_{2m}, v)\big),
  \tag 2.29$$
for $1 \le i < j \le m$ or $m+1 \le i < j \le 2m.$ This together with
$(2.27)$
implies that
$$(u_i + u_j)\; \big | \; N^*_{f_k}(u_1, \ldots,
u_{2m},  v),
  \tag 2.30$$
for $1 \le i \le m$ and $m+1 \le j \le 2m.$
Finally, it follows from $(2.27)$ and $(2.30)$ that for
$\Re(v) > 0,$ the limit
$$\lim_{(u_1, \ldots, u_{2m}) \to (0, \ldots, 0)} \sum_{\gamma\in G'_{2m}}
H_
{f_k}\big(\gamma(u_1, \ldots, u_{2m}, v)\big)$$
exists. Our lemma is proved.

Now, set $u_i = u_{m+i} = i\cdot\epsilon$ (for $i = 1,
2, \ldots, m-1$), $u_m = 0$ and $u_{2m} = m\cdot\epsilon.$ By induction
over
$m,$ it can be checked that
$$\big\{\delta_1u_1 + \cdots
+ \delta_{2m}u_{2m}|\; \delta_{i} = 0, 1;\, \, \delta_1 + \cdots +
\delta_m =
\delta_{m+1} + \cdots + \delta_{2m}\big\} =
\big\{0, 1, \ldots, m ^2 \big\}. \tag 2.31$$
>From Lemma 2.10 and $(2.31),$ it follows that for $k = 0, 1, 2, \ldots,$
$$\sum_{\gamma\in G'_{2m}} H_{f_k}\big(\gamma(u_1, \ldots, u_{2m},
v)\big) = \frac{P_k(\epsilon,
v)}{\overset m ^{2} \to{\underset{\ell = 0}\to{\prod}} (v + \ell\epsilon
)},
\tag 2.32$$
where $P_k(\epsilon, v)$ is a homogeneous polynomial of degree $k$ in
the two variables $\epsilon, v.$

Consequently
$$\lim_{v\to 0}\,\lim_{\epsilon \to 0}\; v^{m ^{2} + 1}
\sum_{\gamma\in G'_{2m}}
H_{f_k}\big(\gamma(u_1, \ldots, u_{2m}, v)\big) = 0$$
if $k > 0,$ and the limit exists if $k = 0.$ Using that $f_0 = 1,$ the
proposition follows by taking the residue at $v=0$ on both sides of
$(2.32).$

\vskip 20pt\noindent
{\bf \S 3. Moments of Quadratic Dirichlet L--Functions}
\vskip 10pt
   Let $$\chi_d(n) = \cases \left(\frac{d}{n}\right) &\;
\text{if $d \equiv 1\pmod{4}$},\\
\left(\frac{4d}{n}\right)\; & \text{if $d \equiv 2, 3
\pmod{4},$}\endcases$$
denote Kronecker's symbol which is precisely the Dirichlet character
associated to the quadratic field $\Bbb Q(\sqrt{d}).$ For $\Re(s) > 1$
we define
$$L(s, \chi_d) \; = \; \sum_{n=1}^\infty \frac{\chi_d(n)}{n^s},$$
to be the classical Dirichlet L--function associated to $\chi_d.$

   We shall always denote by $\sum_{|d|}$ a sum ranging over
fundamental discriminants of quadratic fields.  We shall consider moments
as
$x \to \infty.$ Jutila \cite{J} was the first to obtain the moments
$$\sum_{|d| \le x} L\left({\scriptstyle \frac12},
\,\chi_d\right) \;
\sim
\; a_1 \frac{6}{\pi^2} x\log(x^\frac12) \tag 3.1$$
and
$$\sum_{|d| \le x} L\left({\scriptstyle \frac12},
\,\chi_d\right)^2 \;
\sim
\; 2\cdot \frac{a_2}{3!}  \frac{6}{\pi^2} x\log^3(x^{\frac12})\tag 3.2$$
with
$$a_m \; = \; \prod_p
\frac{\left(1-\frac{1}{p}\right)^{\frac{m(m+1)}{2}}}
{\left(1 +
\frac{1}{p}\right)}
\left(\frac{\left(1-\frac{1}{\sqrt{p}}\right)^{-m} \; + \;
\left(1+\frac{1}{\sqrt{p}}\right)^{-m}}{2} \; + \;
\frac{1}{p}\right),
\quad\qquad (m = 1, 2, \ldots).\tag 3.3$$

  Subsequently, Soundararajan \cite{So} showed that
$$ \sum_{|d| \le x} L\left({\scriptstyle \frac12},
\,\chi_d\right)^3 \;
\sim
\; 16\cdot \frac{a_3}{6!} \frac{6}{\pi^2} x\log^6(x^{\frac12}).\tag 3.4$$
He also conjectured that
$$\sum_{|d| \le x} L\left({\scriptstyle \frac12},
\,\chi_d\right)^4 \;
\sim
\; 768\cdot \frac{a_4}{10!}  \frac{6}{\pi^2} x\log^{10}(x^{\frac12}).\tag
3.5$$

   Motivated by the fundamental work of Katz and Sarnak \cite{K--S},
who introduced symmetry types associated to families of L--functions,
the previous results $(3.1), (3.2), (3.4),
  (3.5),$  and calculations of Keating and Snaith \cite{Ke--Sn--2} based
on
random matrix theory, Conrey and Farmer have made the following
conjecture.

\vskip 10pt
\proclaim{\bf Conjecture 3.1}  For every positive integer $m$,  and $x \to
\infty$,
$$\sum_{|d| \le x} L\left({\scriptstyle \frac12}, \;\chi_d\right)^m \;\;
\sim
\;\;
\frac{6}{\pi^2}\, a_m \cdot \,\prod_{\ell=1}^m \frac{\ell
!}{(2\ell)!}\cdot x
\,(\log x)^M,$$
where $M = {\frac{m(m+1)}{2}}.$
\endproclaim

\vskip 20pt

\noindent
{\bf \S 3.1 The Multiple Dirichlet Series for the Family of Quadratic L-
Functions}
\vskip 10pt

For
$w, s_1, s_2,
\ldots, s_m
\in
\Bbb C$ with
$\Re(w) > 1$ and
$\Re(s_i) > 1\; (i = 1,2,\ldots, m)$,  consider the absolutely convergent
multiple Dirichlet series
$$Z(s_1, s_2, \ldots, s_m, w) = \sum_d \frac{L(s_1, \chi_d)\cdot
L(s_2, \chi_d) \cdots L(s_m, \chi_d)}{|d|^w} \tag 3.6$$
where the  sum ranges over fundamental discriminants of
quadratic fields.

  Recently, (see \cite{B--F--H--1}), for the
special
cases
$m = 1, 2, 3$  a new proof of
Conjecture 3.1, based on the meromorphic continuation of  $Z(s_1,
\ldots, s_m, w)$ , was obtained . Unfortunately, the method of proof
breaks
down when
$m \ge 4$ because there are not enough functional equations of $Z(s_1,
\ldots,
s_m,
w)$ to obtain its meromorphic continuation  slightly beyond
the first significant polar divisor at $w = 1,$ and, $ s_1 \to \frac12,
s_2
\to
\frac12, \ldots, s_m \to \frac12.$

   We shall show that
$Z(s_1,\ldots, s_m, w)$ (suitably modified by breaking it into two
parts and multiplying by appropriate gamma factors) satisfies the
functional equations
$$(s_1,  \ldots, s_m, w)\;\; \overset \alpha_i\to{\longrightarrow}
\;\; (s_1, \ldots, 1-s_i, \ldots, s_m, w + s_i - {\scriptstyle\frac12}),
\qquad\quad (i = 1, 2, \ldots, m). \tag 3.7$$ We then show that for
$\Re(s_i)$
sufficiently large $(i = 1, 2, \ldots, m)$,  that $Z(s_1, \ldots, s_m, w)$
has a simple pole at $w = 1$, and that the residue has analytic
continuation to the region $$\Re(s_i) > \frac12 - \epsilon,
\qquad\qquad (i = 1,
2,
\ldots, m),$$
for any fixed $\epsilon > 0.$ The residue of $Z(s_1, \ldots, s_m, w)$ at
$w = 1
$ and
$s_1 \to \frac12, \;\ldots\;  , s_m \to \frac12$ can be computed exactly
and coincides with the constant in Conjecture 3.1. This is the basis for
Conjecture 3.6 given  in  \S 3.2.

  In order to determine the residues and poles of $Z(s_1, \ldots, s_m,
w)$, it is
necessary to introduce a modified multiple Dirichlet series defined by
$$Z^{\pm}_\nu(s_1, \ldots, s_m, w) =\underset d -\text{sq.free}\to{
\underset  d \equiv\nu
\pmod{4} \to {\underset {\pm d > 0}\to \sum}} \frac{L(s_1, \chi_d)
\cdots L(s_m, \chi_d)}{|d|^w}.\tag 3.8$$
We set
$$Z^\pm(s_1, \ldots, s_m, w) = Z^\pm_1(s_1, \ldots, s_m, w)  +
4^{-w}\Big(Z_2^\pm(s_1, \ldots, s_m, w) + Z_3^\pm(s_1, \ldots,
s_m, w)\Big).\tag 3.9$$

  Further, we define
$$\widehat{Z^+}(s_1, \ldots, s_m, w) = \left(\prod_{i=1}^m
\pi^{-\frac{s_i}{ 2}}\;
\Gamma\left(\frac{s_i}{2}\right)\right)\cdot Z^+(s_1, \ldots, s_m, w)\tag
3.10$$ and
$$\widehat{Z^-}(s_1, \ldots, s_m, w) = \left(\prod_{i=1}^m
\pi^{-\frac{s_i+1}{ 2}}\;
\Gamma\left(\frac{s_i+1}{2}\right) \right)\cdot Z^-(s_1, \ldots,
s_m, w).\tag 3.11$$
   The following two propositions summarize the analytic properties of
the functions $Z^\pm.$
\vskip 10pt
\proclaim{Proposition 3.2}  For $\sigma > 0,$ the functions $Z^\pm$ can
be meromorphically continued
to the domain
$$\Re(s_i) > -\sigma \;
  (i = 1, 2,\ldots, m), \qquad \Re(w) > 1 + m\cdot({\scriptstyle
\frac{1}{2}} +
\sigma
).$$
The only poles in this region are at  $s_i = 1, (i = 1, \ldots,
m)$.  Moreover, both
$\widehat{Z^\pm}$ are invariant under the finite abelian group $G_m$
(of
$2^m$ elements) generated by the involutions
  $$(s_1,  \ldots, s_m, w)\;\; \overset \alpha_i\to{\longrightarrow}
\;\; (s_1, \ldots, 1-s_i,
\ldots, s_m, w + s_i - {\scriptstyle\frac12}), \qquad\quad(i = 1, 2,
\ldots, m).$$ \endproclaim

{\bf Proof:} Note that the term corresponding to $d = 1$ in the definition
of
$Z^\pm$  as a Dirichlet series (see $(3.8), (3.9)$) contributes
$\zeta(s_1)\cdots\zeta(s_m)$ which has poles at $s_i = 1$ for $i = 1,
\ldots,
m$. The functional equation of
$L(s,
\chi_d)$ (see \cite{D})may be written in the form
   $$\align \Lambda(s, \chi_d) &= \pi^{-\frac{s+ a}{2}} \,
\Gamma\left(\frac{s+a}{2}\right) L(s, \chi_d)\tag 3.12\\
&= |D|^{\frac12 - s} \, \Lambda(1 - s, \chi_d), \endalign$$
where $a = 0, 1$ is chosen so that $\chi_d(-1) = (-1)^a,$
and
$$D = \cases d \quad & \text{if}\;\;   d \equiv 1\pmod{4}\\
  4d \quad  & \text{if }\;\, d \equiv 2, 3 \pmod{4}\endcases$$
is the conductor of $\chi_d$. It follows from $(3.12)$ that for $\Re(s)
> -\sigma$, and $d > 1,$
$$|L(s, \chi_d)| = O\Big(|d|^{\frac12 + \sigma}\Big),\tag 3.13$$
where the $O$--constant depends at most on $\Im(s)$. Plugging the
estimate $(3.13)$ into the definition  $(3.8)$ of $Z^\pm_\nu(s_1, s_2,
\ldots, s_m, w)$ (with $\nu = 1, 2, 3$) viewed as an infinite series, we
see
that
the series (with terms $d > 1$) converges absolutely provided $\Re(w)  > 1
+
m\cdot(\frac12 +
\sigma).$ This establishes the first part of Proposition 3.2.

Now, both $\widehat{Z^\pm}$ are invariant under permutations of
the variables
$s_1, s_2, \ldots, s_m$. Therefore, to prove the invariance under the
group $G_m$, it suffices to show the invariance under the
transformation $\alpha_1$, say. To show this invariance, we invoke
the functional equation $(3.12)$ with $s = s_1.$ The invariance under
the transformation $\alpha_1$ immediately follows.

\vskip 10pt\noindent
\proclaim{Proposition 3.3}  The functions $\zeta(2w) Z^\pm$ can be
meromorphically continued for $\Re(w) > 0$ and $\Re(s_i)$
sufficiently large
$(i = 1, 2, \ldots, m)$. They are holomorphic in this region except for
a simple pole at $w = 1$ with residue
$$\align\underset {w = 1}\to{\text{\rm Res}}\Big[\zeta(2w) Z^+(s_1,
\ldots, s_m, w)\Big ]\; &= \;
\underset {w = 1}\to{\text{\rm Res}}\Big[\zeta(2w) Z^-(s_1,
\ldots, s_m, w)\Big]\\
&\\ &= \;   \frac12 \underset {n_1\cdots n_m = \square }\to
{\underset {n_1, \ldots,
n_m }\to
\sum} \frac{\underset{p | n_1\cdots n_m}\to{\prod }\big(1+
p^{-1}\big)^{-1}}{n_1^{s_1}\cdots n_m^{s_m}}.\endalign$$
Here $\square$ denotes any square integer, and the sum ranges over
all $m$--tuples $\{n_1, \ldots, n_m\}$ of positive integers.\endproclaim

\vskip 10pt
{\bf Proof:} It follows from $(3.8)$ that
$$Z_1^\pm(s_1, \ldots, s_m, w) = \sum_{n_1, \ldots, n_m}
\frac{1}{n_1^{s_1}\cdots  n_m^{s_m}}\underset d-\text{sq.free}
\to{\underset{d\,
\equiv\, 1\pmod 4}\to{\sum_{\pm d\; >\; 0}}}
\frac{\chi_d(n_1\cdots n_m)}{|d|^w}. \tag 3.14$$
For any fixed $m$--tuple $\{n_1, \ldots, n_m\}$ of positive integers,
we may write
$$n_1\cdots n_m \; = \; 2^c n N^2 M^2$$
so that
$$\align &\bullet n \;\;\text{is square-free}\\
&\bullet p|N \implies p | n\tag 3.15\\
&\bullet n \;\;\text{and} \;\;  M \;\;\text{are both odd and
coprime}.\endalign$$
It immediately follows from  $(3.15)$ that the inner sum in $(3.14)$
can be rewritten as
$$\align
{\underset d-\text{sq.free} \to{\underset{d\, \equiv\, 1\pmod
4}\to{\sum_{\pm d\; >\; 0}}}}&
\frac{\chi_d(n_1\cdots n_m)}{|d|^w} = \underset (d, M) =
1\to{\underset d-\text{sq.free} \to{\underset{d\, \equiv\, 1\pmod
4}\to{\sum_{\pm d\; >\; 0}}}}
\frac{\chi_d(2)^c\cdot\chi_d(n)}{|d|^w} =
\underset (d, M) =
1\to{\underset d-\text{sq.free} \to{\underset{d\, \equiv\, 1\pmod
4}\to{\sum_{\pm d\; >\; 0}}}}
\frac{\chi_2(d)^c\cdot\chi_n(d)}{|d|^w} \\
& \\
& \tag 3.16\\
& = \frac12 \underset (d, 2M) =
1\to{\underset d -\text{sq.free}\to {\sum_{\pm d \;> \; 0}}}
\frac{\chi_2(d)^c\cdot\chi_n(d)}{|d|^w}  + \frac12 \underset (d, 2M) =
1\to{\underset d -\text{sq.free}\to {\sum_{\pm d \;> \; 0}}}
\frac{\chi_2(d)^c\cdot\chi_{-1}(d)\cdot\chi_n(d)}{|d|^w}.
\endalign$$
Here we have used the law of quadratic reciprocity  $$\chi_d(2)  =
\cases
\chi_2(d) = (-1)^\frac{d^2-1}{8} &\text{if $d \equiv 1 \pmod{4},$} \\
0 & \text{if $d \equiv 2, 3 \pmod{4},$}\endcases$$
and
$$\chi_d(n) = \chi_n(d)\cdot (-1)^\frac{(d-1)(n-1)}{4}, \qquad
\quad (d, n, \;\text{odd}).$$
Further,  for $d$ odd,
$\frac12 \big ( 1 + \chi_{-1}(d)\big )$ is $1$ or $0$ according as
$d \equiv 1 \, \text{or} \; 3\pmod{4}$.  This last assertion follows from
the
identity
$$\chi_{-1}(d) = \left(\frac{-4}{d}\right) = \cases
(-1)^{\frac{|d| - 1}{2}+\frac{\text{sgn}(d) -1}{2}}  & \text{if $d \equiv
1 \pmod{2}$}\\
  0 & \text{if $d \equiv 0 \pmod{2}$}.\endcases$$
In order to complete the proof of Proposition 3.3 we require the following
lemma.
\vskip 10pt
\proclaim{Lemma 3.4}  Let $\chi$ be a primitive quadratic Dirichlet
character of conductor $n$, and let $b$ be any positive integer. If
$L_b(w, \chi)$ is the function defined by
$$L_b(w,\chi) = \underset (d, b) = 1 \to{\underset d-\text{sq.free}
\to{\sum_{d > 0} }} \frac{\chi(d)}{d^w}$$
then $\zeta(2w) L_b(w, \chi)$ can be meromorphically continued to
$\Re (w) > 0.$ It is analytic everywhere in this region, unless $n = 1$
(i.e., $L(w, \chi) = \zeta (w))$, when it has exactly one simple pole at
$w = 1$ with residue
$$\underset w = 1\to{\text{\rm Res}} \,\Big[ \zeta(2w) L_b(w,
\chi)\Big ] \;\; = \;\; \prod_{p | b} \left(1 + \frac{1}{p}\right)^{-1}.$$
\endproclaim
\vskip 10pt
{\bf Proof:} The proof of Lemma 3.4  is a simple consequence of
the elementary identity
$$L_b(w, \chi) = \frac{L(w, \chi) }{\zeta(2w)}\cdot \prod_{p | b}
\left(1 +
\chi(p) p^{-w}\right)^{-1} \prod_{p | n} (1 - p^{2w})^{-1}.$$
  It immediately follows from $(3.16)$ and Lemma 3.4 that
$${\underset d-\text{sq.free} \to{\underset{d\, \equiv\, 1\pmod
4}\to{\sum_{\pm d\; >\; 0}}}}
\frac{\chi_d(n_1\cdots n_m)}{|d|^w} = \frac12 L_{2M}\big (w, \;
\chi_2^c\cdot
\chi_n\big ) \;\; + \;\; \frac12  L_{2M}\big (w, \;
\chi_2^c\cdot \chi_{-1}\cdot
\chi_n\big ),\tag 3.17$$
and that the right hand side of $(3.17)$ has  a meromorphic
continuation to $\Re(w) > 0$. Moreover, it is holomorphic in this
region unless $n = 1$ and $c \equiv 0 \pmod{2}$, in which case
there is exactly one simple pole at $w = 1$ with residue
$$ \underset w = 1\to {\text{Res}}\Bigg[ \zeta(2w)\cdot
\hskip -12pt{\underset d-\text{sq.free}
\to{\underset{d\, \equiv\, 1\pmod 4}\to{\sum_{\pm d\; >\; 0}}}}
\frac{\chi_d(n_1\cdots n_m)}{|d|^w}\Bigg] = \frac12 \prod_{p
| 2M} \left(1 + \frac{1}{p}\right)^{-1}.\tag 3.18$$
Now, if we sum both sides of $(3.18)$ over all $m$--tuples
$\{m_1, \ldots, m_n\}$, it is clear that there will only be a
contribution to the residue coming from $m$--tuples where
$m_1\cdots m_n = \square.$ Combining equations $(3.14)$ and
$(3.18)$, and then removing the factor $1 + 2^{-1}$  when
$n_1\cdots n_m$ is odd gives
$$ \align\underset w = 1\to {\text{Res}} \;\Big[\zeta(2w)&\cdot
Z_1^\pm (s_1,\ldots, s_m, w)
\Big ] = \frac {1}{2} \underset n_1\cdots n_m =
\square\to{\sum_{n_1, \ldots , n_m}} \frac{\underset p|
2n_1\cdots n_m\to {\prod} (1 + p^{-1})^{-1}}{n_1^{s_1}\cdots
n_m^{s_m}}\\
& \tag 3.19\\
&=  \frac {1}{2} \underset n_1\cdots n_m =
\square\to{\sum_{2 \, | \, n_1 \cdots  n_m }}
\frac{\underset p| n_1\cdots n_m \to {\prod} (1 +
p^{-1})^{-1}}{n_1^{s_1}\cdots n_m^{s_m}} \,\; + \;\, \frac {1}{3}
\underset n_1\cdots n_m =
\square\to{\sum_{2\,\nmid \, n_1 \cdots  n_m }}
\frac{\underset p| n_1\cdots n_m\to {\prod} (1 +
p^{-1})^{-1}}{n_1^{s_1}\cdots n_m^{s_m}}\endalign$$
In a completely analogous manner, we can also obtain
$$\underset w = 1\to {\text{Res}} \;\Big[\zeta(2w)\cdot  Z_\nu^\pm
(s_1,\ldots, s_m, w)
\Big ] = \frac {1}{3} \underset n_1\cdots n_m =
\square\to{\underset 2\,\nmid \, n_1 \cdots n_m\to{\sum}}
\frac{\underset p | n_1\cdots n_m\to {\prod} (1 +
p^{-1})^{-1}}{n_1^{s_1}\cdots n_m^{s_m}}\tag 3.20$$
for the cases $\nu = 2, 3.$

  The completion of the proof of Proposition 3.3 now immediately follows
from equations
$(3.9)$,
$(3.19)$ and $(3.20)$ after separating the cases when the product
$n_1\cdots n_m$ is even or odd.

\vskip 10pt

\pagebreak

\proclaim{Proposition 3.5}  Let $\Re (s_i)$ be sufficiently large for $i =
1,
2, \ldots, m.$ Then
$$\underset {w = 1}\to{\text{\rm Res}}\Big[\zeta(2w)\cdot Z^+(s_1,
\ldots, s_m, w)\Big ] \;\; = \;\; \frac12 R(s_1, \ldots, s_m)\cdot
\prod_{i=1}^m \zeta(2s_i) \prod_{1 \le i < j \le m} \zeta(s_i + s_j),$$
where $R(s_1, \ldots, s_m)$ can be holomorphically continued to
the region $\Re(s_i) > \frac12 - \epsilon$ for some fixed $\epsilon > 0$.
Further,
$$ R\left(\frac12,\ldots,
\frac12\right) \;\; = \;\;  a_m,$$
where $a_m$ is the constant given in  $(3.3).$ \endproclaim
\vskip 10pt
{\bf Proof:} If $f(n)$ is a multiplicative function for which the sum
$\overset \infty\to{\underset {n=1}\to{\sum}} f(n)$ converges absolutely,
then
  we have the Euler product identity
$$\overset \infty\to{\underset {n=1}\to{\sum}} f(n) = \prod_p (1 + f(p) +
f
(p^2) +
f(p^3) +
\cdots ).\tag 3.20$$
It now follows from Proposition 3.3 and $(3.20)$ that
$$\underset {w = 1}\to{\text{\rm Res}}\Big[\zeta(2w) Z^+(s_1,
\ldots, s_m, w)\Big ] \;\; = \;\; \frac12 \prod_p\Bigg[1 + \left(1 +
\frac{1}{p}\right)^{-1} \sum_{\mu = 1}^\infty  \;\;\underset  e_i \ge
0, \; (i =
1,\ldots, m)\to  {\underset  e_1+\cdots
+ e_m = 2\mu \to{\sum}} p^{-(e_1s_1 + \cdots + e_ms_m)}\Bigg],$$
where the product converges for $\Re(s_i) > \frac12$, (for $i = 1, 2,
\ldots,m).$
On the other hand, the function $R(s_1, \ldots, s_m)$ defined by
$$\prod_p\left[1 + \left(1 +
\frac{1}{p}\right)^{-1} \sum_{\mu = 1}^\infty \underset  e_i \ge 0, \; (i
=
1,\ldots, m)\to  {\underset  e_1+\cdots
+ e_m = 2\mu \to{\sum}} p^{-(e_1s_1 + \cdots + e_ms_m)}\right ]
\prod_{i=1}^m \zeta(2s_i)^{-1} \prod_{1 \le i < j \le m} \zeta(s_i +
s_j)^{-1}\tag 3.21$$
is holomorphic for $\Re(s_i) > \frac12 - \epsilon$, ($i = 1,
2,
\ldots, m)$ for some fixed small
$\epsilon > 0.$ This establishes the first part of Proposition 3.5.

  Now, the number of terms in the inner sum
$$\underset  e_i \ge 0, \; (i =
1,\ldots, m)\to  {\underset  e_1+\cdots
+ e_m = 2\mu \to{\sum}} p^{-(e_1s_1 + \cdots + e_ms_m)}$$
of  formula $(3.21)$ is precisely
$$d_m\left(p^{2\mu}\right) = \frac{(m + 2\mu  - 1)!}{(m - 1)!\cdot
(2\mu)!}.$$
If we specialize to $s_1 = \cdots = s_m = s$, we get
$$\prod_p\left[ 1 + \left(1 + \frac{1}{p}\right)^{-1} \sum_{\mu = 1}
^\infty\underset
e_i
\ge 0,
\; (i = 1,\ldots, m)\to  {\underset  e_1+\cdots
+ e_m = 2\mu \to{\sum}} p^{-(e_1 + \cdots + e_m)s}\right] = \prod_p\left[1
+
\left(1 + \frac{1}{p}\right)^{-1} \sum_{\mu = 1}^\infty d_m(p^{2\mu})
p^{-2\mu s}\right ].$$
It follows from $(3.21)$ that for $\Re(s) \ge \frac12$,
$$R(s, \ldots, s) =
\prod_p\left[1 +
\left(1 + \frac{1}{p}\right)^{-1} \sum_{\mu = 1}^\infty d_m(p^{2\mu})
p^{-2\mu s}\right ]\cdot \zeta(2s)^{-M},$$
and
$$R\left(\frac12, \ldots, \frac12\right) = \prod_p  \left[\left(1 -
\frac{1}{p}\right)^M
\left(1 +
\left(1 + \frac{1}{p}\right)^{-1} \sum_{\mu = 1}^\infty d_m(p^{2\mu})
p^{-\mu}\right)\right ].$$
If we apply the binomial formula to $\left(1 - p^{-\frac12}\right)^{-m} +
\left
(1
+ p^{-\frac12}\right)^{-m}$ in the definition of $a_m$ given in $(3.3)$
we obtain $R\left(\frac12, \ldots, \frac12\right) = a_m.$ This completes
the
proof of Proposition 3.5.

\vskip 20pt\noindent
{\bf \S 3.2} {\bf A Fundamental Conjecture for the Family of
Quadratic Dirichlet
L--Functions}
  \vskip 10pt
    In view of the invariance of $\widehat{Z^\pm}$  under the
group $G_m$, it follows (as in Section 2.3) from Proposition 3.5 that the
polar
divisors
of $\widehat{Z^\pm}$  must contain the $2^m$ hyperplanes
$$\epsilon_1 s_1 + \cdots + \epsilon_m s_m + w - \frac{\epsilon_1 + \cdots
+
\epsilon_m + 2}{2} \; = \; 0,\tag 3.22$$
where each
$\epsilon_i = 0 \; \text{or}\; 1$ for $i = 1, \ldots, m.$
All the hyperplanes $(3.22)$ pass through the point $\left(\frac12,
\ldots,
\frac12, \,1\right).$  We propose the following conjecture.
\vskip 10pt
\proclaim{Conjecture 3.6} The functions $\widehat{Z^\pm}$  have
meromorphic
continuation to a tube domain in $\Bbb C^{m+1}$ which contains the point
$\left(\frac12, \ldots, \frac12, \, 1\right),$ and both these functions
have
the same polar divisor. The part of the polar divisor passing through
$\left
(\frac12, \ldots, \frac12, \, 1\right)$ consists of all the hyperplanes
$(3.22).$ Moreover, the functions $Z^\pm\left(\frac12, \ldots,
\frac12, w\right)
$ are holomorphic for $\Re(w) > 1.$ \endproclaim
\vskip 10pt
\proclaim{Theorem 3.7}  For $m$ even, Conjecture 3.6 implies the
Keating--Snaith--Conrey--Farmer Conjecture 3.1. \endproclaim
\vskip 10pt
{\bf Proof:} We need to again apply Stark's version of the Wiener--Ikehara
Tauberian theorem as quoted in the proof of Theorem 2.8.
Here we take $S(x) = \sum_{|d| \le x} L(1/2, \chi_d)^m$.
Writing $S(x)$ as a Riemann--Stieltjes integral, it follows by integration by parts,
that
$$\int_1^\infty S(t) \cdot t^{-w} \, \frac{dt}{t} = \frac{1}{w} \sum_{d}
\frac{L\left(\frac12, \chi_d\right)^m}{|d|^w}.$$
Since we have assumed $m$ to be even, it follows from $(3.8), (3.9)$ that
$Z^\pm\left(\frac12, \ldots, \frac12, w\right)$ is a Dirichlet series
satisfying
the conditions of the Tauberian theorem. To prove Conjecture 3.1, it is
enough
to show that
$$\lim_{w\to 1} (w - 1)^{M+1} \; Z^\pm\left(\frac12, \ldots, \frac12,
w\right)
= \frac{3}{\pi^2} g_m \, a_m\, M!,$$
where  $$M = \frac{m(m+1)}{2}, \qquad\qquad g_m = \prod_{\ell = 1}^m
\frac{\ell !}{(2\ell)!},$$
and $a_m$ is the constant given in $(3.3)$.

   Let $T(s_1, \ldots, s_m, w)$ denote the function defined by
$$\frac{1}{2(w-1)} \; R(s_1, \ldots, s_m)\; \prod_{i=1}^m \pi^{-\frac{s_i
+ a}
{2}}\;
\Gamma\left(\frac{s_i + a}{2}\right) \zeta(2s_i) \; \prod_{1\le i < j \le
m}
\zeta(s_i + s_j),\tag 3.23$$
where $a = 0, 1$ is determined by $(-1)^a = \pm 1.$
Then Conjecture 3.6 implies that
$$\zeta(2w) \widehat{Z^\pm}(s_1, \ldots, s_m, w) - \sum_{\alpha\in G_m}
T\big(\alpha(s_1, \ldots, s_m, w)\big)\tag 3.24$$
is holomorphic around $\left(\frac12, \ldots, \frac12, \, 1\right).$
The proof of theorem 3.7 is an immediate consequence of the following
proposition.

\vskip 10pt
\proclaim{Proposition 3.8}  For $m = 1, 2, 3, \ldots,$ let $G_m$ denote
the
direct
product of $m$ groups of order $2$
  generated by the involutions $(3.7)$. Let
$$U(s_1, \ldots, s_m, w) =\frac{1}{w-1}  R(s_1, \ldots,
s_m)\;\prod_{i=1}^m \zeta(2s_i) \; \prod_{1\le i < j \le m}
\zeta(s_i + s_j). $$
  Then we have
$$\lim_{w\to 1} \;\;\lim_{(s_1, \ldots, s_m) \to (\frac12,\ldots,
\frac12)}
\Bigg[ (w - 1)^{M+1}\sum_{\alpha\in G_m}  U\big(\alpha(s_1,
\ldots, s_m, w)\big) \Bigg]   \;\; = \;\; \frac{6}{\pi^2}\; a_m\; g_m\; M!
$$
where  $$M = \frac{m(m+1)}{2}, \qquad\qquad g_m = \prod_{\ell = 1}^m
\frac{\ell !}{(2\ell)!},$$   and $a_m$ is the constant given in $(3.3)
$.\endproclaim

\vskip 10pt
{\bf Proof:} We start by taking the Taylor expansion of
$$U(s_1, \ldots, s_m, w) \; = \; \frac{ a_m}{w-1}\;\cdot \frac{f^*(s_1,
\ldots,
s_m)}{\overset m\to{\underset {i = 1}\to{\prod } }(2s_i-1) \underset {1
\le i <
j \le
m}\to{\prod} (s_i+s_j-1)}\tag 3.25 $$
around $(s_1, \ldots, s_m) = \left(\frac12, \ldots,
\frac12\right).$ Here
$$f^*(s_1, \ldots, s_m) = 1 \; + \; \underset{\ell_1 + \cdots +
\ell_{m} \geq 1}
\to{\sum_{\ell_1=0}^\infty \cdots
\sum_{\ell_m=0}^\infty }\kappa_m(\ell_1, \ldots, \ell_m)\;
(s_1-{\scriptstyle\frac12})^{\ell_1}
\cdots (s_m-{\scriptstyle\frac12})^{\ell_m},$$
(with $\kappa_m(\ell_1, \ldots, \ell_m) \in \Bbb C$),
  will be a holomorphic function
which is  symmetric function with respect to the variables $s_1, \ldots,
s_m.$

Now,  make the change of variables $s_i = \frac12 + \epsilon_i$ for $i
= 1, 2, \ldots, m,$ and $w = v + 1.$ The involutions $(3.7)$ are
transformed to
  $$ (\epsilon_1,  \ldots, \epsilon_i,\ldots, \epsilon_m, v)\;\; \overset
\alpha_i\to{\longrightarrow}
\;\; (\epsilon_1,
\ldots, -\epsilon_i,
\ldots, \epsilon_m, v + \epsilon_i ), \qquad\quad (i = 1, 2, \ldots, m).$$
Henceforth, we denote by $G_m$ the group generated by the above
involutions.

Then by $(3.25)$, it is enough to prove that
$$\lim_{v\to 0} \;\;\lim_{(\epsilon_1, \ldots, \epsilon_m) \to (0,\ldots,
0)}
\Bigg[ v^{M+1}\sum_{\alpha\in G_m}  H_f\big(\alpha(\epsilon_1,
\ldots, \epsilon_m, v)\big) \Bigg]
\;\; = \;\;  2^m\; g_m\; M!, \tag 3.26$$
where
$$H_f(\epsilon_1, \ldots, \epsilon_m, v) \;\; = \;\; \frac{1}{v}\;\cdot
\frac{f(\epsilon_1,
\ldots, \epsilon_m)}{\overset m\to{\underset {i = 1}\to{\prod} }\epsilon_i
\underset {1 \le i < j \le m}\to{\prod} (\epsilon_i + \epsilon_j)},$$
and $f$ (which is simply related to $f^*$) is a certain holomorphic
symmetric
function  with respect to the variables
$\epsilon_1, \ldots, \epsilon_m.$ It satisfies $f(0, \ldots, 0) = 1.$

   The proof of Proposition 3.8 is an immediate consequence of  the
following
two lemmas.
\vskip 10pt
\proclaim{Lemma 3.9} The limit $(3.26)$ exists.\endproclaim

{\bf Proof:}
  Let
$$f = \sum_{k \ge 0} f_k$$
where $f_k$ (for $k = 0, 1, 2, \ldots$) is a symmetric and homogeneous
polynomial of degree $k$. Here $f_0 = 1.$ It follows that
   $$H_f = \sum_{k \ge 0} H_{f_k}.$$
Since the action of the group $G_m$ commutes with permutations of the
variables $\epsilon_1, \ldots, \epsilon_m,$ it easily follows that
$$\sum_{\alpha\in G_m} H_{f_k}\big(\alpha(\epsilon_1, \ldots, \epsilon_m,
v)\big)$$ is also a symmetric function with respect to $\epsilon_1,
\ldots,
\epsilon_m.$

   Define
$$N_{f_k}(\epsilon_1, \ldots, \epsilon_m, v) = \left(\prod_{i=1}^m
\epsilon_i
\prod_{1 \le i < j \le m} \left(\epsilon_i^2  - \epsilon_j^2\right)\right)
\sum_{\alpha\in G_m} H_{f_k}\big(\alpha(\epsilon_1, \ldots, \epsilon_m,
v)\big).$$
Then $N_{f_k}$ is a symmetric function in the variables  $\epsilon_1,
\ldots,
\epsilon_m,$ and it is a rational function
$$N_{f_k} = \frac{N^*_{f_k}}{D^*_{f_k}}\tag 3.27$$
  with denominator of the form
$$D^*_{f_k}(\epsilon_1, \ldots, \epsilon_m, v)  = \prod_{\delta_1=0}^1
\cdots
\prod_{\delta_m=0}^1
\left(v +
\delta_1\epsilon_1 + \cdots + \delta_m\epsilon_m\right). \tag 3.28$$
It follows  that
$$N_{f_k}(\ldots, \epsilon_i, \ldots, \epsilon_j, \ldots , \, v) =
-N_{f_k}
(\ldots,
\epsilon_j,
\ldots, \epsilon_i,
\ldots , \, v),\tag 3.29$$
which implies that
$$N_{f_k}(\ldots, \epsilon_i, \ldots, \epsilon_i, \ldots, v) = 0.$$
This gives
$$(\epsilon_i - \epsilon_j)\; \big | \; N^*_{f_k}(\epsilon_1, \ldots,
\epsilon_m,  v).
\tag 3.30$$
Furthermore, since $\underset {\alpha\in G_m}\to{\sum}
H_{f_k}\big(\alpha(\epsilon_1,
\ldots,
\epsilon_m, v)\big)$ is invariant under the group $G_m$, it follows
that
$$N_{f_k}(\ldots, s_i, \ldots , \, v) = -N_{f_k}(\ldots, -s_i, \ldots , \,
v +
s_i)\tag 3.31$$ which implies that
$$N_{f_k}( \ldots, 0, \ldots , v) = 0.$$
Consequently,
$$s_i \; \big | \; N^*_{f_k}(\epsilon_1, \ldots, \epsilon_m, v). \tag
3.32$$
Also, in the same manner, $(3.29)$ and $(3.31)$ imply that
$$(s_i + s_j) \; \big | \; N^*_{f_k}(\epsilon_1, \ldots, \epsilon_m, v).
\tag
3.33$$
Finally, it follows from $(3.27), (3.28), (3.30), (3.32), $ and $(3.33)$
that
for
$\Re(v) > 0$, the limit
$$\lim_{(\epsilon_1, \ldots, \epsilon_m) \to (0, \ldots, 0)}
\sum_{\alpha\in
G_m} H_{f_k}\big(\alpha(\epsilon_1, \ldots, \epsilon_m,
v)\big)$$
exists. It further follows that if we set $\epsilon_i = i\cdot\epsilon$
(for $i
= 1,
2, \ldots, m$) then for $k = 0, 1, 2, 3 \ldots, $
$$\sum_{\alpha\in G_m} H_{f_k}\big(\alpha(\epsilon_1, \ldots, \epsilon_m,
v)\big) = \frac{P_k(\epsilon,
v)}{\overset M \to{\underset{\ell = 0}\to{\prod}} (v + \ell\epsilon
)},\tag
3.34$$ where $P_k(\epsilon, v)$ is a homogeneous polynomial of degree $k$
in
the two variables $\epsilon, v.$

Consequently
$$\lim_{v\to 0}\,\lim_{\epsilon \to 0}\; v^{M+1} \sum_{\alpha\in G_m}
H_{f_k}\big(\alpha(\epsilon_1, \ldots, \epsilon_m, v)\big) = 0$$
if $k > 0,$ and the limit exists if $k = 0.$ This completes the proof of
Lemma
3.9.
\vskip 10 pt\noindent
\proclaim{Lemma 3.10}  Let
$$H(\epsilon_1, \ldots, \epsilon_m, v) = \frac{1}{v}\cdot
\frac{1}{\overset m\to
{\underset {i = 1}\to{\prod} }\epsilon_i
\underset {1 \le i < j \le m}\to{\prod} (\epsilon_i + \epsilon_j)}.$$
Then
$$\lim_{v\to 0} \, \lim_{(\epsilon_1, \ldots, \epsilon_m) \to (0, \ldots,
0)}
\sum_{\alpha\in G_m} H\big(\alpha(\epsilon_1, \ldots, \epsilon_m, v)\big)
=
2^m\, g_m\, M!.$$\endproclaim
\vskip 10pt
{\bf Proof:} We know  from Lemma 3.9 that the above limit exists, so we
can
compute the limit by setting $\epsilon_j = j\epsilon$ (for $j = 1, 2,
\ldots,
m)$ and letting $\epsilon \to 0.$ It follows from $(3.34)$ that
  $$ \sum_{\alpha\in G_m} H\big(\alpha(\epsilon, 2\epsilon,  \ldots,
m\epsilon, v)\big) = \frac{\kappa_m}{\overset M \to{\underset{\ell =
0}\to{\prod}} (v + \ell\epsilon )}$$
for some constant $\kappa_m$. By taking the residue at $w=0$ on both
sides, we
have
$$ \frac{1}{m!}
\prod_{1 \le i < j \le m}\frac{1}{\left(i  + j \right)}=
\frac{\kappa_m}{M!}.$$

By induction over $m,$ one can show that $\kappa_m = 2^m\,
g_m\, M!,$ and the lemma follows.

\vskip 20pt\noindent
{\bf \S 4. Cubic moments of quadratic $L$-series}
\vskip 10pt
As mentioned in the introduction, in the particular cases when $m \le
3$ it is possible to
define an analog of the multiple Dirichlet series given in (3.6).  In
this analog the sum
is not restricted to fundamental discriminants, but ranges over all
integers $d$.  When
an appropriate definition is given for $\prod_{i=1}^m L(s_i,\chi_d)$ for general $d$
one can extend the
multiple Dirichlet series to a meromorphic function of
$s_1,s_2,\dots,s_m,w$ in $\C^{m+1}$.  In this section we will
explicitly provide this
continuation in the case $m = 3$ and $s_1 = s_2 = s_3 = s$.  This
work relies heavily on the
results of \cite{B--F--H--1}. We will then develop a sieving
method analogous to that used in \cite{G--H} to isolate fundamental
discriminants and will
prove as a consequence Theorem 1.1.

\vskip 20pt

\noindent
{\bf  \S 4.1 Some foundations}
\vskip  10pt
The $L$ series $\zeta (s)^3$ can actually be associated
to a certain Eisenstein series $F$ on $GL(3)$, and $L(s,F) = \zeta (s)^3$.

For future convenience, we will write
$$
L(s,F) = \sum_{1}^\infty \frac{c(n)}{n^s},
\tag 4.1
$$
where $c(n) = \sum_{d_1d_2d_3 = n}1$, and we have the
Euler product decomposition
$$
L(s,F) = \prod_p \big( 1 - p^{-s}\big)^{-3},
\tag 4.2
$$
the product being over all primes $p$ of $\Q$.

As in the previous sections, let $\chi_d$ denote the primitive
quadratic character
associated to the quadratic field $\Q(\sqrt{d})$.  If $F$ is twisted
by $\chi_d$,
then the associated $L$-series becomes
$$
L(s,F,\chi_d) = L(s,\chi_d)^3 = \prod_p \big( 1 -
\chi_d(p)p^{-s}\big)^{-3},
\tag 4.3
$$
and by (3.12) the functional equation is given by
$$
(|D|^3)^{s/2}G_{d}(s) L(s,F,\chi_d)  =
   (|D|^3)^{(1-s)/2}G_{d}(1-s) L(1-s,F,\chi_d).
\tag 4.4
$$
Here $D = 4d$ or $D=d$ is the conductor of $\chi_d$ and $G_d(s)$ denotes the
product of gamma factors.

The gamma factors of $(4.4),$ described in $(3.12),$ depend only on the
sign of $d$.
Although we will not require many explicit properties of the gamma
factors, the following
upper bound will be convenient.  For $\sigma_1>\sigma_2$ and $t$ real,
it follows from
Stirling's formula that for large $|t|$, independent of $d$,
$$
\frac{|G_{d}(\sigma_1+it)|}{|G_{d}(\sigma_2-it)|} \ll
(|t| +1)^{3(\sigma_1 - \sigma_2)/2}.
\tag 4.5
$$

When all primes are included in the product $(4.3)$ the functional
equation $(4.4)$ has its optimal form.  However, it is often convenient to
omit factors
corresponding to ``bad" primes, for example those contained in  $S$,
a finite set of
primes including $2$. Let $M=\prod_{p\in S}p$.    For such $M, S,$ we denote
the
$L$-series with Euler factors corresponding to primes dividing $M$
removed as follows:
$$
L_M(s,F) = \prod_{p\notin S} \big( 1 - p^{-s}\big)^{-3}
= L(s,F) \prod_{p\in S} \big( 1 - p^{-s}\big)^3.
\tag 4.6
$$

When twisted by $\chi_d$,  the $L$-series $L(s,F,\chi_d)$ will have a
perfect
functional equation of the form $(4.4)$ when $\chi_d$ is a primitive
character.
This corresponds to the case where $d$ is square free.  It is very
interesting to note that
often, when $d$ is {\it not} square free, it is possible to complete
$L(s,F,\chi_d)$ by
multiplying by a certain Dirichlet polynomial in such a way that the
resulting product has a
functional equation of precisely the same form $(4.4),$ with $D$
replaced by $|d|$
or $|4d|$.  For the simplest example, with $m=1$, see \cite{G--H}.
What is more
remarkable is the fact that some very stringent additional conditions can
be
imposed on the Dirichlet polynomial.

To be more precise, let $l_1, l_2>0, \; l_1,l_2|M$, and $a_1,a_2 \in \{1,-1\}$
and
let $\chi_{a_1l_1},\chi_{a_2l_2}$ be the quadratic characters
corresponding to
$a_1l_1,a_2l_2$ as defined above.  We then formulate the following
collection of properties for two classes of Dirichlet polynomials
associated to $F$.

\pagebreak

\proclaim{Property 4.1}
For $n,d$ positive integers, $(nd,M)=1$, we write  $d=d_0
d_1^2,
\, n = n_0 n_1^2$, with $d_0,n_0$ square free and $d_1,n_1$ positive.  Let $c(n)$
denote the coefficients of
$L(s,F)$ as  defined earlier.

  For complex numbers $A_{d,p^e}^{(\alpha)}, \;
B_{d, p^e}^{(\alpha)}$ (depending on  $d, \alpha \in \Bbb Z, 1 \le
e \le
\alpha$),
let
$P_{d_0,d_1}^{(a_1l_1)}(s)$, $Q_{n_0,n_1}^{(a_2l_2)}(w)$ be Dirichlet
polynomials
defined by
$$
P_{d_0,d_1}^{(a_1l_1)}(s) = \prod_{p^\a ||d_1}
\big(1 + A_{d_0\cdot a_1 l_1, \; p}^{(\a)}\; p^{-s} + \dots +
A_{d_0\cdot a_1 l_1, \; p^{6\a}}^{(\a)} \; p^{-{6\a s}}\big)
$$
and
$$
c(n_0n_1^2)Q_{n_0,n_1}^{(a_2l_2)}(w)= c(n_0n_1^2)\prod_{p^\b ||n_1}
\big(1 + B_{n_0\cdot a_2 l_2, \; p}^{(\b)} \; p^{-w} + \dots +
B_{n_0\cdot a_2 l_2, \; p^{2\b}}^{(\b)} \; p^{-{2\b w}}\big).
$$
We say that $P,Q$ satisfy the conditions of Property 4.1 if the following identities
hold:
$$
d_1^{3s} P_{d_0,d_1}^{(a_1l_1)}(s) = d_1^{3(1-s)}
P_{d_0,d_1}^{(a_1l_1)}(1-s),
\tag 4.7
$$
$$
n_1^wc(n_0n_1^2)Q_{n_0,n_1}^{(a_2l_2)}(w)=
n_1^{1-w}c(n_0n_1^2)Q_{n_0,n_1}^{(a_2l_2)}(1-w)
\tag 4.8
$$
$$
P_{d_0l_3, d_1} ^{(a_1l_1)}(s) = P_{d_0, d_1} ^{(a_1l_1l_3)}(s), \qquad
Q_{n_0l_3, n_1} ^{(a_2l_2)}(w) = Q_{n_0, n_1} ^{(a_2l_2l_3)}(w),
\tag 4.9
$$
(where $d_0l_3,$ $n_0l_3$ are positive square free numbers), and if in
addition, the following interchange of summation is valid for
$s$ and $w$ having sufficiently large real parts:
$$
\multline
\sum_{(d,M)=1} \frac{L_M(s,F,\chi_{d_0}\chi_{a_1l_1})\chi_{a_2l_2}(d_0)
P_{d_0,d_1}^{(a_1l_1)}(s)}{d^w} \\
= \sum_{(n,M)=1}
\frac{L_M(w,\tilde\chi_{n_0}\chi_{a_2l_2})\chi_{a_1l_1}(n_0)
c(n_0n_1^2)Q_{n_0,n_1}^{(a_2l_2)}(w)}{n^s}.
\endmultline
\tag 4.10
$$
Here $\tilde\chi_{n_0}$ denotes the quadratic character with
conductor $n_0$ defined by $\tilde\chi_{n_0}(*) =
\left(\frac{*}{n_0}\right)$. (Recall $2|M$,  so $(2,n_0)=1$.)
\endproclaim

It was observed in \cite{B--F--H--1} that the three properties $(4.7),$
$(4.8)$ and $(4.10)$ were sufficient to determine the polynomials $P$ and
$Q$,
precisely, in the cases of
$GL(1), GL(2), GL(3)$.  This unique determination of $P$ and $Q$ corresponded to a
finite
group of functional equations of the double Dirichlet series given in
(4.10) and
this in turn made it possible to obtain an analytic continuation of the
double
Dirichlet series in these three cases.   It was also noted that for
$m\ge 4$ the corresponding group of functional equations becomes
infinite and that
simultaneously the polynomials $P,$ $Q$ are no longer uniquely
determined by the properties
(4.7), (4.8), and (4.10).  The space of local solutions becomes 1
dimensional in the case
$m=4$, and higher for $m>4$.

In \cite{B--F--H--1} a complete description of certain factors of the
polynomials $P,$ $Q$ was obtained for the case of $m=3$ and an arbitrary
automorphic form $f$
on $GL(3)$. These were the factors corresponding to the ``good"
primes, i.e., primes not
dividing 2 or the level of $f$.  It was also verified that for sums over
positive integers $n,$ $d$ relatively prime to the ``bad" primes, the relations
(4.7), (4.8), (4.9), and (4.10) hold.   In addition, it was verified that for
fixed $d = d_0d_1^2, n
= n_0n_1^2$ and $\epsilon >0$, $\Re s \ge \frac12, \Re w \ge \frac12,$
$$
P_{d_0,d_1}^{(a_1l_1)}(s) \ll |d|^\epsilon \quad \hbox{and} \quad
c(n_0n_1^2)Q_{n_0,n_1}^{(a_2l_2)}(w) \ll |c(n)||n|^\epsilon.
\tag 4.11
$$
In both cases the implied constant depends only on $\epsilon$.   This
information was
then used to obtain the analytic continuation of the double Dirichlet
series on the left
hand side of (4.10).  As a consequence, non vanishing results for
quadratic twists of
$L(\frac12,f,\chi_d)$ were obtained and also, after taking a residue at
$w=1$, a new proof was
obtained for the analytic continuation of the symmetric square of $f$.

As the technique is new, there may be some advantage to presenting the
details of the
analytic continuation argument specialized to the very concrete case where
$L(s,f,\chi_d) = L(s,F,\chi_d) = L(s, \chi_d)^3$,
and we will do so below.
\vskip 20pt\noindent
{\bf  \S 4.2 The cubic moment, continued}
\vskip  10pt
Our object will be to obtain the analytic continuation in $(s, w),$ with
$\Re(s)\ge \frac{1}{2},$ $\Re (w) > \frac{4}{5},$ and an estimate for the
growth in vertical strips $w = \nu+ it$ (for fixed $\nu$ and $s$) of the
double Dirichlet series
$$
Z(s, w) = \sum_{D = \,{\scriptscriptstyle\text{fund.
disc.}}}\frac{L(s,\chi_D)^3}{|D|^w}.
\tag 4.12
$$

To accomplish this, we will obtain the analytic properties of a building
block: For $l_1,$ $l_2>0,$ $l_1,$ $l_2|M$ and $a_1,$ $a_2 \in \{1,-1\},$
we define
$$
Z_M(s,w;\chi_{a_2l_2},\chi_{a_1l_1}) = \sum_{(d,M)=1}
\frac{L_M(s,F,\chi_{d_0}\chi_{a_1l_1})\chi_{a_2l_2}(d_0)
P_{d_0,d_1}^{(a_1l_1)}(s)}{d^w},
\tag 4.13
$$
where we recall that we sum over $d \ge 1$ and use the decomposition
$d = d_0d_1^2,$ with $d_0$ square free and $d_1$ positive.

The following proposition will provide a useful way of collecting the
properties of the multiple Dirichlet series
$Z_M(s, w; \chi_{a_2l_2}, \chi_{a_1l_1}).$
For a positive integer $M,$ define
$$
\text{Div}(M) = \Big\{ a\cdot l \; \Big | \; a = \pm 1, \; 1 \le l, \; l |
M\Big\},
$$
which has cardinality
$2d(M)=
2\sum_{d|M}1.$ Let
$\overarrow \Z_M(s,w;\chi_{a_2l_2}, \chi_{\text{Div}(M)})$ denote the
$2d(M)
$ by 1
column vector whose $j^{th}$ entry is
$Z_M(s, w; \chi_{a_2l_2},  \chi^{(j)}),$ where $\chi ^{(j)}$ $(j = 1, 2,
\ldots, 2d(M))$ ranges over the characters $\chi_{a_1 l_1}$ with $a_1  =
\pm 1, \; 1 \le  l_1, \; l_1 | M.$
Then, we will prove

\vskip 10pt
\proclaim{Proposition 4.2}
There exists a $2d(M)$ by $2d(M)$ matrix $\Phi^{(a_2l_2)} (w)$ such that
for any fixed $w,$ $w\ne 1,$ and for any $s$ with sufficiently large real
part
(depending on $w$)
$$
\prod_{p|(M/l_2)}\big(1-p^{-2+2w}\big)\cdot
\overarrow\Z_M(s,w;\chi_{a_2l_2},\chi_{\text{Div}(M)}) =\Phi^{(a_2l_2)}
(w)
\overarrow\Z_M(s+w-1/2,1-w;\chi_{a_2l_2},\chi_{\text{Div}(M)}).
$$
The entries of $\Phi^{(a_2l_2)} (w),$ denoted by $\Phi_{i,j}^{(a_2l_2)}
(w),$
are meromorphic functions in $\Bbb C.$
\endproclaim
\vskip 10pt
{\bf Proof:} By Property 4.1,
$$
Z_M(s,w;\chi_{a_2l_2},\chi_{a_1l_1}) = \sum_{(n,M)=1}
\frac{L_M(w,\tilde\chi_{n_0}\chi_{a_2l_2})\chi_{a_1l_1}(n_0)
c(n_0n_1^2)Q_{n_0,n_1}^{(a_2l_2)}(w)}{n^s}.
\tag 4.14
$$
Now
$$
L_M(w,\tilde\chi_{n_0}\chi_{a_2l_2}) =
L(w,\tilde\chi_{n_0}\chi_{a_2l_2})\cdot
\prod_{p|M}\big(1- \tilde\chi_{n_0}\chi_{a_2l_2}(p)p^{-w}\big),
\tag 4.15
$$
where $L(w,\tilde\chi_{n_0}\chi_{a_2l_2})$ satisfies the functional
equation
$$
\multline
G_\epsilon(w) (n_0l_2D_{a_2l_2})^{w/2}L(w,\tilde\chi_{n_0}\chi_{a_2l_2})=
G_\epsilon(1-w)
(n_0l_2D_{a_2l_2})^{(1-w)/2}L(1-w,\tilde\chi_{n_0}\chi_{a_2l_2}).
\endmultline
\tag 4.16
$$
Here $\epsilon = \tilde\chi_{n_0}\chi_{a_2l_2}(-1),$
$$
G_\epsilon(w)=\cases
\pi^{-w/2}\Gamma (w/2)&\text{if \; \;  $\epsilon = 1$}\\
\pi^{-(w + 1)/2}\Gamma ((w+1)/2)&\text{if\; \; $\epsilon = -1$,}\endcases
\tag 4.17
$$
and
$$
D_{a_2l_2}=\cases
1&\text{if \; \; $a_2l_2 \equiv 1 \pmod 4$}\\
4 &\text{otherwise.}\endcases$$

Combining this with the functional equation for $Q$ given in $(4.8),$ we
obtain
$$
\multline
Z_M(s,w;\chi_{a_2l_2},\chi_{a_1l_1})\\ =\sum_{a_3 = 1,-1}
\sum_{(n,M)=1,\,n\equiv a_3\,(4) }\frac{G_{\epsilon
(a_3a_2l_2)}(1-w)(l_2D_{a_2l_2})^{1/2
-w}}{G_{\epsilon(a_3a_2l_2)}(w)n^{s+w-1/2}}\\ \times
\chi_{a_1l_1}(n_0) L_M(1-w,\tilde\chi_{n_0}\chi_{a_2l_2})
c(n_0n_1^2)Q_{n_0,n_1}^{(a_2l_2)}(1-w) \cdot
\prod_{p|(M/l_2)}\big(1- \tilde\chi_{n_0}\chi_{a_2l_2}(p)p^{-w}\big)\\
\times
\prod_{p|(M/l_2)}\big(1-
\tilde\chi_{n_0}\chi_{a_2l_2}(p)p^{-1+w}\big)^{-1}.
\endmultline
$$
Here $\epsilon(a)$ denotes the sign of $a$.
Note that we are leaving out terms in the product where $p|l_2$ as the
character vanishes here.

Multiplying by $\prod_{p|(M/l_2)}\big( 1-p^{-2+2w}\big)$ and reorganizing,
we
obtain
$$
\multline
\prod_{p|(M/l_2)}\big( 1-p^{-2+2w}\big) \cdot
Z_M(s,w;\chi_{a_2l_2},\chi_{a_1l_1})\\ =
\sum_{a_3=1,-1}\frac{G_{\epsilon(a_3a_2l_2)}(1-w)}
{G_{\epsilon(a_3a_2l_2)}(w)(l_2D_{a_2l_2})^{w-1/2}}
\sum_{l_3,l_4|(M/l_2), \; (l_4, 2) = 1}\mu
(l_3)\chi_{a_2l_2}(l_3l_4)l_3^{-w}l_4^{-1+w}\\ \times
\sum_{(n,M)=1,\,n\equiv a_3\,(4)}
\frac{A_{2}(1 - w,
\tilde\chi_{n_0}\chi_{a_2l_2})L_M(1-w,\tilde\chi_{n_0}\chi_{a_2l_2})
c(n_0n_1^2)Q_{n_0,n_1}^{(a_2l_2)}(1-w)\chi_{a_1l_1l_3l_4}(n_0)}{n^{s+w-1/2}},
\endmultline
$$
where
$$
A_2(w, \tilde\chi_{n_0}\chi_{a_2l_2}) = \cases
1&\text{if\; \; $2 | l_2,$}\\
1 + \tilde\chi_{n_0}\chi_{a_2l_2}(2)2 ^ {-w}&\text{if\; \; $a_2l_2\equiv 1
\pmod{4},$}\\
1 - 2 ^ {-2w}&\text{if\; \; $a_2l_2\equiv -1 \pmod{4}.$}\endcases$$
We have used here the fact that
$\tilde\chi_{n_0}(l_3)\tilde\chi_{n_0}(l_4)=\chi_{l_3l_4}(n_0),$
and the identity
$$
\big(1 - 2 ^{-2+2w} \big) \big(1- \tilde\chi_{n_0}\chi_{a_2l_2}(2)2^{-1 +
w}\big) ^{-1}
= A_2(1 - w, \tilde\chi_{n_0}\chi_{a_2l_2}),
$$
for $a_2l_2\equiv -1,\;1 \pmod{4}.$

Using $\chi_{-1}$ to sieve congruence classes of $n \pmod 4:$
$$
\frac{1}{2}\big(1 + a_3 \chi_{-1}(n_0) \big)=\cases
1&\text{if \; \; $n_0 \equiv a_3 \pmod 4$}\\
0&\text{if \; \; $n_0 \equiv -a_3 \pmod 4,$}\endcases
$$
we finally obtain (in the case of $a_2l_2\equiv 1 \pmod{4}$)
$$
\multline
\prod_{p|(M/l_2)}\big( 1-p^{-2+2w}\big) \cdot
Z_M(s,w;\chi_{a_2l_2},\chi_{a_1l_1})\\ =
\frac{1}{2}\cdot l_2^{1/2-w} \cdot \sum_{l_3,l_4|(M/l_2)}\mu
(l_3)\chi_{a_2l_2}(l_3l_4)l_3^{-w}l_4^{-1+w}\sum_{a_3=1,-1}
\frac{G_{\epsilon(a_3a_2l_2)}(1-w)}
{G_{\epsilon(a_3a_2l_2)}(w)}
\\ \times
\big( Z_M(s+w-1/2,1-w;\chi_{a_2l_2},\chi_{a_1l_1l_3l_4}) +
a_3  Z_M(s+w-1/2,1-w;\chi_{a_2l_2},\chi_{-a_1l_1l_3l_4}) \big).
\endmultline
\tag 4.18
$$

If $a_2l_2\equiv -1,\; 2 \pmod{4},$ we have a similar expression.
Actually,
it can be easily observed that just the behavior at the finite place $2$
changes.

This completes the proof of Proposition 4.2.

The function $Z_M(s,w;\chi_{a_2l_2},\chi_{a_1l_1})$ defined in $(4.13)$
also
possesses a functional equation as $s \rightarrow 1-s.$ To describe this,
let $d(M)$ be as before, and let
$\overarrow\Z_M(s,w;\chi_{\text{Div}(M)},\chi_{a_1l_1})$ denote the
$2d(M)$
by 1 column vector
whose $j^{th}$ entry is $Z_M(s,w; \chi ^{(j)}, \chi_{a_1l_1}),$ where
$\chi^{(j)} \;(j = 1, 2, \ldots,  2d(M))$ ranges over the characters
$\chi_
{a_2 l_2}$ with $a_2  = \pm 1, \; 1 \le  l_2, \; l_2 | M.$

Then we have the following.
\vskip 10pt
\proclaim{Proposition 4.3}
There exists a $2d(M)$ by $2d(M)$ matrix $\Psi^{(a_1l_1)} (s)$ such that
for any fixed $s,$ $s\ne 1,$ and for any $w$ with sufficiently large real
part
(depending on $s$)
$$
\overarrow\Z_M(s,w;\chi_{\text{Div}(M)},\chi_{a_1l_1})\cdot\prod_{p|(M/l_1)}
\big(1-p^{-2+2s}\big)^3
=\Psi^{(a_1l_1)} (s)
\overarrow\Z_M(1-s,w+3s-3/2;\chi_{\text{Div}(M)},\chi_{a_1l_1}).
$$
The entries of $\Psi^{(a_1l_1)} (s),$ denoted by $\Psi_{i,j}^{(a_1l_1)}
(s),$ are meromorphic functions in $\Bbb C.$
\endproclaim
\vskip 10pt
{\bf Proof:}
First, write
$$
\multline
L_M(s,F,\chi_{d_0}\chi_{a_1l_1})=
L(s,F,\chi_{a_1d_0l_1}) \cdot \prod_{p|(M/l_1)}
\big(1-\chi_{a_1d_0l_1}(p)p^{-s}\big)^3 \\ =
L(s,F,\chi_{a_1d_0l_1}) \cdot
\Bigg( \sum_{l|(M/l_1)}\mu(l)
\chi_{a_1d_0l_1}(l)l^{-s} \Bigg)^3.
\endmultline
\tag 4.19
$$
By $(4.4)$
$$
L(s,F,\chi_{a_1d_0l_1}) =
(d_0l_1D_{a_1d_0l_1})^{3/2-3s}\frac{G_{\epsilon}(1-s) ^{3}}
{G_{\epsilon}(s) ^{3}}L(1-s,F,\chi_{a_1d_0l_1}),
\tag 4.20
$$
where $G_{\epsilon}$ and $D_{a_1d_0l_1}$ is given by $(4.17)$ and
$\epsilon$ equals the sign of $a_1d_0l_1,$.

On the other side of the functional equation  (4.20), we have,
$$
L(1-s,F,\chi_{a_1d_0l_1}) = L_M(1-s,F,\chi_{a_1d_0l_1})\cdot
\prod_{p|(M/l_1)}
\big(1-\chi_{a_1d_0l_1}(p)p^{-1+s}\big)^{-3}.
$$
In view of the elementary identity
$$
\prod_{p|(M/l_1)}\big(1-p^{-2+2s}\big) = A_2(1-s, \chi_{a_1d_0l_1})
\underset \; p\ne 2\to{\prod_{p|(M/l_1)}}
\big(1+\chi_{a_1d_0l_1}(p)p^{-1+s}\big) \prod_{p|(M/l_1)}
\big(1-\chi_{a_1d_0l_1}(p)p^{-1+s}\big)
$$
where
$$
A_2(s, \chi_{a_1d_0l_1})=\cases
1&\text{if\; \; $2 | l_1,$}\\
1 + \chi_{a_1d_0l_1}(2)2 ^ {-s}&\text{if\; \; $a_1d_0l_1\equiv 1
\pmod{4},$}\\
1 - 2 ^ {-2s}&\text{if\; \; $a_1d_0l_1\equiv -1 \pmod{4},$}\endcases\tag
4.21
$$
it immediately follows that
$$
\align
L(1-s,F, & \,\chi_{a_1d_0l_1})\cdot\prod_{p|(M/l_1)}
\big(1-p^{-2+2s}\big)^3
    = \\
& \\
&=
 L_M(1-s,F,\chi_{a_1d_0l_1})\cdot A_2(1-s, \chi_{a_1d_0l_1})
^3
\cdot{\underset \; p\ne 2\to{\prod_{p|(M/l_1)}}}
\big(1+\chi_{a_1d_0l_1}(p)p^{-1+s}\big)^3\\
&\phantom{xxxxxxxx}= L_M(1-s,F,\chi_{a_1d_0l_1})\cdot A_2(1-s,
\chi_{a_1d_0l_1}) ^3\cdot
\left({\underset \; (l, 2)
= 1\to{\sum_{l|(M/l_1)}}}
\chi_{a_1d_0l_1}(l)l^{-1+s}\right)^3.
\endalign
$$

Combining the above with $(4.7),$ $(4.13),$ $(4.20),$ we obtain
$$
\multline
Z_M(s,w;\chi_{a_2l_2},\chi_{a_1l_1})\cdot\prod_{p|(M/l_1)}
\big(1-p^{-2+2s}\big)^3
= \\ \sum_{(d,M)=1}
(l_1D_{a_1d_0l_1})^{3/2-3s}\frac{G_{\epsilon(a_{1})}(1-s) ^{3}}
{G_{\epsilon(a_{1})}(s) ^{3}}\cdot \frac{L_M(1-s, F,
\chi_{a_1d_0l_1})}{d ^{w+ 3s -3/2}}\cdot
P_{d_0,d_1}^{(a_1l_1)}(1-s)\chi_{a_2l_2}(d_0)\\ \cdot
\Bigg(\sum_{l|(M/l_1)}\mu(l)
\chi_{a_1d_0l_1}(l)l^{-s}\Bigg)^3
\cdot A_2(1-s, \chi_{a_1d_0l_1}) ^3\cdot \Bigg({\underset \; (l, 2) =
1\to{\sum_{l|(M/l_1)}}}\chi_{a_1d_0l_1}(l)l^{-1+s}\Bigg)^3.
\endmultline
\tag 4.22
$$
Write
$$
\multline
\Bigg(\sum_{l|(M/l_1)}\mu(l)
\chi_{a_1d_0l_1}(l)l^{-s}\Bigg)^3=\\ \sum_{l_\a|(M/l_1)}\mu(l_\a)
\chi_{a_1d_0l_1}(l_\a)l_\a^{-s} \cdot \sum_{l_\b|(M/l_1)}\mu(l_\b)
\chi_{a_1d_0l_1}(l_\b)l_\b^{-s} \cdot\sum_{l_\g|(M/l_1)}\mu(l_\g)
\chi_{a_1d_0l_1}(l_\g)l_\g^{-s},
\endmultline
$$
and similarly, write
$$
\multline
\Bigg({\underset \; (l, 2) = 1\to{\sum_{l|(M/l_{1})}}}
 \chi_{a_1d_0l_1}(l)l^{-1+s}\Bigg)^3= \\
{\underset \; (l_{\tilde\a}, 2) =
1\to{\sum_{l_{\tilde\a}|(M/l_1)}}}\chi_{a_1d_0l_1}(l_{\tilde\a})l_{\tilde\a}^{-1
+ s}
\cdot {\underset \; (l_{\tilde\b}, 2) = 1\to{\sum_{l_{\tilde\b}|(M/l_1)}}}
\chi_{a_1d_0l_1}(l_{\tilde\b})l_{\tilde\b}^{-1 + s}
\cdot {\underset \; (l_{\tilde\b}, 2) =
1\to{\sum_{l_{\tilde\g}|(M/l_1)}}}\chi_{a_1d_0l_1}(l_{\tilde\g})l_{\tilde\g}^{-1
+ s}.
\endmultline
$$

It is quite clear that $(4.22)$ decomposes into a linear combination of
the functions $$Z_M(1-s,w+3s-3/2;\chi^{(*)},\chi_{a_1l_1})$$ depending upon
the
congruence class of $a_1l_1$ modulo $4.$ Since the shape of the final
result
is very similar in all the three cases (as in the previous proposition,
just the behavior at the finite place $2$ changes), we will just consider
the case of $a_1l_1 \equiv -1 \pmod 4,$ say. The character $\chi^*$ takes
one of the two forms $\chi_{l_\a l_\b l_\g
l_{\tilde\a}l_{\tilde\b}l_{\tilde\g}}\chi_{a_2l_2},$ $\chi_{-1}\chi_{l_\a
l_\b l_\g
l_{\tilde\a}l_{\tilde\b}l_{\tilde\g}}\chi_{a_2l_2}.$ Note that for
$d_0\equiv 1 \pmod 4,$ $\chi_{a_1d_0l_1}(2) = 0$ and
$\chi_{a_1d_0l_1}(l') = \chi_{a_1l_1}(l')\chi_{d_0}(l')=
\chi_{a_1l_1}(l')\chi_{l'}(d_0),$ for $(l',2)=1.$ For $d_0\equiv -1 \pmod
4$ and any $l>0,$
$\chi_{a_1d_0l_1}(l) = \chi_{l}(a_1l_1)\chi_{l}(d_0).$ Using this and the
character $\chi_{-1}$ to separate the congruence classes $1,-1 \pmod 4,$
we
combine $(4.22)$ with the definition of $Z_M$ in $(4.13)$ to obtain
$$
\multline
Z_M(s,w;\chi_{a_2l_2},\chi_{a_1l_1})\cdot\prod_{p|(M/l_1)}
\big(1-p^{-2+2s}\big)^3
= l_1^{3/2-3s}\frac{G_{\epsilon(a_{1})}(1-s) ^{3}}
{G_{\epsilon(a_{1})}(s) ^{3}} \cdot \frac{1}{2} \Big[4 ^{3/2-3s}(1 - 2 ^
{-2 + 2s}) ^3 \\
\cdot\sum_{l_{\a}|(M/l_1), \; (2, l_{\a})=1}\mu(l_{\a})
\chi_{a_1l_1}(l_{\a})l_{\a}^{-s}
\cdot\sum_{l_{\b}|(M/l_1), \; (2, l_{\b})=1}\mu(l_{\b})
\chi_{a_1l_1}(l_{\b})l_{\b}^{-s}\\
\cdot \sum_{l_{\g}|(M/l_1), \; (2, l_{\g})=1}\mu(l_{\g})
\chi_{a_1l_1}(l_{\g})l_{\g}^{-s}
\cdot \sum_{l_{\tilde\a}|(M/l_1), \; (2, l_{\tilde\a})=1}
\chi_{a_1l_1}(l_{\tilde\a})l_{\tilde\a}^{-1+s}
\\ \cdot\sum_{l_{\tilde\b}|(M/l_1), \; (2, l_{\tilde\b})=1}
\chi_{a_1l_1}(l_{\tilde\b})l_{\tilde\b}^{-1+s}
\cdot \sum_{l_{\tilde\g}|(M/l_1), \; (2, l_{\tilde\g})=1}
\chi_{a_1l_1}(l_{\tilde\g})l_{\tilde\g}^{-1+s}\\
\times \Big(Z_M(1-s,w+3s-3/2;
\chi_{l_\a l_\b l_\g
l_{\tilde\a}l_{\tilde\b}l_{\tilde\g}}\chi_{a_2l_2},\chi_{a_1l_1})
+ Z_M(1-s,w+3s-3/2;
\chi_{-1}\chi_{l_\a l_\b l_\g
l_{\tilde\a}l_{\tilde\b}l_{\tilde\g}}\chi_{a_2l_2},\chi_{a_1l_1})
\Big)\\ +
\sum_{l_{\a}|(M/l_1)}\mu(l_{\a})
\chi_{l_{\a}}(a_1l_1)l_{\a}^{-s}
\cdot\sum_{l_{\b}|(M/l_1)}\mu(l_{\b})
\chi_{l_{\b}}(a_1l_1)l_{\b}^{-s}
\cdot \sum_{l_{\g}|(M/l_1)}\mu(l_{\g})
\chi_{l_{\g}}(a_1l_1)l_{\g}^{-s}
\\ \cdot \sum_{l_{\tilde\a}|(M/l_1)}
\chi_{l_{\tilde\a}}(a_1l_1)l_{\tilde\a}^{-1+s}
\cdot\sum_{l_{\tilde\b}|(M/l_1)}
\chi_{l_{\tilde\b}}(a_1l_1)l_{\tilde\b}^{-1+s}
\cdot \sum_{l_{\tilde\g}|(M/l_1)}
\chi_{l_{\tilde\g}}(a_1l_1)l_{\tilde\g}^{-1+s}
\\ \times
\Big(Z_M(1-s,w+3s-3/2;
\chi_{l_\a l_\b l_\g
l_{\tilde\a}l_{\tilde\b}l_{\tilde\g}}\chi_{a_2l_2},\chi_{a_1l_1})
- Z_M(1-s, w+3s-3/2;
\chi_{-1}\chi_{l_\a l_\b l_\g
l_{\tilde\a}l_{\tilde\b}l_{\tilde\g}}\chi_{a_2l_2},\chi_{a_1l_1})
\big) \Big].
\endmultline
\tag 4.23
$$
This rather complicated formula is the content of Proposition 4.3, where
it is
expressed in a considerably more compact way.

This completes the proof of Proposition 4.3.

\vskip 20pt

\noindent{\bf \S 4.3 The analytic continuation of
$Z_M(s, w; \chi_{a_2l_2}, \chi_{a_1l_1})$}
\vskip 10pt

We begin by recalling some fundamental concepts from the theory of several
complex variables. Our basic reference is H\"ormander \cite{H\"o}.

\proclaim{Definition 4.4} An open set $R$ in $\Bbb C^{m}$ is called a
{\bf domain of holomorphy} if there are no open sets $R_1$ and $R_2$ in
$\Bbb C^{m}$ such that $\varnothing \ne R_1\subset R_2\cap R,$ $R_2$
is connected and not contained in $R,$ and for any holomorphic function
$f$ in $R$ there exists a holomorphic function $f_2$ in $R_2$ satisfying
$f = f_2$ in $R_1.$\endproclaim

\vskip 10pt
\proclaim{Definition 4.5} An open set $\Omega$ in $\Bbb C^{m}$ is called a
{\bf tube} if there is an open set $\omega$ in $\Bbb R^{m},$ called
{\bf the base} of $\Omega,$ such that
$\Omega = \{s \; | \; \Re(s) \in \omega\}.$
\endproclaim

We will denote by $\hat {R},$ the convex hull of a subset
$R \subset \Bbb R^{m}\; \text{or} \; \Bbb C^{m}.$ It is easy to see that
the
convex hull $\hat {\Omega}$ of a tube $\Omega$ is a tube with base $\hat
{\omega}.$
\vskip 10pt
\proclaim{Proposition 4.6} If $\Omega$ is a connected tube, then any
holomorphic function in $\Omega$ can be extended to a holomorphic function
$\hat f$ in $\hat{\Omega}.$\endproclaim
\vskip 10pt
\proclaim{Proposition 4.7} Let $R$ and $R'$ be domains of holomorphy in
$\Bbb C^{m}$ and $\Bbb C^{n},$ respectively, and let $f$ be an analytic
map of
$R$ into $\Bbb C^{n}.$ Then the set
$$R_f = \{s\in R \; | \; f(s)\in R'\}$$
is a domain of holomorphy.\endproclaim

  In order to analytically continue $Z_M(s, w; \chi_{a_2l_2}, \chi_{a_1l_1})$ as a
function of two complex variables $s, w$, we repeatedly apply the functional
equations given in Propositions $4.2, 4.3$.

Accordingly, we define two involutions on $\C \times \C:$
$$
\a:(s, w) \rightarrow (1-s, w+3s-3/2) \quad \hbox{and} \quad \b:(s,w)
\rightarrow (s+w-1/2, 1-w).
$$
Then $\a,$ $\b$ generate $D_{12},$ the dihedral group of order 12, and
$\a^2 = \b^2 = 1,$ $(\a\b)^6 = (\b \a)^6 = 1.$ Note that $\a \b \ne \b
\a.$

We will find it useful in the following to define three regions $R_1,$
$R_2,$
$R_3$ as follows: Write $s,$ $w$ as $s = \s + it,$ $w =  \nu + i\g.$

\pagebreak

The tube region $R_1$
\vskip -18pt
\centerline{\Insert{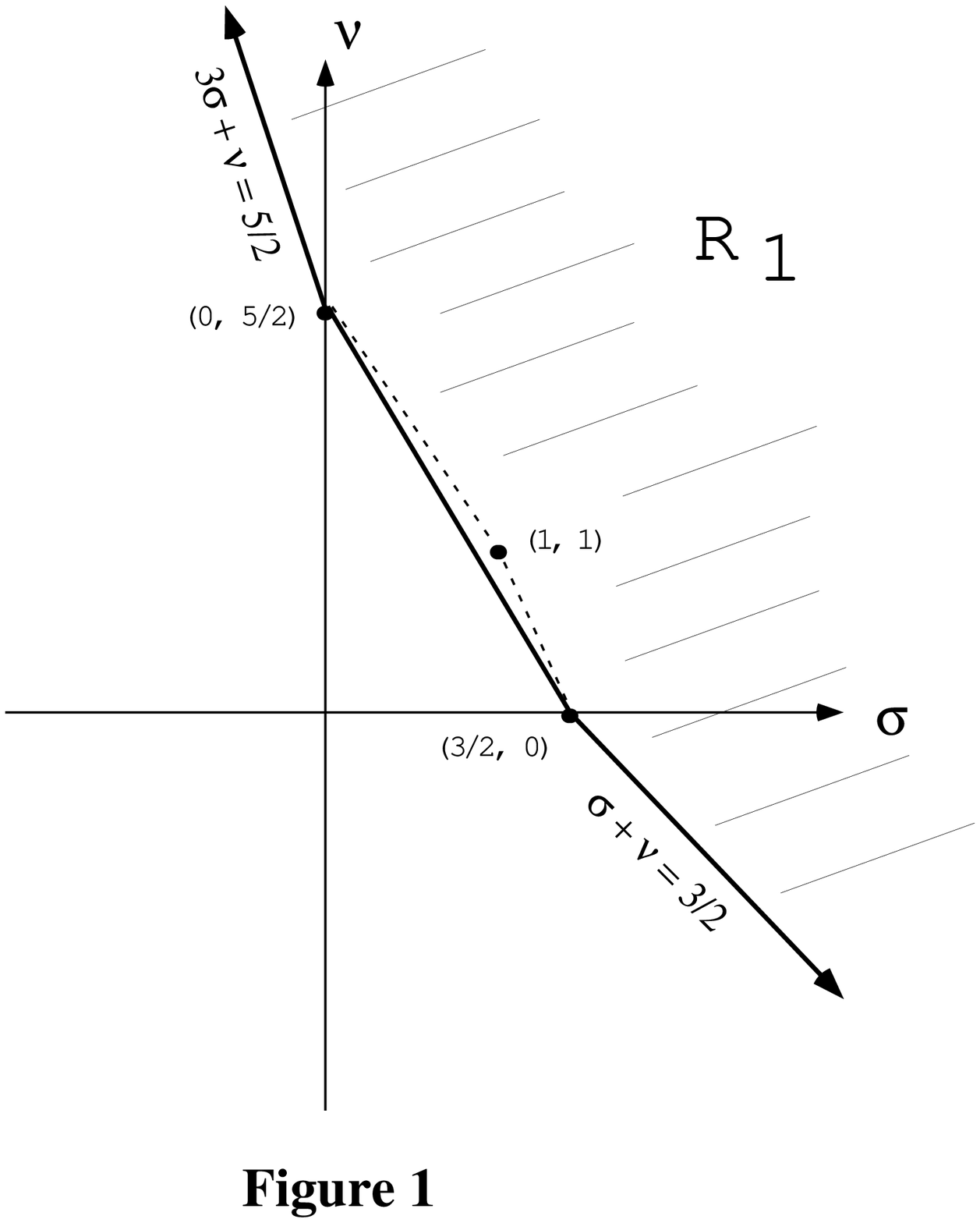}{3.75in}}
\vskip -35pt
\noindent
is defined to be the set of all points $(s, w)$ such
that
$(\s, \nu)$ lie strictly above the polygon determined by
 $(0, 5/2),$ $(3/2, 0),$  and the rays $\nu = -3\s
+
5/2$ for $\s \le 0$ and $\nu = -\s + 3/2$ for $\s \ge 3/2.$
Note that $R_1$ is the convex closure of the region given in Figure 1
which is bounded by the dotted lines and the two rays $\nu = -3\s
+
5/2$ for $\s \le 0$ and $\nu = -\s + 3/2$
for $\s \ge 3/2,$ which is the actual region that comes up in the proof of
Propositions $4.8, 4.9.$

The tube region $R_2$
\vskip -12pt
\centerline{\Insert{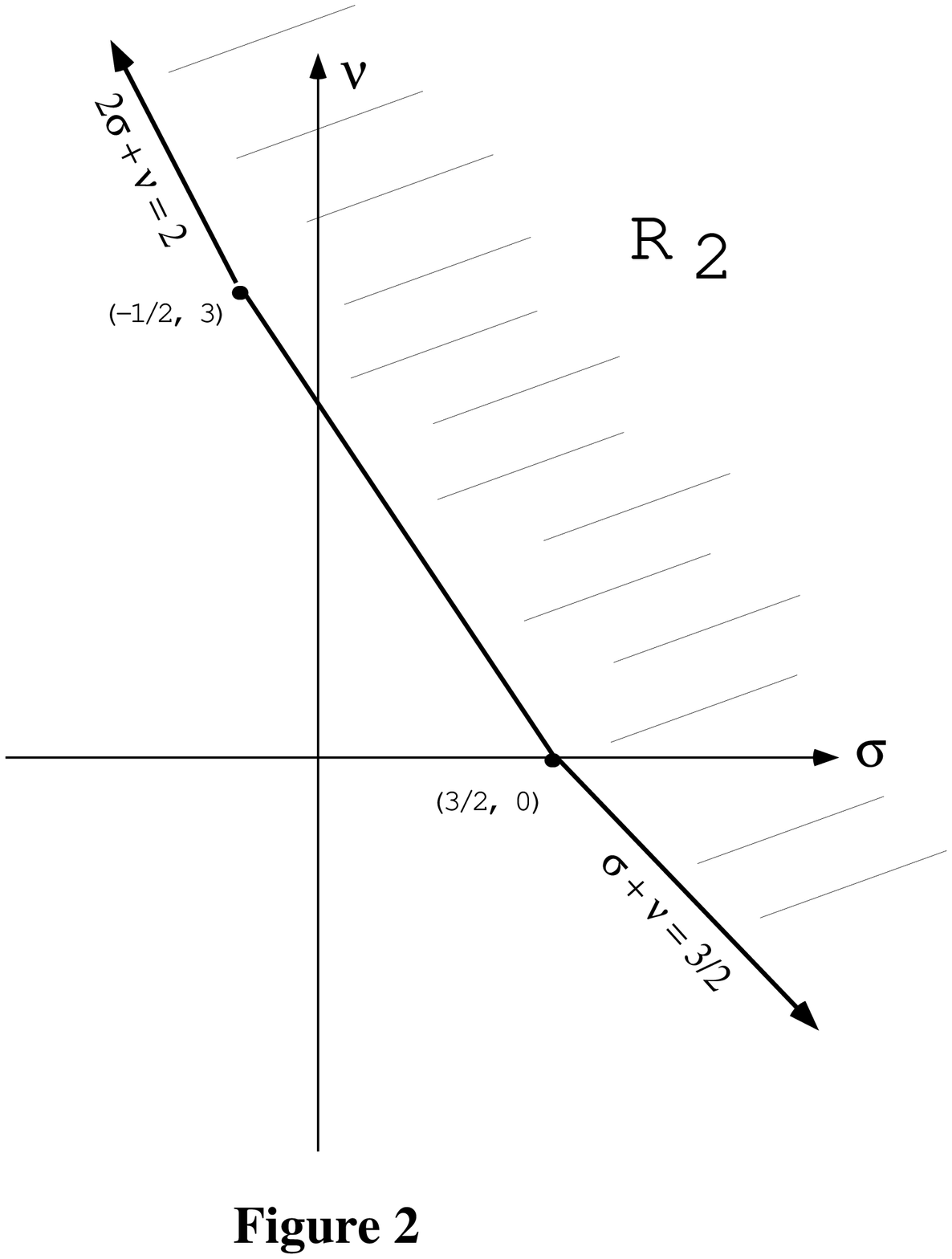}{3.75in}}
\vskip -24pt

\pagebreak
\noindent
is defined to be the set of all points $(s, w)$ such
that
$(\s, \nu)$ lie strictly above the line segment connecting $(-1/2, 3)$
and $(3/2, 0)$ and the rays $\nu = -2\s + 2$ for $\s \le -1/2,$ and
$\nu = -\s +3/2$ for $\s \ge 3/2.$

The tube region $R_3$
\vskip -12pt
\centerline{\Insert{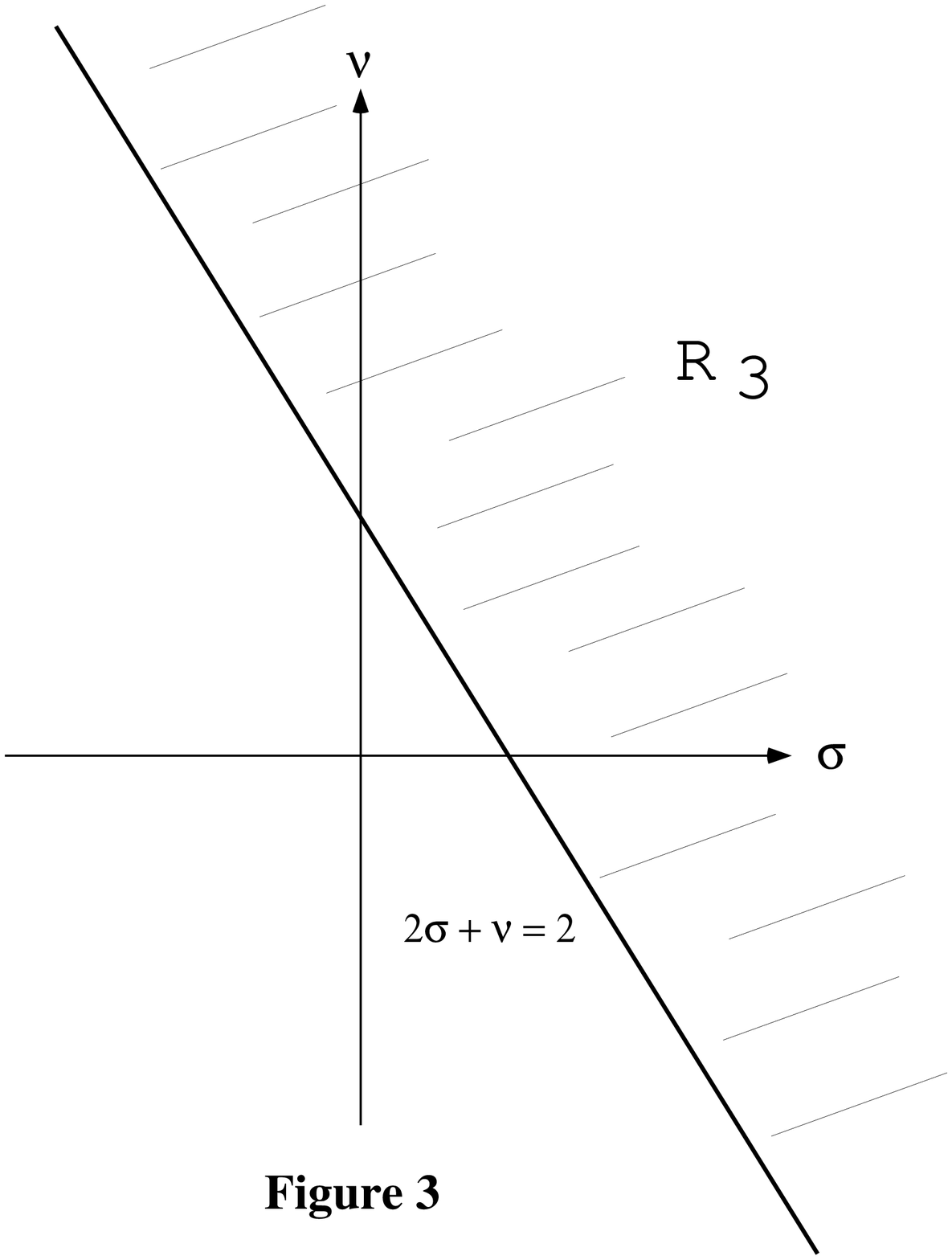}{3.75in}}
\vskip -24pt
\noindent
 is defined to be the set of all points $(s, w)$ such
that
$(\s, \nu)$ lie strictly above the line $\nu = -2\s + 2.$

These regions are related by the involutions $\a,$ $\b$ as described in
the following proposition. The proof, a simple exercise, is omitted.
\vskip 10pt
\proclaim{Proposition 4.8}
The regions $R_1$ and $\a (R_1)$ have a non-empty intersection, and
the convex hull of $R_1 \cup \,\a (R_1)$ equals $R_2.$ Similarly, $R_2$
and
$\b (R_2)$ have a non-empty intersection and the convex hull of $R_2 \cup
\,\b (R_2)$ equals $R_3.$ Finally, $R_3$ and $\a (R_3)$ have a non-empty
intersection and the convex hull of $R_3 \cup \,\a (R_3)$ equals $\C^2.$
\endproclaim
\vskip 10pt
Let
$$
P(s, w) = (s - 1)^3 (w - 1).
\tag 4.24
$$
We will begin by demonstrating
\vskip 10pt
\proclaim{Proposition 4.9}
Let $R_1$ be the tube region defined above. The function
$$
P(s, w)Z_M(s, w; \chi_{a_2l_2}, \chi_{a_1l_1})
$$
is analytic in $R_1.$
\endproclaim
\vskip 10pt
{\bf Proof:}
Consider first the left hand side of the expression for
$Z_M(s,w;\chi_{a_2l_2},\chi_{a_1l_1})$ given in $(4.10)$. If
the sum were restricted only to square free $d = d_0,$ then the usual
Phragmen-Lindel\"of bounds for $L(s,\chi_{d_0})$ would imply absolute
convergence for $\nu >1$ when $\s >1,$ for $\nu > (-3/2) \s + 5/2$ when
$0\le\s\le1$ and for $\nu > -3 \s + 5/2$ when $\s < 0.$ Because we have
the
bound $(4.11)$ and functional equation $(4.7)$ applied to
$P_{d_0,d_1}^{(a_1l_1)}(s),$ precisely the same estimates apply as we sum
over all $d.$ Consequently, $Z_M(s, w; \chi_{a_2l_2}, \chi_{a_1l_1})$
converges above the given lines, and the factor $(s - 1)^3$ in $P(s,w)$
cancels the pole at $s = 1.$

Noting that both sides of the expression converge when $\nu,$ $\s >1,$ we
now change the order of summation and examine the right hand side. Here
the coefficients $c(n)$ are order 3 divisor functions and are bounded
above
by $n^\epsilon$ for any $\epsilon >0.$ Consequently, applying
Phragmen-Lindel\"of again to $L(w, \chi_{n_0})$ and the corresponding
estimate and functional equations for
$c(n_0n_1^2)Q_{n_0,n_1}^{(a_2l_2)}(w),$
we obtain convergence of $Z_M(s, w; \chi_{a_2l_2}, \chi_{a_1l_1})$ for
$\s >1$ when $\nu >1,$ for $\s > (-1/2) \nu + 3/2$ when $0 \le \nu \le 1$
and $\s > - \nu + 3/2$ when $\nu < 0.$ The factor $w - 1$ in $P(s, w)$
cancels the pole at $w=1.$ These regions overlap when $\nu,$ $\s >1,$ and
thus by Proposition 4.6, $Z_M(s, w; \chi_{a_2l_2}, \chi_{a_1l_1})P(s,w)$
has an analytic continuation to the convex closure of the regions, which
is $R_1$ described above.

This completes the proof of Proposition 4.14.

Our plan is now to apply the involutions $\a,$ $\b,$ $\a$ in that order to
$R_1,$ and use Propositions 4.2 and 4.3 to extend the analytic
continuation
to $\Bbb C ^2.$
To aid in this, it will be useful to introduce some additional notation
to make the content of these propositions a bit clearer and easier to
apply.
Let

$$
A(s, w) \equiv A_M(s, w) =\prod_{p|M}(1 - p ^{-2 + 2s})^3 \quad \hbox{and}
\quad B(s, w) \equiv B_M(s, w) = \prod_{p|M}(1 - p ^{-2 + 2w}),
\tag 4.25
$$
and let $\tilde \Psi ^{(a_1l_1)} (s, w) = \Psi^ {(a_1l_1)}
(s)\prod_{p|l_1}(1 - p ^{-2 + 2s})^3,$
$\tilde\Phi ^{(a_2l_2)} (s, w) = \Phi^{(a_2l_2)} (w)
\prod_{p|l_2}(1 - p ^{-2 + 2w}).$

The following is a reformulation of the content we require now from
Propositions 4.2 and 4.3. For $(s, w)$ such that both sides are contained
in a connected region of analytic continuation for
$P(s, w)Z_M(s, w; \chi_{a_2l_2}, \chi_{a_1l_1})$
$$
A(s, w)\overarrow\Z_M(s, w; \chi_{\text{Div}(M)}, \chi_{a_1l_1})
=\tilde\Psi^{(a_1l_1)} (s, w)
\overarrow\Z_M(\a(s, w); \chi_{\text{Div}(M)}, \chi_{a_1l_1})
\tag 4.26
$$
and
$$
B(s, w)\overarrow\Z_M(s, w; \chi_{a_2l_2}, \chi_{\text{Div}(M)})
=\tilde\Phi ^{(a_2l_2)} (s, w)
\overarrow\Z_M(\b(s, w); \chi_{a_2l_2}, \chi_{\text{Div}(M)}).
\tag 4.27
$$

The following proposition will now complete the analytic continuation of
$Z_M(s, w; \chi_{a_2l_2}, \chi_{a_1l_1}).$

\vskip 10pt
\proclaim{Proposition 4.10}
Let
$$
\multline
\Pcal(s,w) =  s^3(s-1)^3(s+w-3/2)^3(2s+w-1)^3(s+w-1/2)^3(2s+w-2)^3\\
\times w(w-1)(3s+w-5/2)(3s+2w-3)(3s+w-3/2).
\endmultline
$$
Then the following product has an analytic continuation to an entire
function in $\C^2:$
$$
\multline
{\widetilde Z}_M(s, w; \chi_{a_2l_2}, \chi_{a_1l_1}) :=
A(s,w)A(\a(s,w))A(\b(s,w))
A(\b\a(s,w))B(s,w)B(\a(s,w))\Pcal(s,w)\\
\times Z_M(s, w; \chi_{a_2l_2}, \chi_{a_1l_1}).
\endmultline
$$
\endproclaim
\vskip 10pt
{\bf Proof:}
In Proposition 4.9 we established the continuation of
$Z_M(s, w; \chi_{a_2l_2}, \chi_{a_1l_1})P(s,w)$ in $R_1.$ As $\a ^2 = 1$
and
$\tilde\Psi^{(a_1l_1)} (s, w)$ is
meromorphic in $\C ^2,$ it follows that
$$\tilde\Psi^{(a_1l_1)}(s, w)\overarrow\Z_M(\a(s, w);
\chi_{\text{Div}(M)},\chi_{a_1l_1})P(\a(s, w))$$ is a meromorphic function
in $\a  (R_1).$
 From $(4.23),$ we observe that poles can just occur at the points $s  = 1, 3,
5, \ldots$ or $s = 2, 4, 6, \ldots$ (depending on $\epsilon(a_1)$).
However,
except for the possible pole at $s = 1,$ all the others are canceled by the
trivial zeros of $L(1 - s,\chi_{d_0}).$ We can conclude from Proposition
4.9 and $(4.26)$ that
$A(s, w)P(s, w)P(\a(s, w))\overarrow\Z_M(s,
w;\chi_{\text{Div}(M)},\chi_{a_1l_1})$ is analytic in $R_1 \cup \;
\a(R_1),$ $R_1$ and $\a  (R_1)$ having a
substantial intersection (containing $\Re (s),$ $\Re (w) >1$). Thus by
Proposition 4.6, this function is analytic in $R_2,$ the convex hull of
the
union.

Since $\b ^2 = 1$ and $\tilde\Phi ^{(a_2l_2)}(s, w)$ is meromorphic in
$\C ^2,$ it follows from what we have just proved that
$$
\tilde\Phi ^{(a_2l_2)} (s, w)A(\b (s, w))P(\b (s, w))P(\a\b(s, w))
\overarrow\Z_M(\b(s, w); \chi_{a_2l_2}, \chi_{\text{Div}(M)})
$$
is a meromorphic function in $\b(R_2).$ As before, all the poles,
except the possible one at $w = 1,$ of $\tilde\Phi ^{(a_2l_2)} (s, w)$ are
canceled by trivial zeros of $L$--functions. From $(4.27),$ we conclude
that
$$
A(s,w)A(\b(s,w))B(s,w)P(s,w)P(\a(s,w))P(\b(s,w))P(\a\b(s,w))
\overarrow\Z_M(s, w; \chi_{a_2l_2}, \chi_{\text{Div}(M)})
\tag 4.28
$$
is an analytic function in $R_2 \cup \b(R_2).$ As this has a non-empty
intersection, it follows from Proposition 4.6 again that $(4.28)$ is
analytic
in $R_3,$ the convex hull of $R_2 \cup \b(R_2).$

To complete the argument, apply $\a$ to $(4.26),$ obtaining
$$
A(\a(s, w))\overarrow\Z_M(\a(s, w); \chi_{\text{Div}(M)}, \chi_{a_1l_1})
=\tilde\Psi^{(a_1l_1)} (\a(s, w))
\overarrow\Z_M(s, w; \chi_{\text{Div}(M)}, \chi_{a_1l_1}).
$$
Multiplying the above by
$A(s, w)A(\b(s, w))B(s, w)P(s, w)P(\a(s, w))P(\b(s, w))P(\a\b(s, w))$
and applying $(4.28),$ we see that
$$
A(s, w)A(\a(s, w))A(\b(s, w))B(s, w)P(s, w)P(\a(s, w))P(\b(s, w))P(\a\b(s,
w))
\overarrow\Z_M(\a(s, w); \chi_{\text{Div}(M)}, \chi_{a_1l_1})
$$
is analytic for $(s, w)\in R_3.$ Replacing $(s, w)$ by $\a(s, w),$ we
obtain
$$
\multline
A(s, w)A(\a(s, w))
A(\b\a(s, w))B(\a(s, w))P(s, w)P(\a(s, w))P(\b\a(s, w))P(\a\b\a(s, w))\\
\times \overarrow\Z_M(s, w; \chi_{\text{Div}(M)}, \chi_{a_1l_1})
\endmultline
$$
is analytic for $(s, w)\in \a(R_3).$ Combining this with the fact that
$(4.28)$ is analytic in $R_3,$ we obtain the analyticity of
$$
\multline
A(s, w)A(\a(s, w))A(\b(s, w))A(\b\a(s, w))B(s, w)
B(\a(s, w))P(s, w)P(\a(s, w))P(\b(s, w))\\ \times
P(\b\a(s, w))P(\a\b(s, w))P(\a\b\a(s, w))
\overarrow\Z_M(s, w; \chi_{\text{Div}(M)}, \chi_{a_1l_1})
\endmultline
$$
in $R_3 \cup \a(R_3).$ As this has a non-empty intersection, it follows
from
Proposition 4.6 again that the above is analytic in $\C^2,$ the convex
hull
of $R_3 \cup \b(R_3).$

In fact, $P(\a\b(s, w)),$ $P(\a\b\a(s, w))$ have one factor in common:
$2w + 3s - 3,$ and so in the last step we included one unnecessary
multiple of
$2w + 3s - 3.$ Removing this, we complete the proof of Proposition 4.10.

\vskip 20pt\noindent
\noindent{\bf \S 4.4 An estimate for
$Z_M\left(\frac12, w; \chi_{a_2l_2}, \chi_{a_1l_1}\right)$ in vertical strips.}
\vskip 10pt
In this section we will use the analytic continuation and functional
equations $(4.26),$ $(4.27)$ for
$\overarrow\Z_M(s, w; \chi_{\text{Div}(M)}, \chi_{a_1l_1})$ to locate
poles
and obtain an estimate for the growth of this function in a vertical strip. Before
doing this, however,
we need some additional notation.

Let $\overarrow\Z_M(s,w)$ denote the $4d(M)^2$--dimensional column vector
consisting of the
concatenation of the $2d(M)$ column vectors
$\overarrow\Z_M(s,w;\chi_{a_2l_2},\chi_{\text{Div}(M)})$ for $a_2 \in
\{1,-1\}$
and all $l_2|M.$
Then by Propositions 4.2 and 4.3, combined with $(4.26),$ $(4.27),$ there
exist $4d(M)^2$ by $4d(M)^2$ matrices
$\Phi_M(s, w), \Psi_M(s, w)$ such that
$$
A_M(s, w) \overarrow\Z_M(s,w) = \Psi_M(s,w) \overarrow\Z_M(\a(s,w))
\tag 4.29
$$
and
$$
B_M(s,w) \overarrow\Z_M(s,w) = \Phi_M(s,w) \overarrow\Z_M(\b(s,w)).
\tag 4.30
$$
Here $A_M(s,w),$ $B_M(s,w)$ are given by $(4.25).$ The matrices
$\Phi_M(s, w),$ $\Psi_M(s, w)$ are constructed from blocks of
$\tilde\Phi^{(a_2l_2)} (s, w)$ and
$\tilde\Psi^{(a_1l_1)} (s, w)$ on the diagonal.

Next, we use Proposition 4.7 to show that the function
${\widetilde Z}_M(1/2, w; \chi_{a_2l_2}, \chi_{a_1l_1}),$
defined in Proposition 4.10, is of finite order. Although it seems to be a
one--variable problem, the theory of several complex variables is still
needed in the proof.
\vskip 10pt
\proclaim{Proposition 4.11} The entire function
$$
{\widetilde Z}_M\left(\frac{1}{2}, w; \chi_{a_2l_2}, \chi_{a_1l_1}\right)
$$
is of the first order.
\endproclaim
\vskip 10pt
{\bf Proof:} First, the convexity bound $L(1/2, \chi_{d_0}) \ll_{\epsilon}
d_0 ^{\frac{1}{4} + \epsilon}$ together with $(4.11),$ implies that
$$
Z_M\bigg(\frac{1}{2}, w; \chi_{a_2l_2}, \chi_{a_1l_1}\bigg)\ll_{\epsilon}
1,
$$
for $\Re(w) = \nu > \frac{7}{4} + \epsilon.$
Applying $(4.29)$ and $(4.30)$ several times in succession,
we obtain
$$
\multline
\overarrow\Z_M(s, w)  = B_M(s, w)^{-1}\Phi_M(s, w)A_M(\b(s,
w))^{-1}\Psi_M(\b(s, w))
B_M(\a\b(s, w))^{-1}\Phi_M(\a\b(s, w))\\
\times A_M(\b\a\b(s, w))^{-1}\Psi_M(\b\a\b(s, w))B_M((\a\b)^2(s, w))^{-1}
\Phi_M((\a\b)^2(s, w)) \overarrow\Z_M(s, 5/2 - 3s - w).
\endmultline
\tag 4.31$$

For $s = 1/2,$ we observe that $\overarrow\Z_M(1/2, w)$ is related to
$\overarrow\Z_M(1/2, 1 - w)$ by the functional equation $(4.31).$ Using
Stirling's formula, we can bound from above the entries of the right hand side
matrices in
$(4.31),$ obtaining
$$
Z_M\bigg(\frac{1}{2}, \nu + it; \chi_{a_2l_2}, \chi_{a_1l_1}\bigg)
\ll_{\epsilon} (1 + |t|) ^C,
$$
where $C$ is an absolute positive constant and $\nu < -\frac{3}{4} -
\epsilon.$

The proof of Proposition $4.11$ is based on an application of Proposition
4.7 to the function $f:\Bbb C^2 \to \Bbb C,$ defined by
$$
f(s, w) = \Gamma(s + 5)\Gamma(w + 5){\widetilde Z}_M(s, w;
\chi_{a_2l_2},\chi_{a_1l_1}).
$$

\vbox{Now let $\Omega_0$ be the tube region whose
base  is given in Figure 1. This tube
already appeared at the end of the proof of
Proposition 4.9 (its convex hull is $R_1$).
Reflecting several times under $\a,$ $\b,$
$\a,$ $\b\ldots,$ until it stabilizes and then
taking the union, we obtain a tube whose base
is $\Bbb R^{2}$ with a hole in the middle (see Figure 4 below).
\vskip -18pt
\centerline{\Insert{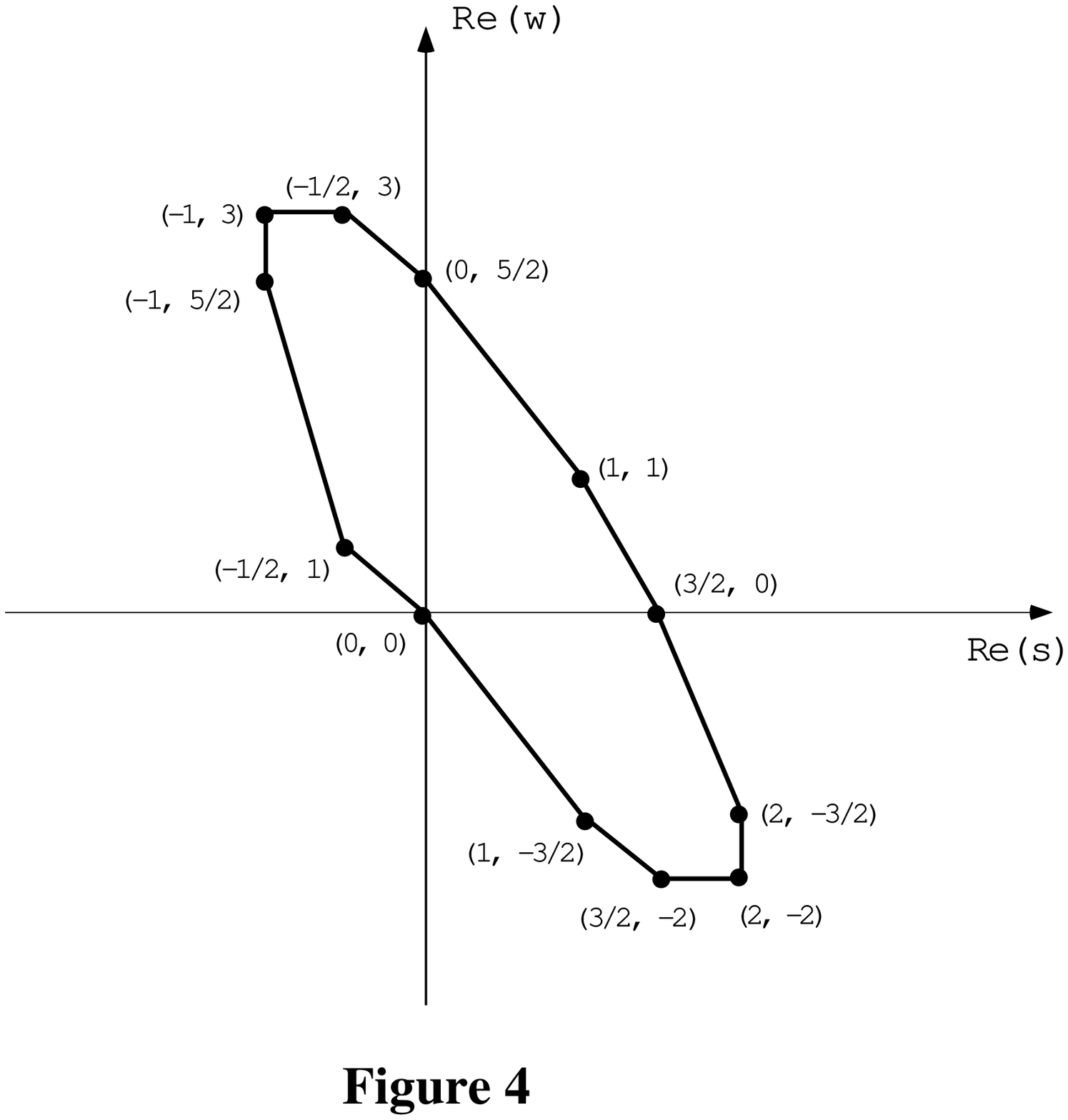}{4in}}
\vskip -48pt
\noindent
This hole is a tube with base a polygon,
which lies inside the
open ball  $B(0, 4)$ (of radius $4$ centered at the origin)  in $\Bbb R^{2}$. The
function ${\widetilde Z}_M (s, w; \chi_{a_2l_2},
\chi_{a_1l_1})$
is obviously of polynomial growth in $\Im(s)$ and $\Im(w)$ as long as
$(s, w) \in \Omega_0,$ and $\s,$ $\nu$ are both bounded. Applying
Stirling's
formula in equations $(4.18)$ and $(4.23),$ we observe that the same holds
when $\a,$ $\b$ are applied. Combining this with Stirling's formula,
we conclude that the function
$f(s, w)$
is bounded in the tube $\Omega '$ with base the annulus $\omega ' = \{(\s,
\nu)
\in \Bbb R^{2}\; |\; 16 < \s ^{2} + \nu ^{2} < 25 \}.$  See Figure 5 below.
\vskip -28pt
\centerline{\Insert{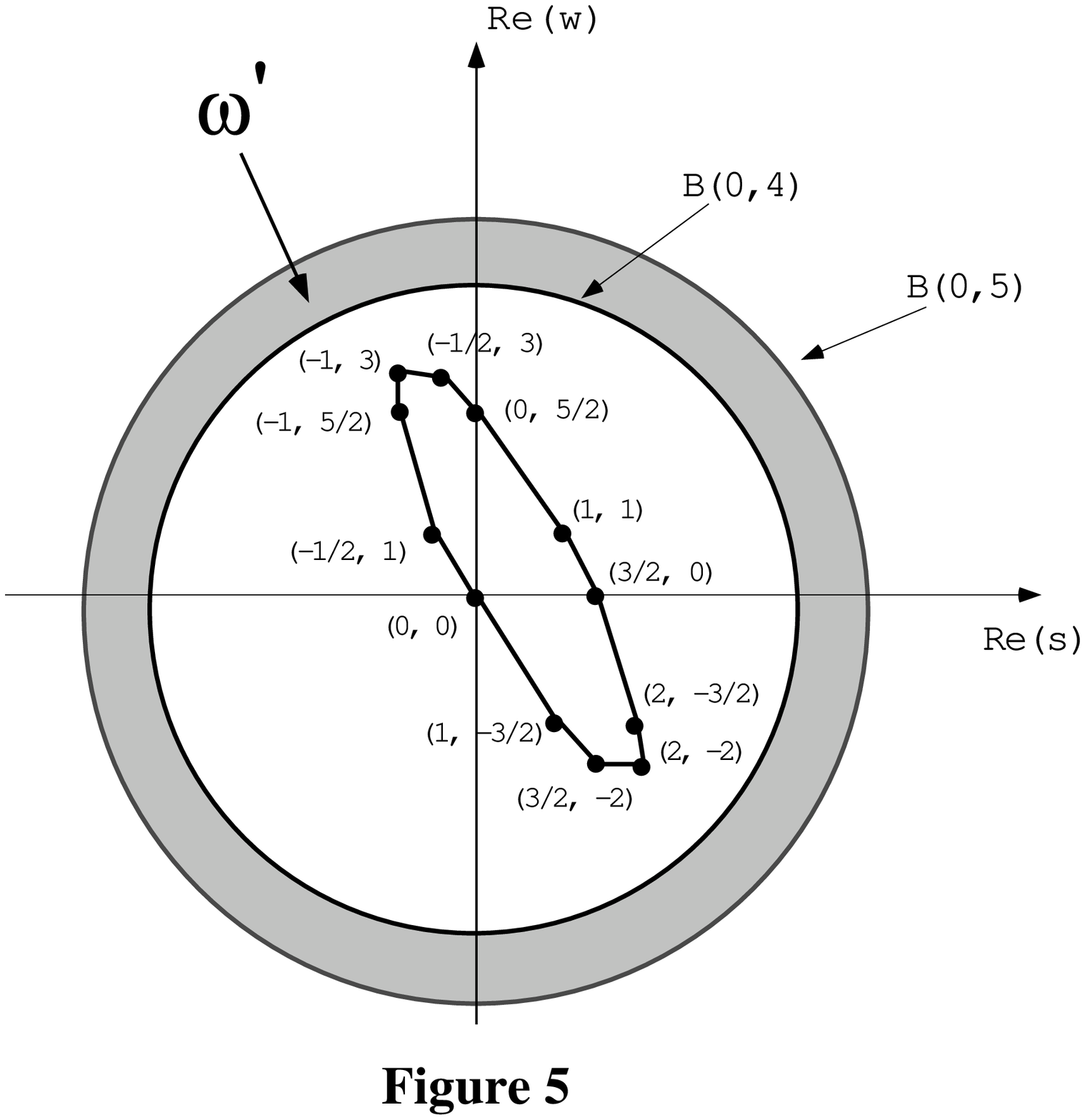}{4in}}
\vskip -224pt
}
\pagebreak

 Let $R \subset \Bbb C^2$
be the
tube whose base is $B(0, 5)$ in $\Bbb R^{2},$ and let $R' = B(0, m) \in
\Bbb C,$ where
$m$ is an upper bound for $f$ on the annulus
$\Omega '.$ Since
$B(0, 5)$ in $\Bbb R^{2}$ is a convex set, it follows that $R$ is a
domain of holomorphy. Obviously, $R'$ is also a domain of holomorphy.
Applying
Proposition 4.7, it follows that
$$
\Gamma(s + 5)\Gamma(w + 5){\widetilde Z}_M(s, w; \chi_{a_2l_2},
\chi_{a_1l_1})
$$
is bounded in $R,$ since in this case, the set $R_f$ contains the annulus $\Omega ' $
whose convex closure contains $R$.
In particular, this function is bounded in the tube with
base given by the polygon in Figure 4.
Proposition 4.11 immediately follows.

One of the key ingredients in what follows, is that the series
$$
\sum_{d_0}L\left| \left(\frac{1}{2} + it, \chi_{d_0} \right) \right| ^4
|d_0|^{-\nu}
\tag 4.32
$$
is convergent, for $\nu = \Re (w) >1.$ Here the summation is over all
positive or negative square free integers. This follows from the work of
Heath--Brown \cite{H--B}. Applying the Cauchy--Schwartz inequality, we deduce
that
$$
\sum_{d}|c_{d}|L\left| \left(\frac{1}{2} + it, \chi_{d} \right) \right| ^3
|d|^{-\nu}
\tag 4.33
$$
is convergent, for $\nu = \Re (w) >1,$ and any sequence $c_{d}$ such that
$c_{d}\ll_\epsilon d ^{\epsilon}.$ Here the summation is over all
integers.

We now show:
\vskip 10pt
\proclaim{Proposition 4.12}
Let $w = \nu + it.$ For $\epsilon > 0,$ $-\epsilon \le \nu ,$ and
any $a_1,$ $a_2 \in \{1, -1\},$ $l_1, l_2|M$ the function
$Z_M(1/2, w ;\chi_{a_2l_2}, \chi_{a_1l_1})$
is an analytic function of $w,$ except for possible poles at $w = \frac{3}{4}$
and $w = 1.$ If $(l_1,\; l_2) = 1$ or $2$ and $|t| > 1,$ then it satisfies the
upper
bounds
$$
Z_M\left(\frac12, \nu + it; \chi_{a_2l_2}, \chi_{a_1l_1}\right)
\ll_\epsilon 1,
$$
for $1+\epsilon < \nu,$ and
$$
Z_M\left(\frac12, \nu + it; \chi_{a_2l_2}, \chi_{a_1l_1}\right)
\ll_\epsilon M ^{3(1 - \nu) + v_1(\epsilon)}|t| ^{5(1 - \nu) +
v_2(\epsilon)}
\sum_{a = 1,\; -1}\;
\sum_{l|M}\;  \sum_{(d_{0}, M) = 1}
\frac{\big|L\big(\frac{1}{2}, \chi_{d_0} \chi_{a l}\big) \big| ^{3}}
{d_{0} ^ {1 + \epsilon}},
$$
for $-\epsilon \le \nu \le 1 + \epsilon.$ The functions $v_1(\epsilon),$
$v_2(\epsilon)$ are some explicitly computable functions satisfying
$$
\lim_{\epsilon\rightarrow 0}v_1(\epsilon) = \lim_{\epsilon\rightarrow 0}
v_2(\epsilon) = 0.
$$
\endproclaim
\vskip 10pt
{\bf Proof:} The first bound in the region $1+\epsilon < \nu$ is immediate by the
remarks concerning $(4.33).$ The bound for $-\epsilon \le \nu \le 1+\epsilon$ is
more difficult to obtain. We shall first obtain a bound for
$Z_M\left(\frac12, \nu + it; \chi_{a_2l_2}, \chi_{a_1l_1}\right)$, (i.e., for $\nu =
-\epsilon$), and then apply a convexity argument to complete the proof for $-\epsilon
< \nu < 1+\epsilon$.

Recall the functional equations
$$\align &\alpha(s, w) = \left(1-s, \; 3s+w-\frac32\right) \qquad(\text{see equation
(4.23)})\\ &\beta(s, w) = \left(s+w-\frac12, \; 1-w\right) \qquad(\text{see equation
(4.18)}).\endalign$$
Fix $(s, w) = (\frac12, -\epsilon + it).$ We then have
$$
\b(s, w)=(-\epsilon + it, 1 + \epsilon - it), \quad
\a\b(s, w)=(1 + \epsilon - it, -1/2- 2\epsilon + 2it),
$$
$$
\b\a\b(s, w)=(-\epsilon + it, 3/2 + 2\epsilon - 2it), \quad
\a\b\a\b(s, w)=(1 + \epsilon - it, -\epsilon + it),
$$
and
$$ \b \a\b\a\b(s, w) = \left(\frac12, 1+\epsilon-it\right).$$
We shall estimate $Z_M\left(\frac12, \nu + it; \chi_{a_2l_2}, \chi_{a_1l_1}\right)$
by alternately applying the functional equations $\b, \a$ as above. Note that each
time we apply
$\b$ the value of $w$ is either $-\epsilon+it$ or $-\frac12-2\epsilon +2it,$
and each time we apply $\a$, the value of $s$ is $-\epsilon+it.$ It is thus sufficient to
obtain upper bounds in only these cases. We proceed to do this.

Now, it immediately follows from $(4.18)$ and Stirling's asymptotic formula for
the Gamma function that away from poles,
$$\multline Z_M\left(s, -\epsilon + it; \chi_{a_2l_2}, \chi_{a_1l_1}\right)
\; \ll_\epsilon  \;     l_2^{\frac12+\epsilon}\sum_{l_3, l_4 |  M/l_2} M^\epsilon
\sum_{a_3 = 1, -1} |t|^{\frac12+\epsilon} \cdot \\  \cdot\Bigg(\Big
|Z_M\big(s-{\scriptstyle \frac12}-\epsilon+it, 1+\epsilon-it;
\chi_{a_2l_2}, \chi_{a_1l_1l_3l_4}\big) \Big |\; + \;  \Big |
Z_M\big(s-{\scriptstyle \frac12}-\epsilon+it, 1+\epsilon-it;
\chi_{a_2l_2}, \chi_{-a_1l_1l_3l_4}\big)\Big|\Bigg).\endmultline$$
Since $M$ is even and squarefree, we also have
$$(l_2, l_1l_3l_4) = 1 \; \text{or} \; 2.$$
The characters $\chi_{a_1l_1l_3l_4}$ and $\chi_{-a_1l_1l_3l_4}$ can be replaced
by $\chi_{a_1d_2}, \; \chi_{-a_1d_2}$ with $d_2$ squarefree.

Similarly, for $w = -\frac12-\epsilon+it$, we have, after replacing $l_2$ by $d_3$
and $l_1$ by $d_2$ that
$$\multline Z_M\left(s,-{\scriptstyle \frac12} -\epsilon + it; \chi_{a_2d_3},
\chi_{a_1d_2}\right)
\; \ll_\epsilon  \;     d_3^{1+\epsilon}\sum_{l_3, l_4 |  M/d_3} M^\epsilon
l_3^\frac12 l_4^{-\frac32}
\sum_{a_3 = 1, -1} |t|^{1+\epsilon}  \\  \cdot\Bigg(\Big
|Z_M\big(s-1-\epsilon+it, {\scriptstyle \frac32}+\epsilon-it;
\chi_{a_2d_3}, \chi_{a_1d_4}\big) \Big |\; + \;  \Big |
Z_M\big(s-1-\epsilon+it, {\scriptstyle \frac32}+\epsilon-it;
\chi_{a_2d_3}, \chi_{-a_1d_4}\big)\Big|\Bigg),\endmultline$$
where we have denoted by $d_4$, the squarefree part of $d_2l_3l_4.$ Note that
$(d_3, d_4) = 1 \; \text{or}\; 2.$

 In a similar manner, we consider $s = -\epsilon+it$ in $(4.23).$ It follows from
Stirling's formula that away from poles,

$$\multline Z_M\left( -\epsilon + it,w; \chi_{a_2l_2}, \chi_{a_1d_2}\right)
\; \ll_\epsilon  \;     |d_2\cdot t|^{\frac32+3\epsilon} \sum_{l_\alpha,
l_\beta,l_\gamma,l_{\bar\alpha},l_{\bar\beta},l_{\bar\gamma} |  (M/l_1)}
M^{3\epsilon }
 \cdot \\  \cdot\Bigg(\Big
|Z_M\big(1+\epsilon-it, w-3\epsilon+3it-{\scriptstyle \frac32};
\chi_{l_\alpha
l_\beta l_\gamma l_{\bar\alpha} l_{\bar\beta} l_{\bar\gamma}}\cdot \chi_{a_2l_2},
\;\chi_{a_1d_2}\big)
\Big |\; + \\   \; +\;  \Big | Z_M\big(1+\epsilon-it, w-3\epsilon+3it-{\scriptstyle
\frac32};
\chi_{-l_\alpha
l_\beta l_\gamma l_{\bar\alpha} l_{\bar\beta} l_{\bar\gamma}}\cdot \chi_{a_2l_2},
\;\chi_{a_1d_2}\big)\Big|\Bigg).\endmultline$$
As before, $(l_\alpha
l_\beta l_\gamma l_{\bar\alpha} l_{\bar\beta} l_{\bar\gamma}, d_2) = 1 \; \text{or}
\; 2.$ We can replace $l_\alpha
l_\beta l_\gamma l_{\bar\alpha} l_{\bar\beta} l_{\bar\gamma}$ by $d_3,$
squarefree. We again obtain that $(d_3, d_2) = 1 \; \text{or}
\; 2.$

It now follows from the previous estimates and remarks that
$$\align Z_M\left(\frac12, -\epsilon + it; \chi_{a_2l_2}, \chi_{a_1l_1}\right)
&\ll_{\epsilon} |t| ^{5 + 10\epsilon}M ^{10\epsilon}
d_{1} ^{\frac{1}{2} + \epsilon}
d_{2} ^{\frac{3}{2} + 3\epsilon}
d_{3} ^{1 + 2\epsilon}l_{3} ^{\frac{1}{2}}l_{4} ^{- \frac{3}{2}}
d_{4} ^{\frac{3}{2} + 3\epsilon}d_{5} ^{\frac{1}{2} + \epsilon} \cdot S\\
&= |t| ^{5 + 10\epsilon}M ^{10\epsilon}(d_{1}d_{2})  ^{\frac{1}{2} +
\epsilon}
(d_{2}d_{3})  ^{\frac{1}{2} + \epsilon}
(d_{3}d_{4})  ^{\frac{1}{2} + \epsilon}(d_{4}d_{5})  ^{\frac{1}{2} +
\epsilon}
d_{2} ^{\frac{1}{2} + \epsilon}l_{3} ^{\frac{1}{2}}l_{4} ^{- \frac{3}{2}}
d_{4} ^{\frac{1}{2} + \epsilon} \cdot S,
\endalign
$$
where $d_{1} = l_{2},$ $d_{j} = 2 ^{\a _{j}} b_{j},$ $\a _{j} = 0$ or $1,$
and $b_{j}|\frac{M}{2},$ $(b_{j},\; b_{j + 1}) = 1$
$(j = 1,\; 2,\; \ldots,\; 5),$  and $S$ is a sum of absolute values of the multiple
Dirichlet series $Z_M$ at various arguments of the characters. We can take
$$\align S &= \sum_{a = 1,\; -1}\;
\sum_{l|M}Z_{M}\left(\frac12, 1 + \epsilon; \chi_{a l}\right)  \\
&=  \sum_{a = 1,\; -1}\;
\sum_{l|M} {\underset\; (d,\;
M)=1\to{\sum_{d =
d_0 d_{1} ^{2} }}}
\frac{|L\left(\frac12, F, \chi_{d_0} \chi_{a l}\right) P_{d_0,d_1}^{(a
l)}(1/2)|} {d ^ {1 + \epsilon}}.  \endalign$$

The positive integer $d_{4}$ is such that
$d_{4} = d_{2}l_{3}l_{4}$ modulo squares, and $l_{3},$ $l_{4}| M.$
Since $M$ is square free, it follows that
$$
\text{ord}_{p}\left(\frac{d_{2}l_{3}d_{4}}{l_{4} ^{3}} \right)\le 2,
$$
for any prime dividing $\frac{M}{2}.$ Consequently,
$$|t| ^{5 + 10\epsilon}M ^{10\epsilon}(d_{1}d_{2})  ^{\frac{1}{2} +
\epsilon}
(d_{2}d_{3})  ^{\frac{1}{2} + \epsilon}
(d_{3}d_{4})  ^{\frac{1}{2} + \epsilon}(d_{4}d_{5})  ^{\frac{1}{2} +
\epsilon}
d_{2} ^{\frac{1}{2} + \epsilon}l_{3} ^{\frac{1}{2}}l_{4} ^{- \frac{3}{2}}
d_{4} ^{\frac{1}{2} + \epsilon} \ll_{\epsilon}M ^{3 + 16\epsilon}|t| ^{5 +
10\epsilon}.$$
We finally arrive
at the bound

$$
Z_M\left(\frac12, -\epsilon + it; \chi_{a_2l_2}, \chi_{a_1l_1}\right) \ll_\epsilon
M ^{3 + 30 \epsilon}|t| ^{5 + 10\epsilon}\sum_{a = 1,\; -1}\;
\sum_{l|M} {\underset\; (d,\;
M)=1\to{\sum_{d =
d_0 d_{1} ^{2} }}}
\frac{|L\left(\frac12, F, \chi_{d_0} \chi_{a l}\right) P_{d_0,d_1}^{(a
l)}(1/2)|} {d ^ {1 + \epsilon}}.
\tag 4.34
$$

We now need to establish that
$Z_M(\frac12,w;\chi_{a_2l_2},\chi_{a_1l_1})$ is analytic for $w$
in the region described in the proposition. We have already shown,
in Proposition 4.10 that
the product
$$
A(s,w)A(\a(s,w))A(\b(s,w))
A(\b\a(s,w))B(s,w)B(\a(s,w))\Pcal(s,w)
\overarrow\Z_M(s,w;\chi_{\text{Div}(M)},\chi_{a_1l_1}).
$$
is an entire function of $s,w$. Specializing to $s = \frac{1}{2},$ we see
that the only possible poles of $\overarrow\Z_M(\frac12, w;
\chi_{\text{Div}(M)}, \chi_{a_1l_1})$ could occur at zeros of
$$A(1/2,w)A(\a(1/2,w))A(\b(1/2,w))
A(\b\a(1/2,w))B(1/2,w)B(\a(1/2,w))\Pcal(1/2,w).$$
Zeros of $\Pcal(1/2,w)$ can only occur on the real line, at $w = 0,$
$\frac34,$
$1.$ The other terms in the product have factors of the form $(1 - p ^{-2
+ 2w})$ for $p|M.$ Thus the only potential locations for poles in the
region under
consideration are $w = 1 + it,$ for a discrete sequence of $t\ne 0.$ Such
poles cannot occur, however, for the following reason.

For any $s,$ $w$ with $\Re (s) \ge \frac{1}{2}$ and $\Re (w) >1,$
$\overarrow\Z_M(s, w; \chi_{\text{Div}(M)}, \chi_{a_1l_1})$ is an analytic
function of $s$ and $w.$ Suppose $\overarrow\Z_M(\frac12, w;
\chi_{\text{Div}(M)},\chi_{a_1l_1})$ has a pole of order $\g>0$ at $w = 1
+ it_0.$
Then
$$
\lim_{(s, w)\rightarrow (\frac12,
1 +
it_0)}\Pcal_0(s,w)\overarrow\Z_M(s,w;\chi_{\text{Div}(M)},\chi_{a_1l_1})
\ne 0,
$$
where $\Pcal_0(s, w)$ is a product of $\g$ linear factors of the form
$w-1-it_0,$ $s + w - 3/2 - it_0,$ $2s + w -2 - it_0$ or $3s + w - 5/2 -
it_0.$ These correspond to potential zeros of the products
$A(\b(s, w)),$ $A(\b\a(s, w)),$ $B(s, w)$ and $B(\a(s, w)).$
By the analyticity in $s,$ $w,$ we can interchange the limits:
$$
\lim_{w \rightarrow 1 + it_0} \Pcal_0(s, w) \lim_{s \rightarrow \frac12}
\overarrow\Z_M(s, w; \chi_{\text{Div}(M)}, \chi_{a_1l_1})
= \lim_{s \rightarrow \frac12} \lim_{w \rightarrow 1 + it_0} \Pcal_0(s, w)
\overarrow\Z_M(s, w; \chi_{\text{Div}(M)}, \chi_{a_1l_1}).
$$
On the right hand side, for any $s$ with $\Re (s) \ge \frac12,$ let
$$
T(s) = \lim_{w \rightarrow
1 + it_0} \Pcal_0(s, w) \overarrow\Z_M(s, w; \chi_{\text{Div}(M)},
\chi_{a_1l_1}).
$$
Then $T(s)$ is an analytic function  around $s = \frac12.$ Since for
$\Re (s)$ sufficiently large the right hand side of $(4.10)$ converges
absolutely, it is clear that if $\Pcal_0(s, w)$ contains a factor of the
form $w - 1 - it_0$ then $T(s) = 0$ for all such $s.$ This would imply
that the left hand side equals zero, which
contradicts our assumption. In a similar way we will eliminate the
possibility of the other three factors dividing $\Pcal_0(s, w).$

By applying $(4.30)$ to $\b(s, w)$ and setting $w = 3/2 + it_0 -s,$
we obtain the relation
$$
\prod_{p|M}\big(1 - p ^{-2(3/2 + it_0 -s)}\big)
\overarrow\Z_M(1 + it_0, s - 1/2 - it_0) = \Phi_M (s - 1/2 - it_0)
\overarrow\Z_M(s, 3/2 + it_0 - s).
$$
For $\Re (s)$ sufficiently large and $t_0 \ne 0,$ the left hand side of
the
above converges absolutely, and hence the right hand side is an analytic
function of $s.$ Consequently, $\Pcal_0(s, 3/2 + it_0 -s)$ times the right
hand side will vanish identically if $\Pcal_0(s, w)$ contains a factor of
$s + w - 3/2 - it_0.$ As $\Phi (s - 1/2 - it_0)$ does not vanish
identically, it follows
that the right hand side of $(4.36)$ equals zero if we approach along the
line $w = 3/2 + it_0 - s.$ This is a contradiction, so $\Pcal_0(s, w)$
does not
contain a factor of $s + w - 3/2 - it_0.$

Similarly, applying $(4.29),$ $(4.30)$ and setting $w = 2 + it_0 -2s,$
we obtain the relation
$$
\multline
\prod_{p|M}\big(1 - p ^{-4 + 4s - 2it_0 )}\big) \prod_{p|M}\big(1 - p ^{-4
+ 2s - 4it_0 )}\big)
\prod_{p|M}\big(1 - p ^{-3 + 2s - 2it_0 )}\big) ^3
\overarrow\Z_M(1 + it_0,s - 1 - 2it_0)\\ =\Phi ^{(a_2l_2)} (2s - 1 - it_0)
\Psi_M (2s - 1 - it_0)\Phi_M (2s - 1 - it_0)
\overline\Z_M(s, 2 + it_0 - 2s).
\endmultline
$$
By the same argument as above, $\Pcal_0(s, w)$ does not contain a factor
of
$2s + w -2 - it_0.$

Finally, applying $(4.29)$ to $\a(s, w)$ and setting $w = 5/2 + it_0 -
3s,$ we obtain the relation
$$
\prod_{p|M}\big(1 - p ^{-2s}\big) ^3
\overarrow\Z_M(1 - s, 1 + it_0) = \Psi_M (1 - s)
\overarrow\Z_M(s, 5/2 + it_0 - 3s),
$$
from which it follows that $\Pcal_0(s, w)$ does not contain a factor of
$3s + w - 5/2 - it_0.$

The possibility of a pole at $w = 0$ can be eliminated in the same way.

To see that there may, actually, be a pole at $w = \frac34,$ observe that the
transformation $\a\b$ relates the hyperplane $w = 1$ to $3s + 2w -3 = 0.$
Since $w = 1$ may certainly be a pole, it follows from $(4.18)$ and $(4.23)$
that $3s + 2w - 3 = 0$ is a pole.

This establishes the analyticity of $Z(\frac12, w)$ for $-\epsilon < \Re (w) <
1 + \epsilon,$ except possibly at $w = \frac34, 1.$

The upper bound follows from $(4.11),$ $(4.34)$ and the
Phragmen--Lindel\"of
principle.

This completes the proof of Proposition 4.12.

\vskip 20pt

\noindent{\bf \S 4.5 The sieving process}
\vskip 10pt
In this section we will use the series $Z_M$ as building blocks to
construct
$$
Z(s, w) = \sum_d \frac{L(s, \chi_{d_0})^3}{|d|^w},
\tag 4.35
$$
where the sum ranges over square free integers $d_0$ and for each $d_0,$
$d$ is the associated fundamental discriminant. This is simply the series
$(4.12),$ as $\chi_{d_0} = \chi_d.$ The series $Z(s, w)$ will then inherit
its analytic properties from those of $Z_M.$

Our object is to prove

\proclaim{Theorem 4.13} Let the series $Z(s, w)$ be as defined above, and
choose any $\epsilon >0.$
When the specialization $s = \frac{1}{2}$ is made,
$Z\big(\frac{1}{2}, w \big)$ is an analytic function of $w$ for $\Re(w) >
\frac{4}{5}$ except for a pole of order 7 at $w=1.$ For $w = \nu + it,$
with
$\nu > \frac{4}{5},$ $Z\big(\frac{1}{2}, w \big)$ satisfies the upper
bound
$$
Z\Big(\frac{1}{2}, w \Big) \ll_\epsilon \cases \; 1 & \text{if \; $1+
\epsilon < \nu$}, \\ (1 + |t|)^{5(1 - \nu) + v(\epsilon)} & \text{if \;
$\frac{4}{5} < \nu \le
1+\epsilon$}, \endcases
$$
where $v(\epsilon)$ is an explicitly computable function satisfying
$\lim_{\epsilon \rightarrow 0}v(\epsilon) = 0.$

Also,
$$
\lim_{w \rightarrow 1}(w-1) ^7 Z\Big(\frac{1}{2}, w \Big) = \frac{6a_3}
{4 \pi ^2},
$$
where $a_3$ is given by $(3.3).$
\endproclaim
\vskip 10pt

In this section let $r$ denote a positive square free integer with
$(r, 2)=1.$ We also fix the notation $a_1,$ $a_2 \in \{1, -1\}$ and $l_1,$
$l_2 \in \{1,2\}.$ Let $F,$ as before, be the $GL(3)$ Eisenstein series
associated to $L(s, \chi_{d_0}) ^3,$ so
$L(s, F, \chi_{d_0}) = L(s, \chi_{d_0})^3.$ For any $l|r,$ define
$$
Z_{a_1l_1, a_2l_2} ^{(l)}(s, w) =
\underset d=d_0d_1^2\to{\sum_{(d_0, 2)=1,\;(d_1, 2l)=1}} \frac{L_2(s, F,
\chi_{d_0}\chi_{a_1l_1})\chi_{a_2l_2}(d_0)P_{d_0,d_1}^{(a_1l_1)}(s)}{d ^w}
\tag 4.36
$$
and as usual $d_0$ varies over positive square free integers and $d_1$
varies over positive integers.

If we then define
$$
Z_{a_1l_1, a_2l_2}(s, w; r)= \sum_{l|r} \mu (l) Z_{a_1l_1, a_2l_2}
^{(l)}(s, w),
\tag 4.37
$$
where $\mu$ denotes the usual M\"obius function, it is easy to check that
$$
Z_{a_1l_1, a_2l_2}(s, w; r) =
\underset d=d_0d_1^2\to{\sum_{(d_0d_1, 2)=1,\;d_1 \equiv 0 \mod r}} \frac{L_2(s,
F,
\chi_{d_0}\chi_{a_1l_1})\chi_{a_2l_2}(d_0)P_{d_0, d_1}^{(a_1l_1)}(s)}{d
^w}.
\tag 4.38
$$
In the next proposition we demonstrate that
$Z_{a_1l_1, a_2l_2}^{(l)}(s, w),$ and hence
$Z_{a_1l_1, a_2l_2}(s, w; r)$ can be written as a linear combination of
the functions
$Z_M(s, w; \chi_{a_2l_2}, \chi_{a_1l_1})$ whose analytic properties have
already been studied in the preceding sections.
\proclaim{Proposition 4.14} We have
$$
\multline
Z_{a_1l_1, a_2l_2} ^{(l)}(s, w)\cdot\prod_{p|l}(1 - p ^{-2s})^3 = \frac
{1}{2}\sum_{l_3|l}l_3 ^{-w}
\prod_{p|l_3}(1 - p ^{-2s})^3\cdot\sum_{m_1, m_2,
m_3|(l/l_3)}\frac{\chi_{a_1l_1
l_3}(m_1m_2m_3)
\chi_{a_2l_2}(l_3)}{(m_1m_2m_3) ^s}\\
\times \big( Z_{2l}(s, w; \chi_{a_2l_2}\chi_{m_1m_2m_3}, \chi_{a_1l_1l_3})
+ Z_{2l}(s, w; \chi_{a_2l_2}\chi_{-m_1m_2m_3}, \chi_{a_1l_1l_3}) \\
+ \chi_{-1}(m_1m_2m_3)Z_{2l}(s, w; \chi_{a_2l_2}\chi_{m_1m_2m_3},
\chi_{a_1
l_1l_3})
- \chi_{-1}(m_1m_2m_3)Z_{2l}(s, w; \chi_{a_2l_2}\chi_{-m_1m_2m_3},\chi_{a
_1l_1l_3})\big).
\endmultline
$$
\endproclaim
\vskip 10pt
{\bf Proof:} Referring to $(4.36)$ and $(4.9),$ write
$$
Z_{a_1l_1, a_2l_2} ^{(l)}(s, w) = \sum_{l_3|l}\sum_{(d_0d_1, 2l) = 1}
\frac{L_2(s, F, \chi_{d_0l_3}\chi_{a_1l_1})\chi_{a_2l_2}(d_0l_3)P_{d_0,
d_1}
^{(a_1l_1l_3)}(s)}{d_0 ^w l_3 ^w d_1 ^{2w}}.
$$
Replacing $L_2(s, F, \chi_{d_0l_3}\chi_{a_1l_1})$ by
$L_{2l}(s, F, \chi_{d_0l_3}\chi_{a_1l_1})\cdot\prod_{p|l}
(1 - \chi_{d_0l_3}\chi_{a_1l_1}(p)p ^{-s})^{-3}$ and multiplying both
sides by $\prod_{p|l}(1 - p ^{-2s})^3,$ the result follows after some
simple
manipulations, and the use of $\chi_{-1}$ to distinguish the cases
$m_1 m_2 m_3 \equiv 1 \mod 4$ and $m_1 m_2 m_3 \equiv 3 \mod 4.$

This completes the proof of Proposition 4.14.

It follows from Propositions 4.12, 4.14, and the definition of
$Z_{a_1l_1, a_2l_2}(s, w; r)$ in $(4.37)$ that for $\epsilon > 0,$ if
$w = \nu + it,$ with $\nu > -\epsilon,$ then
$Z_{a_1l_1,a_2l_2}(1/2, w; r)$ is analytic except for possible poles
at $w = \frac34,$ $1,$ and satisfies the upper bound
$$
Z_{a_1l_1, a_2l_2}\Big(\frac{1}{2}, -\epsilon + it; r \Big)
\ll_\epsilon r ^{3 + v_3(\epsilon)}|t| ^{5 + v_4(\epsilon)}
\sum_{a = 1,\; -1}\;
\sum_{l| 2r}\;  \sum_{d_{0}}
\frac{\big|L\big(\frac{1}{2}, \chi_{d_0} \chi_{a l}\big) \big| ^{3}}
{d_{0} ^ {1 + \epsilon}},
$$
with $v_3(\epsilon),$ $v_4(\epsilon)$ some explicitly computable functions
satisfying $\lim_{\epsilon\rightarrow 0}v_3(\epsilon) =
\lim_{\epsilon\rightarrow 0}
v_4(\epsilon) = 0.$
For $\nu >1,$ the series $Z_{a_1l_1, a_2l_2}(1/2, w; r)$ converges
absolutely, by $(4.11)$ and $(4.33),$ and a factor of $r ^{2\nu}$ factors
out of
the denominator. Thus $Z_{a_1l_1, a_2l_2}(1/2, 1 + \epsilon + it; r)
\ll_\epsilon r ^{-2 - 2\epsilon}.$ Combining these bounds and applying
Phragmen--Lindel\"of, we obtain, for $-\epsilon < \nu < 1 + \epsilon$ and
$|t| > 1,$
$$
Z_{a_1l_1, a_2l_2}\Big(\frac{1}{2}, \nu + it; r \Big) \ll_\epsilon
r ^{3 - 5 \nu + v_3(\epsilon)}
|t|^{5 - 5 \nu  +  v_4(\epsilon)}\sum_{a = 1,\; -1}\;\sum_{l| 2r}\;
\sum_{d_{0}}
\frac{\big|L\big(\frac{1}{2}, \chi_{d_0} \chi_{a l}\big) \big| ^{3}}
{d_{0} ^ {1 + \epsilon}}.
\tag 4.39
$$
We now define
$$
Z_{a_1l_1, a_2l_2}(s, w)= \sum_{(r, 2) = 1}\mu(r) Z_{a_1l_1, a_2l_2}(s, w;
r),
$$
and observe that
$$
Z_{a_1l_1, a_2l_2}(s, w)=
\sum_{(d_0, 2)=1}\frac{L_2(s, \chi_{d_0}\chi_{a_1l_1}) ^3
\chi_{a_2l_2}(d_0)}{d_0 ^w},
$$
where the sum is over odd, square free positive integers $d_0.$ The sum
over
$r$ has removed all $d_1 \ne 1$ from the sum. Applying the bound of
$(4.39)$
and taking $\nu > \nu_{0} > \frac{4}{5},$ we have
$$
\multline
Z_{a_1l_1, a_2l_2}\Big(\frac{1}{2}, \nu + it \Big)\ll_\epsilon
|t| ^{5 - 5 \nu  +  v_4(\epsilon)} \sum_{(r, 2) = 1}
(2r) ^{3 - 5 \nu + v_3(\epsilon)}\sum_{a = 1,\; -1}\;\sum_{l| 2r}\;
\sum_{d_{0}}
\frac{\big|L\big(\frac{1}{2}, \chi_{d_0} \chi_{a l}\big) \big| ^{3}}
{d_{0} ^ {1 + \epsilon}}\\
\ll_\epsilon |t| ^{5 - 5 \nu  +  v_4(\epsilon)}
\sum_{a = 1,\; -1}\;\sum_{l}\;
\sum_{d_{0}}
\frac{\big|L\big(\frac{1}{2}, \chi_{d_0} \chi_{a l}\big) \big| ^{3}}
{d_{0} ^ {1 + \epsilon} l ^{5 \nu - 3 - v_3(\epsilon)} } \; \sum_{r' \ge
1}
\frac{1}{{r'} ^{5 \nu -3 - v_3(\epsilon)}}\ll_{\nu_{0}, \epsilon}
|t| ^{5 - 5 \nu  +  v_4(\epsilon)},
\endmultline
\tag 4.40
$$
if $\epsilon$ is chosen sufficiently small. In $(4.40),$ the last estimate
follows from $(4.33).$

We have thus proved

\proclaim{Proposition 4.15}
For any $a_1, a_2 \in \{ 1, -1\}$ and $l_1, l_2 \in \{1, 2\},$ the series $Z_{a_1l_1,
a_2l_2}\big (\frac{1}{2}, w \big)$ is analytic for $w = \nu + it$ when
$\nu > \frac{4}{5},$ except possibly for a pole at $w=1.$ For $|t|>1$
it satisfies the upper bound
$$Z_{a_1l_1, a_2l_2}\left(\frac{1}{2}, \nu + it\right)\ll_\epsilon
|t| ^{5 - 5 \nu  +  v_4(\epsilon)}.$$
\endproclaim
\vskip 10pt

To complete the proof of the first part of Theorem 4.13, we make choices of
$1,-1,2,-2$ for $a_1l_1$ and $a_2l_2$ and take linear combinations of
$Z_{a_1l_1,a_2l_2}(1/2,w)$ to isolate sums over $d_0>0,d_0<0$, and for
each sign, sums over
$d_0 \equiv 1 \mod 8$, $d_0
\equiv 5 \mod 8$, $d_0 \equiv 3 \mod 4$ and $d_0 \equiv 1 \mod 4$.
After these sums are isolated, the $2$-factor of the $L$-series can
be restored, and
the analyticity of $Z(\frac12,w)$ for $w \ne 1$ together with the upper
bound stated in Theorem
4.13 follows.

It now remains to calculate the order of the pole and compute the
leading coefficient in
the Laurent expansion at $w=1$.  This can be done directly from the
analytic information
and functional equations we have accumulated about
$Z_{a_1l_1,a_2l_2}(s,w)$.  However, it is
an intricate computation, and so we will instead make use of the
computations already
performed in Section 3 for a general multiple Dirichlet series.

In the notation of Section 3, taking $m = 3$, $Z(s,w) = Z(s,s,s,w)$, where
$Z(s_1,s_2,s_3,w)$ is defined by (3.6).  In the previous work of this
section we considered
the $L$-series $L(s,F) = \zeta (s)^3$.  Here $F$ was an Eisenstein
series on $GL(3)$
specialized to the center of the critical strip.  We could just have
easily have considered
the $L$--series associated to $F'$, a general minimal parabolic
Eisenstein series.  In the
case of $F$, the Euler product parameters at a prime $p$ were $\a_p =
\b_p = \g_p =1$ and
the corresponding local factor of the Euler product was $(1 -
p^{-s})^{-3}$.  For the more
general $F'$, we can take $\a_p = p^{-\epsilon_1}, \; \b_p = p^{-\epsilon_2}, \;
\g_p = p^{\epsilon_1+\epsilon_2}$.  The corresponding local factor of
$L(s,F')$ is then
equal to $\big((1-p^{-s - \epsilon_1})(1-p^{-s - \epsilon_2})(1-p^{-s
+\epsilon_1 +\epsilon_2})\big)^{-1}$.
Applying exactly the same arguments as before, we may obtain the
analytic continuation of
the more general object
$$
Z(s+\epsilon_1,s+\epsilon_2,s-\epsilon_1-\epsilon_2, w) = \sum_d
\frac{L(s+\epsilon_1, \chi_{d_0})  L(s+\epsilon_2,
\chi_{d_0})  L(s-\epsilon_1-\epsilon_2,
\chi_{d_0})   }{|d|^w}
$$
in a neighborhood of $s = 1/2$ and $\epsilon_1 = \epsilon_2 = 0$.
Setting $s_1 =
s+\epsilon_1$, $s_2=s+\epsilon_2$ and $s_3 = s-\epsilon_1-\epsilon_2$,
we are in a position
to take advantage of the calculations done in Section 3, as we have
established the
conjectured analytic continuation.  This completes the proof of Theorem
4.13.

It is worth remarking that we could just as easily have proved the
more general analytic
continuation of $Z(s_1,s_2,s_3,w)$.  However, our intent was to make
the outlines of the
technique as clear as possible.  Writing out the explicit details in
greater generality
would have made it significantly harder to distinguish the ideas
through the notation.

We now have only a small additional piece of work to do to complete
the proof of the first part of Theorem
1.1.  Applying the integral transform
$$
\frac{1}{2 \pi i}\int_{2 - i \infty}^{2 + i \infty}\frac{x^w
dw}{w(w+1)} = \cases (1-1/x)
&\;
\text{if $x>1$},\\
0\; & \text{if $0<x\le 1$,}\endcases
$$

we obtain first

$$
\frac{1}{2 \pi i}\int_{2 - i \infty}^{2 + i \infty}\frac{Z(1/2,w)x^w
dw}{w(w+1)} =
\sum_{|d|<x} L\left(\frac12, \chi_{d}\right) ^3 \left( 1 - \frac{|d|}{x} \right).
$$

Moving the line of integration to $\Re(w) = \frac{4}{5} + \epsilon$, for
$\epsilon >0$, we pick
up from the pole at $w=1$ a polynomial type expression of the form $x(A_6 (\log x)
^6 + A_5 (\log x) ^5 +
\dots + A_0),$ where the constants $A_6, ..., A_0$ are computable and
$$
A_6 = \frac{6a_3}{8 \pi^2 6!},
$$
i.e., $1/2$ the constant of Theorem 4.13, divided by $6!$. The integral
at $\Re (w) = \frac{4}{5} +
\epsilon$ converges absolutely by the upper bound estimate of Theorem
4.13, and contributes
an error on the order of $x^{\frac{4}{5} + \epsilon}.$ This completes the
proof of the
first
part of Theorem 1.1

\vskip 20pt

\noindent{\bf \S 4.6
An unweighted estimate}
\vskip 10pt
In this section we will prove the second part of Theorem 1.1.
An essential ingredient of an estimate for such a theorem, and, more generally, an
estimate for an unweighted sum $\sum_{d<x}
a_d$ when $a_d$ is not known to be non-negative, is  an estimate
for sums of $a_d$ over
short intervals.  In our case, if $d$ is square free then $a_d =
L(1/2,\chi_d)^3$, while
if $d = d_0d_1^2$ with
$d_0$ square free, then
$$
a_d = L(1/2,\chi_{d_0})^3P_{d_0,d_1}(1/2),
\tag 4.41
$$
where
$d^{-\epsilon }\ll P_{d_0,d_1}(1/2)\ll d^\epsilon$.   Here
$P_{d_0,d_1}(1/2)$ is a linear
combination of
$P_{d_0,d_1}^{(a_1l_1)}(1/2)$.  As a first step we will require the
following.
\proclaim{Proposition 4.16}
For $x>0$ sufficiently large,  $\epsilon >0$, and $\frac35 < \theta_0 \le 1$,
$$
\sum_{|d-x|< x^{\theta_0}} L(1/2,\chi_{d_0})^2 \ll_\epsilon
x^{\theta_0 + \epsilon}.
$$
The sum here is over $d$ of the form $d = d_0 m^2$ for some $m$, with
$d_0$ square free and either positive or negative.
\endproclaim
\vskip 10pt
{\bf Proof:}
 The easiest
way to prove the Proposition is to apply Theorem 4.1 of \cite{C--N} to
the analog of
$Z_M(s,w;\chi_{1},\chi_{1})$ of (4.13) in the case of $GL(2)$, i.e., when
$L_M(s,F,\chi_{d_0})= L_M(s,\chi_{d_0})^2$ for $d_0$ square free.
Then all coefficients
are non-negative. There are four gamma factors, so $A=2$ in their
notation, and the
result with exponent $3/5$ follows immediately, by ignoring all but
the square free terms.
(The sum over $m$ does not affect the exponent.)  The derivation of
the analytic
continuation and functional equation of $Z_M(s,w;\chi_{1},\chi_{1})$
is done precisely as
in the preceding sections and is omitted.  Alternatively, and more
traditionally, one could
obtain this analytic continuation by considering the Rankin-Selberg
convolution of a half--integral weight Eisenstein series with itself.  The analysis,
however, is considerably more
complicated.

Fix an $x$, and an $r<\sqrt{x}$.  The following Proposition will
begin the proof of our
estimate for unweighted sums of coefficients of $Z_{a_1l_1,a_2l_2}(s,
w;r)$.  To simplify
notation we will suppress $a_1,a_2,l_1,l_2$ and write
$$
a(d) = L_2(1/2,\chi_{d_0}\chi_{a_1l_1})^3\chi_{a_2l_2}(d_0)P_{d_0,d_1}^
{(a_1l_1)}(1/2).
\tag 4.42
$$
Thus
$$
Z_{a_1l_1,a_2l_2}(1/2, w; r)=
\underset d = d_0d_1^2 \to {\sum_{(d_0d_1,2)=1,d_1 \equiv 0 \mod r}}
\frac{a(d)}{d^w}.
\tag 4.43
$$
\proclaim{Proposition 4.17}
Fix $x,T >0$, $r$ square free, $a_1, a_2 \in \{1, -1\}, \; l_1, l_2\in \{1, 2\}, $ and
$\epsilon >0$.  Let
$$
I_1(r) = \frac{1}{2 \pi i}\int_{1+\epsilon - i T}^{1+\epsilon + i T}
\frac{Z_{a_1l_1,a_2l_2}(1/2, w;r)x^w dw}{w}.
$$
Then for any $1 \ge \theta_0 >3/5$
$$
I_1(r) = \sum_{d<x,d \equiv 0 \mod {r^2}} a(d) \; + \;
\O_\epsilon \big(x^\epsilon r^\epsilon
\big(\frac{x}{r^2}\big)^{(1+\theta_0)/2}\big)
\; + \;\O_\epsilon \big(x^\epsilon r^\epsilon
\frac{1}{T}\big(\frac{x}{r^2}\big)^{(3-\theta_0)/2}\big).
$$

\endproclaim
\vskip 10pt
{\bf Proof:}
Applying the integral transform
$$
\frac{1}{2 \pi i}\int_{1+\epsilon - i T}^{1+\epsilon + i T}\frac{x^w
dw}{w} = \cases
1&\;
\text{if $x > 1$},\\
0\; & \text{if $0 < x < 1$}\endcases
+ \O_\epsilon \left( x^{1+\epsilon}\min \left( 1,\;\; \frac{1}{T|\log
(x)|}\right) \right)
$$
to $Z_{a_1l_1,a_2l_2}(1/2, w;r)$ and interchanging the order of
summation and integration,
as we are in a region of absolute convergence, we obtain
$$
I_1(r) = \sum_{d<x,d \equiv 0 \mod {r^2}} a(d) + E_1,
$$
where
$$
E_1 \ll_\epsilon \sum_{d \equiv 0 \mod {r^2},d\ne 0}|a(d)|
\left(\frac{x}{d}\right)^{1 + \epsilon}
\min \left( 1, \;\; \frac{1}{T|\log (x/d)|}\right).
$$
Break the sum $E_1$ into three pieces: $E_1 = E_2 + E_3 + E_4$, where
the sums are over
$d < \frac12 x$, $d > 2x$ and $\frac12 x < d < 2x,$
respectively.
Write $ d = d_0 m^2 r^2$, with $d_0$ square free.  By its definition
in (4.42), together
with the bound of (4.11), we have the bound
$$
a(d) \; \ll_\epsilon  \; \big |L(1/2,\chi_{d_0}\chi_{a_1l_1}) \big |^3 \cdot d^\epsilon.
\tag 4.44
$$
Applying (4.44) to $E_2,E_3$, we see that
$E_2,E_3 \ll_\epsilon  x^{1+\epsilon}r^{-2 - 2\epsilon} T^{-1}$
follows immediately from the absolute convergence of $\sum
L(1/2,\chi_{d_0})^3
|d_0|^{-1-\epsilon}$  (which follows, as remarked before, from
Heath-Brown's results
\cite{H-B}).

To analyze $E_4$, note that we are summing over the range $\frac12 x
r^{-2} < d_0 m^2<2 x r^{-2}$, so
$$
E_4 \ll_\epsilon \sum_{d \equiv 0 \mod {r^2}, \; \frac12 x<d<2x}|a(d)|
   \cdot\min \left( 1, \;\;
\frac{1}{T|\log
(x/d)|}\right).
\tag 4.45
$$
We are summing over the range $\frac12 x r^{-2} < dr^{-2} =d_0 m^2<2 x
r^{-2}$.
Consequently, for any $\theta_0 >0$ we may write $ d_0 m^2 = [x r^{-2} +
d'(x
r^{-2})^{\theta_0}
+d'']$.  As $d',d''$ vary over the ranges $0 \le |d'| \ll (x
r^{-2})^{1-\theta_0}$ and
$0 \le d'' \ll (x r^{-2})^{\theta_0}$, the full range of values of
$d_0 m^2$ will be hit.
We will treat the cases $d' = 0, -1$ and $d' \ne 0, -1$ separately.

Write $E_4 = E_5 + E_6$ where $E_5$ is the sum over $d$ with $d' =
0, -1$.  Then choosing 1 in
the minimum of (4.45) we have
$$
E_5 \ll \sum_{d'=0,-1}\;\; \sum_{0 \le d'' \ll (x r^{-2})^{\theta_0}}
|a(d)| \;
 = \;  {\sum}^\ast     |a(d)|,
$$
where ${\sum}^\ast$ denotes the sum ranging over $d', d_0, m$ satisfying
$d' =
0,-1$ and
$$
0 \le |d_0m^2 - x
r^{-2} - d'(x r^{-2})^{\theta_0}| \ll (x r^{-2})^{\theta_0}.
$$
Also, by (4.44)
$$
a(d) \ll_\epsilon r^\epsilon x^\epsilon |L(1/2,\chi_{d_0}\chi_{a_1l_1})|^3.
$$
It follows by the Cauchy--Schwartz inequality that
$$
E_5 \ll_\epsilon r^\epsilon x^\epsilon
\left({\sum}^{\ast\ast}  \big |L(1/2,\chi_{d_0}\chi_{a_1l_1}) \big |^4
\right)^{1/2}\left( {\sum}^{\ast\ast} \big |L(1/2,\chi_{d_0}\chi_{a_1l_1})\big |^2
\right)^{1/2},
$$
where ${\sum}^{\ast\ast}$ denotes the sum ranging over $d_0, m$
satisfying the condition
$$\left |d_0 - \frac{x r^{-2}}{m^2}\right | \ll \frac{(x
r^{-2})^{\theta_0}}{m^2}.$$ Using \cite{H-B} to bound the sum of fourth
powers
by $x$, and using  Proposition 4.16 to
bound the sum over squares we obtain
$$
E_5 \; \ll_\epsilon \; r^\epsilon x^\epsilon \big(\frac{x}{r^2} \big)^{(1
+
\theta_0)/2} \;\;
\sum_{m=1}^\infty m^{-1-\theta_0} \; \ll_\epsilon \; r^\epsilon x^\epsilon
\big(\frac{x}{r^2} \big)^{(1 +
\theta_0)/2}.
\tag 4.46
$$
To bound $E_6$ we first use the same argument as above to bound the
sum over $d''$ for
fixed $d'$.   We then observe that for $d' \ne 0, -1$ and any $d''$ we
have
$|\log(d/x)|^{-1} \ll (x r^{-2})^{1-\theta_0}/|d'|$. Taking the
$\log$ term in the minimum
of (4.45)  and  summing over $d'\ne 0$ we obtain

$$
E_6 \ll_\epsilon r^\epsilon x^\epsilon \left(\frac{x}{r^2} \right)^{(1 +
\theta_0)/2}
T^{-1}\sum_{d' \ne 0, -1}\left (x r^{-2}\right )^{1-\theta_0}/|d'|
\ll_\epsilon
r^\epsilon x^\epsilon
T^{-1}\left(\frac{x}{r^2} \right)^{(3 - \theta_0)/2}.
\tag 4.47$$
This completes the proof of Proposition 4.17.

Continuing with the proof of the Theorem, we now define, for $\epsilon >
0$,
and any $ - \epsilon \le \sigma \le 1-\epsilon$
$$
I_2(r, \sigma) = \frac{1}{2 \pi i}\int_{\sigma - i T}^{\sigma + i T}
\frac{Z_{a_1l_1,a_2l_2}(1/2, w;r)x^w dw}{w}
\tag 4.48
$$
and
$$
I_3(r, \sigma) = \frac{1}{2 \pi i}\int_{\sigma + i T}^{1+ \epsilon  + i T}
\frac{Z_{a_1l_1,a_2l_2}(1/2, w;r)x^w dw}{w},
\quad I_4(r, \sigma) = \frac{1}{2 \pi i}\int_{1 + \epsilon - i T}^{\sigma - i T}
\frac{Z_{a_1l_1,a_2l_2}(1/2, w;r)x^w dw}{w}.
$$
Thus,
$$\multline
I_1(r) \; = \; x \sum_{i=0}^6 d_i (r) (\log x)^i \, + \, I_2(r, \sigma) +I_3(r,
\sigma)  +\\ + I_4(r, \sigma) \, + \,
\delta_\sigma\cdot {\scriptstyle \frac43} x^\frac34 \cdot\underset w =
\frac34\to{\text{Res}}\Big(  Z_{a_1l_1,a_2l_2}(1/2, w;r)\Big),
\endmultline \tag 4.49
$$
for some computable constants $d_i(r)$.   The main
term is contributed by the seventh  order pole at
$w=1$ and the  residue term comes from the  possible
pole at $w = \frac34$, provided $-\epsilon < \sigma < \frac34-\epsilon$  for some
sufficiently small
$\epsilon>0$. Here $\delta_\sigma = 1$ if $-\epsilon<\sigma<\frac34-\epsilon$ and
$\delta_\sigma=0,$ otherwise. Note that there is no pole at
$w = 0$, so there are no additional error terms.

 It immediately follows from Proposition $4.17$ and  $(4.49)$ that
$$\multline\underset d \;{\scriptscriptstyle\text{ squarefree} }\to{\sum_{d < x}}
a_d
\; =
\; \sum_{r \le \sqrt{x}} \mu(r)\Bigg[x\sum_{i=0}^6 d_i(r)(\log x)^i \; + \;  I_2(r,
\sigma) \; + \;   I_3(r, \sigma) \;  + \; I_4(r, \sigma) \\
 + \; \delta_\sigma\cdot {\scriptstyle \frac43} x^\frac34 \cdot\underset w =
\frac34\to{\text{Res}}\Big(  Z_{a_1l_1,a_2l_2}(1/2, w;r)\Big)  \; + \;
\O_\epsilon
\left(x^\epsilon r^\epsilon
\left(\frac{x}{r^2}\right)^{(1+\theta_0)/2}\right)
\; + \;  \O_\epsilon \left(x^\epsilon r^\epsilon
\frac{1}{T}\left(\frac{x}{r^2}\right)^{(3-\theta_0)/2}\right)
\Bigg]\endmultline\tag 4.50$$
The sum $\sum_{r \le \sqrt{x}} \mu(r) x\sum_{i=0}^6 d_i(r)(\log x)^i$ will give
the main term of the second part of Theorem $1.1$ with a negligible error of
$O(x^{\frac12 + \epsilon})$. Thus, to complete the proof of Theorem $1.1$ it
remains to estimate the integrals and error terms in $(4.50)$.  These will be estimated
by breaking the sum over $r$ into $1 \le r \le x^\gamma$ and $x^\gamma < r \le
\sqrt{x}$ for some $0 < \gamma \le \frac12$ to be chosen later. We note that we will make
different choices of $T$ and $\sigma$ depending on whether $1 \le r \le x^\gamma$ or
$x^\gamma < r \le \sqrt{x}$. After computing all the error terms, we will make an optimal
choice of the variables $\gamma, \sigma, T, \theta_0.$
\vskip 10pt
  In order to estimate the integrals in $(4.50),$ we make use of the upper bound
$(4.39)$.  It follows that for $-\epsilon \le \nu < 1$,
$$
Z_{a_1l_1, a_2l_2}\Big(\frac{1}{2}, \nu + it; r \Big) \ll_\epsilon
r ^{3 - 5 \nu + v_3(\epsilon)}
\big(1+|t|\big)^{5 - 5 \nu  +  v_4(\epsilon)}\sum_{a = 1,\; -1}\;\sum_{l| 2r}\;
\sum_{d_{0}}
\frac{\big|L\big(\frac{1}{2}, \chi_{d_0} \chi_{a l}\big) \big| ^{3}}
{d_{0} ^ {1 + \epsilon}}.
\tag 4.51
$$

\proclaim{Proposition 4.18}
Let $x,$ $T >0,$ $r$ square free, and $\epsilon >0$. The integral $I_2(r, -\epsilon)$
given in $(4.48)$ satisfies
$$
\multline
 I_2(r, -\epsilon) = \frac{1}{2 \pi i}\int_{-\epsilon - i T}^{-\epsilon + i T}
\frac{Z_{a_1l_1,a_2l_2}(1/2, w;r)x^w \;dw}{w}\\
\ll_\epsilon r ^{3 + v_5(\epsilon)}T ^{\frac{9}{2} + v_6(\epsilon)}
\sum_{a = 1,\; -1}\;
\sum_{l|2r}\; \sum_{d_{0}}
\frac{\big|L\big(\frac{1}{2}, \chi_{d_0} \chi_{a l}\big) \big| ^{3}}
{d_{0} ^ {1 + \epsilon}},
\endmultline
$$
where $v_5(\epsilon)$ and $v_6(\epsilon)$ are some explicitly computable
functions satisfying
$$
\lim_{\epsilon\rightarrow 0}v_5(\epsilon) = \lim_{\epsilon\rightarrow 0}
v_6(\epsilon) = 0.
$$
\endproclaim
\vskip 10pt
{\bf Proof:}
The ultimate effect of this proposition is to save a power of $T ^{1/2}$
in the estimate for $I_2.$ To accomplish this, our goal is to apply the
functional equation $(4.31)$ to $Z_{a_1l_1, a_2l_2}(1/2, -\epsilon + it;
r),$
reflecting it into a region where it converges absolutely. This functional
equation reflects $Z$ into a new series which is actually a linear
combination of convergent series. This combination is summed over divisors
of
$2r$ and also over ratios of gamma factors corresponding to $L$--series
with
both positive and negative conductors. The easiest way to deal with this
is
to use the following notation:

Let $\overarrow \beta = (\beta_1, \beta_2, \beta_3, \beta_4, \beta_5),$
where
each $\beta_i \in \{0, 1 \}.$ Let $\Delta_{\overarrow \beta}$ denote the
product of gamma factors
$$
\Delta_{\overarrow \beta}(w) = G(w + \beta_1)G(w + \beta_2) ^3 G(2w - 1/2
+ \beta_3)G(w + \beta_4) ^3 G(w+\beta_5),
$$
where $G(w) = \pi ^{-w/2}\Gamma(w/2).$

Then for fixed $x$ and $T,$ it follows from $(4.31)$ and the explicit
forms
of the functional equations of Propositions 4.2 and 4.3 given by $(4.18)$
and
$(4.23)$ that we may reflect $Z_{a_1l_1, a_2l_2}(1/2, w; r)$ into a
complicated sum of Dirichlet series evaluated at $1-w.$ By a similar
argument
to the one given in the proof of Proposition 4.13, it can be observed that
the bound for the integral $I_2$ follows, if we show the estimate
$$
I_{\overarrow \beta}(y, T, \epsilon) := \int_{-T} ^{T}
\frac{\Delta_{\overarrow \beta}(1 + \epsilon + it)}
{\Delta_{\overarrow \beta}(-\epsilon - it)}\cdot
\frac{y ^{it}}{\epsilon + it}\;dt\ll_\epsilon T ^{\frac{9}{2} +
10\epsilon},
\tag 4.52
$$
where $y$ is any positive number.

To prove the estimate $(4.52),$ we first observe that from Stirling's
formula,
we have
$$
\frac{\Delta_{\overarrow \beta}(1 + \epsilon + it)}
{\Delta_{\overarrow \beta}(-\epsilon - it)}
= |t| ^{5 + 10 \epsilon + 10 it}e ^{cit}c'(\epsilon, \overarrow \beta)\left \{1 +
\O\left(\frac{1}{|t|}\right) \right\},
\tag 4.53
$$
for certain constants $c,$ $c'(\epsilon, \overarrow \beta).$

Replacing the ratio on the left hand side of $(4.53)$ with the main term,
the contribution from the error term is easily seen to be bounded above by
$\O \big( T ^{4 +\epsilon}\big),$ and using the expansion
$$
\frac{1}{\epsilon + it} = - \frac{i}{t}\left(1 +
\left(\frac{i\epsilon}{t}\right) + \left( \frac{i\epsilon}{t}\right) ^2 +
\cdots \right),
$$
it is enough to prove that
$$
\int_{1} ^{T}t ^ {u + 10 it}y ^{it}\;dt \ll
\cases T ^{u + \frac{1}{2}}
&\; \text{if $u\ge 0$},\\
T ^{\frac{1}{2}} \; & \text{if $u < 0.$}\endcases
\tag 4.54
$$
This is a simple consequence of the following lemma \cite{T}.
\proclaim{Lemma 4.19}
Let $F(x)$ be a real function, twice differentiable, and let
$F''(x) \ge m > 0,$ or $F''(x) \le -m < 0,$ for any $x,$ $a \le  x \le b.$
Let $G(x)/F'(x)$ be monotonic, and $|G(x)|\le M.$ Then
$$
\left|\int_{a} ^{b} G(x) e ^{iF(x)}\;dx\right|\le \frac{8M}{\sqrt m}.
$$
\endproclaim
\vskip 10pt

Choosing $F(t) = t(10\log t + \log y)$ and $G(t) = t ^u,$ we can divide
the interval $[1, T]$ in several subintervals such that the conditions in
the Lemma 4.19 are satisfied in each subinterval. The bound $(4.54)$
follows.

This completes the proof of Proposition 4.18.

\proclaim{Lemma 4.20} Let $0 < \gamma < \rho$ and $x \to
\infty$. Then for any $\epsilon > 0$,
$$
\sum_{x^\gamma \le r \le  x^\rho} r^u
\sum_{a = 1,\; -1}\;
\sum_{l|2r}\; \sum_{d_{0}}
\frac{\big|L\big(\frac{1}{2}, \chi_{d_0} \chi_{a
l}\big)
\big| ^{3}} {d_{0} ^ {1 + \epsilon}} \ll_\epsilon \cases
x^{\rho(u+1)+\epsilon} &\text{if $u >
-1$}\\x^{\gamma(u+1)+\epsilon} &\text{if $u <
-1$}.\endcases
$$
\endproclaim

{\bf Proof:} Let $S$ denote the quadruple sum given above.
 By interchanging sums and writing $2r = l\cdot r_1$,  we
easily see that

$$\align S \;  &=  \; \sum_{a = 1,\; -1}\; \sum_{l\le
2x^\rho}
\; \sum_{d_{0}}\; 2^{-u} \sum_{2r \equiv 0 (l)} (2r)^u \cdot
\frac{\big|L\big(\frac{1}{2}, \chi_{d_0} \chi_{a
l}\big)
\big| ^{3}} {d_{0} ^ {1 + \epsilon}} \\
&\ll \;
\sum_{a = 1,\; -1}\; \sum_{l\le 2x^\rho}
\; \sum_{d_{0}} \sum_{\frac{2}{l} x^\gamma\le r_1 \le
\frac{2}{l} x^\rho} l^{u+1+\epsilon} {r_1}^u
\cdot
\frac{\big|L\big(\frac{1}{2}, \chi_{d_0} \chi_{a
l}\big)
\big| ^{3}} {(l\cdot d_{0} )^ {1 + \epsilon}}.\endalign
$$
Now, if $u < -1,$ the inner sum over $r_1$ is a convergent
series which is bounded by $x^{\gamma(u+1)+\epsilon}.$
The remaining sums are absolutely convergent and bounded
by $(4.33)$. This establishes the first case of the Lemma.

If $u > -1,$ then the inner sum over $r_1$ is bounded by
$\left(\frac{2}{l} x^\rho\right)^{u+1+\epsilon}.$ The result
then again immediately follows from $(4.33).$ This
completes the proof of Lemma $4.20$.

We now proceed to systematically estimate the integrals and error terms
in $(4.50).$ Consider first the case $r > x^\gamma$ for some $\gamma$ to be determined
later.
Choosing
$T= x ^{(3 - \theta_0)/2}$, $
\sigma = 1 - \epsilon$, and summing over
$x ^\gamma < \gamma \le x ^{\frac{1}{2}},$ we find that the error
contributions
$$
\O_\epsilon \left(x ^\epsilon r ^\epsilon
\left(\frac{x}{r ^2}\right) ^{\frac{1 + \theta_0}{2}}\right),\;\quad\quad
\O_\epsilon \left(x ^\epsilon r ^\epsilon
\frac{1}{T}\left(\frac{x}{r ^2}\right)^{\frac{3-\theta_0}{2}}\right)
\tag 4.55
$$
are dominated by the first, which contributes (changing $\epsilon$
as appropriate)
$$
\sum_{x ^\gamma \le r \le x ^\frac{1}{2}} x ^\epsilon r ^\epsilon
\left(\frac{x}{r ^2}\right) ^{\frac{1 + \theta_0}{2}}
\ll_\epsilon x ^{\frac{1+ \theta_0}{2} - \gamma  \theta_0 + \epsilon}.
\tag 4.56
$$
Applying $4.51$ and Lemma $4.20$ to the definition of
$I_2(r, \sigma)$ given in $(4.48)$,  it follows that
$$
\sum_{x ^\gamma \le r \le x ^{\frac{1}{2}}}\left|I_2(r, 1 - \epsilon)\right|\ll_\epsilon
x ^{1 - \gamma + \epsilon},\tag 4.57
$$
again changing $\epsilon$ as appropriate.
Similarly, using $(4.51)$ and Lemma $4.20$, the integrals $I_3(r, 1-\epsilon)$ and $I_4(r,
1-\epsilon)$ contribute a smaller amount than the above error terms.

  Finally, we consider the case when $r < x^\gamma$. For this case, we choose $\sigma =
-\epsilon, \; T = \frac{x^\alpha}{r^\beta}$ with $\alpha - \beta\gamma > 0$ where $0
<
\alpha, \beta$ will be chosen later. First, we consider the error from the pole at $w =
\frac34.$ It follows from $(4.51)$ and Lemma $4.20$ that the
contribution is bounded by
$$\sum_{r < x^\gamma} r^{-\frac34 + \epsilon} x^\frac34 \ll
x^{\frac{\gamma}{4} + \frac34  +
\epsilon}.\tag 4.58$$
This error will be negligible compare to the others and can be discarded. The error coming
from the $I_2$ integral can be estimated using Proposition $4.18$ and Lemma $4.20.$ We
obtain
$$\align \sum_{r < x^\gamma} I_2(r, -\epsilon) \; &\ll \; x^{\frac92\alpha+\epsilon}
\sum_{r<x^\gamma} r^{3 - \frac92\beta+\epsilon} \sum_{a = 1,\; -1}\;\sum_{l| 2r}\;
\sum_{d_{0}}
\frac{\big|L\big(\frac{1}{2}, \chi_{d_0} \chi_{a l}\big) \big| ^{3}}
{d_{0} ^ {1 + \epsilon}}\tag 4.59\\
&\ll x^{\frac92\alpha} \cases x^{\gamma(4-\frac92\beta+\epsilon)} &\text{if $\beta <
\frac89$}\\
x^\epsilon &\text{if $\beta > \frac89$.}\endcases
\endalign$$
We now estimate the errors contributed by $(4.55).$ First
$$\sum_{r < x^\gamma} x^\epsilon
r^\epsilon \left(\frac{x}{r^2}\right)^{\frac{1+\theta_0}{2}}
\;
\ll
\; x^{\frac{1+\theta_0}{2} + \epsilon}. \tag 4.60$$
Secondly, we have
$$\align\sum_{r < x^\gamma} x^\epsilon
r^\epsilon \frac{1}{T} \left(\frac{x}{r^2}\right)^{\frac{3-\theta_0}{2}} \; &\ll
\; x^{-\alpha + \frac{3-\theta_0}{2} + \epsilon} \sum_{r < x^\gamma}
r^{\beta-3+\theta_0} \tag 4.61\\
& \ll \;
\cases x^{\frac{3-\theta_0}{2} - \alpha + \epsilon} & \text{if
$3-\theta_0-\beta>1$,}\\
x^{\frac{3-\theta_0}{2}-a+\gamma\theta_0-2\gamma+\epsilon}
&\text{if $3-\theta_0-\beta<1$,}\endcases\endalign$$
where $a = \alpha - \gamma\beta.$
All the other error terms contribute a smaller amount.  We leave
them as an exercise.

Collecting all the error terms in $(4.56)$, $(4.57)$, $(4.58)$,  $(4.59)$, $(4.60)$,
and
$(4.61)$, we see that if  $\beta > \frac89$ and $3 - \theta_0-\beta < 1,$ then
the total error is
 $$\O
\left(x^{1-\gamma+\epsilon} \; + \; x^{\frac{1+\theta_0}{2} + \epsilon} \; + \;
x^{\frac92\alpha + \epsilon} \; + \;
x^{\frac{3-\theta_0}{2}-a+\gamma\theta_0-2\gamma+\epsilon}\right).\tag 4.62$$
If we equalize these four error terms above, and solve in terms of $\theta_0$, it follows that
$$\gamma = \frac{1-\theta_0}{2}, \;\quad \alpha = \frac{1 + \theta_0}{9},
\;\quad a = 0 \; \implies \alpha = \gamma\beta.$$
The condition $3 - \theta_0 - \beta < 1$ implies that $\beta > 2 - \theta_0$ which implies
that $\alpha = \gamma\beta > \gamma(2 - \theta_0)$ which gives
$$\frac{1+\theta_0}{9} \; > \; \frac{1-\theta_0}{2}(2-\theta_0).$$
These inequalities imply that
$$\theta_0 > \frac{1}{18}(29 - \sqrt{265}).$$ With this choice, the total error in $(4.62)$
is $$\O\left(x^{\frac{1}{36}\big(47 - \sqrt{265}\big) + \epsilon}\right),$$
where $\frac{1}{36}(47 - \sqrt{265}) \sim 0.853366...$
This completes the proof of Theorem 1.1.


\Refs

\widestnumber\key{B--F--H--1}


\ref \key B--H \by E. Brezin and S. Hikami
\paper Characteristic polynomials of random matrices
\jour Comm. Math. Phys.
\vol 214
\yr 2000
\pages 111--135
\endref



\ref \key B--F--H--1 \by D. Bump, S. Friedberg and J. Hoffstein
\paper Sums of twisted $GL(3)$ automorphic $L$--functions
\jour Preprint
\endref



\ref \key B--F--H--2 \by D. Bump, S. Friedberg
and J. Hoffstein
\paper On some applications of automorphic forms to
number theory
\jour  Bull. A.M.S \vol 33\yr
1996 \pages 157--175
\endref


\ref \key C--F \by J. B. Conrey and D. W. Farmer
\paper Mean values of $L$--functions and symmetry
\jour Internat. Math. Res. Notices \vol 17
\yr 2000
\pages 883--908
\endref


\ref \key  CFKRS \by J.B. Conrey, D.W. Farmer,
J.P. Keating, M.O. Rubenstein, and N.C. Snaith
\paper Moments of zeta and L-functions
\jour
\yr
\pages
\endref



\ref \key C--G \by J. B. Conrey and S. M. Gonek
\paper High moments of the Riemann zeta-function
\jour Duke Math. J. \vol 107
\yr 2001
\pages 577--604
\endref


\ref \key C--Gh--1 \by J. B. Conrey and A. Ghosh
\paper A conjecture for the sixth power moment of the Riemann
zeta-function
\jour Internat. Math. Res. Notices \vol 15
\yr 1998
\pages 775--780
\endref

\ref \key C--Gh--2 \by J. B. Conrey and A. Ghosh
\paper On mean values of the zeta--function
\jour Internat. Mathematika
\vol 31
\yr 1984
\pages 159--161
\endref

\ref \key C--N \by K. Chandrasekharan and R. Narasimhan
\paper Functional equations with multiple gamma factors and the average
order
of arithmetical functions
\jour Annals of Math. \vol 76
\yr 1962
\pages 93--136
\endref


\ref \key D \by H. Davenport
\book Multiplicative number theory. Third edition. Graduate Texts in
Mathematics
\vol 74
\publ Springer-Verlag, New York
\yr 2000
\endref


\ref \key D--I \by J. M. Deshouillers and H. Iwaniec
\paper Kloosterman sums and Fourier coefficients of cusp forms
\jour Inventiones Math. \vol 70
\yr 1982/83
\pages 219--288
\endref


\ref \key F--F \by B. Fisher and S. Friedberg
\paper Double Dirichlet Series Over Function Fields
\jour
\vol
\yr
\pages
\endref



\ref \key F--H \by S. Friedberg and J. Hoffstein
\paper Nonvanishing theorems for automorphic
$L$--functions on $GL(2)$
\jour Annals of Math.
\vol 142
\yr 1995
\pages 385--423
\endref


\ref \key G \by A. Good
\paper The convolution method for Dirichlet series. The Selberg trace
formula
and related topics (Brunswick, Maine, 1984)
\pages 207--214
\inbook Contemp. Math., 53
\publ American Mathematical Society
\publaddr Providence, RI
\yr 1986
\endref



\ref \key G--H \by D. Goldfeld and J. Hoffstein
\paper Eisenstein series of $1/2$--integral weight
and the mean value of real Dirichlet $L$--series
\jour Inventiones Math. \vol 80 \yr 1985
\pages 185--208
\endref

\ref \key H \by J. Hoffstein
\paper
Eisenstein series and theta functions on the metaplectic group
\inbook  Theta functions: from the classical to the modern,
{\rm CRM Proc. Lect. Notes. 1}
\ed M. Ram Murty
\publaddr Providence, RI
\publ American Mathematical Society
\pages 65-104
\yr 1993
\endref

\ref \key H--B \by D.R. Heath-Brown
\paper A mean value estimate for real character sums
\jour Acta Arithmetica
\vol 72
\yr 1995
\pages 235--275
\endref


\ref \key H--L \by G. H. Hardy and J. E. Littlewood
\paper Contributions to the theory of the Riemann zeta-function and the
theory
of the distributions of primes
\jour Acta Mathematica
\vol 41
\yr 1918
\pages 119--196
\endref




\ref\key H\"o \by L. H\"ormander
\book An introduction to complex analysis in several
 variables
\publ Van Nostrand
\publaddr Princeton, N.J.
\yr 1966
\endref

\ref \key H--R \by J. Hoffstein and M. Rosen
\paper Average values of $L$--series in function fields
\jour J. Reine Angew. Math.
\vol 426
\yr 1992
\pages 117--150
\endref


\ref \key I \by A. E. Ingham
\paper Mean-value theorems in the theory of the Riemann zeta-function
\jour Proceedings of the London Mathematical Society
\vol 27
\yr 1926
\pages 273--300
\endref

\ref \key J \by M. Jutila
\paper On the mean value of $L({1\over 2},\,\chi )$ for real characters
\jour Analysis
\vol 1
\yr 1981
\pages 149--161
\endref



\ref \key K--S \by N. M. Katz and P. Sarnak
\book Random matrices, Frobenius eigenvalues, and monodromy
\publ American Mathematical Society Colloquium Publications
\vol 45
\publaddr American Mathematical Society, Providence, RI
\yr 1999
\endref



\ref \key Ke--Sn--1 \by J. P. Keating and N. C. Snaith
\paper Random matrix theory and $\zeta(1/2+it)$
\jour Comm. Math. Phys.
\vol 214
\yr 2000
\pages 57--89
\endref



\ref \key Ke--Sn--2 \by J. P. Keating and N. C. Snaith
\paper Random matrix theory and $L$--functions at $s=1/2$
\jour Comm. Math. Phys.
\vol 214
\yr 2000
\pages 91--110
\endref




\ref \key K \by T. Kubota
\book On automorphic forms and the reciprocity law in a number field
\publ Kinokuniya Book Store Co.
\publaddr Tokyo
\yr 1969
\endref



\ref \key Mot1 \by Y. Motohashi
\paper An explicit formula for the fourth power mean of the Riemann zeta-
function
\jour Acta Math
\vol 170
\yr 1993
\pages 181--220
\endref



\ref \key Mot2 \by Y. Motohashi
\paper A relation between the Riemann zeta-function and the
hyperbolic Laplacian
\jour Ann. Scuola Norm. Sup. Pisa Cl. Sci. (4)
\vol 22
\yr 1995
\pages 299--313
\endref



\ref \key S \by C. L. Siegel
\paper Die Funktionalgleichungen einiger Dirichletscher Reihen
\jour Math. Zeitschrift \vol 63 \yr 1956
\pages 363--373
\endref


\ref \key So \by K. Soundararajan
\paper Nonvanishing of quadratic Dirichlet $L$--functions at $s=\frac12$
\jour Annals of Math
\vol 152
\yr 2000
\pages 447--488
\endref


\ref \key St \by H. M. Stark
\
\ unpublished notes
\endref


\ref \key T \by E. C. Titchmarsh
\book The theory of the Riemann zeta-function. Second edition.
\ed D. R. Heath-Brown
\publ The Clarendon Press, Oxford University Press, New York
\yr 1986
\endref


\endRefs
\enddocument